\documentclass[preprint]{elsarticle}
\usepackage{hyperref}
\journal{Journal of Computational and Applied Mathematics}
\biboptions{sort&compress}
\usepackage[utf8]{inputenc}
\usepackage[top=2cm, bottom=2cm, left=2cm, right=2cm]{geometry}

\usepackage{graphicx}
\usepackage[export]{adjustbox}
\usepackage{amsmath}
\usepackage[version=4]{mhchem}
\usepackage{siunitx}
\usepackage{longtable,tabularx}
\usepackage{cleveref} 
\crefname{equation}{}{} 
\usepackage{setspace}
\usepackage{amssymb}
\usepackage{booktabs}
\usepackage{mathrsfs}
\usepackage{caption}
\DeclareCaptionLabelSeparator{none}{ } 
\usepackage{calc}
\usepackage{multicol}
\usepackage{pdfpages}
\usepackage{color}
\usepackage{listings}
\usepackage{lscape}
\usepackage{siunitx}
\usepackage{empheq}
\usepackage{mathtools,cases}
\usepackage{mdframed}
\usepackage{tikz}
\usepackage{framed}
\usepackage{amsfonts}
\usepackage{environ}
\usepackage{tabularx}
\usepackage{multirow}
\usepackage{enumitem}
\usepackage{caption}
\usepackage{subcaption}
\usepackage{makecell} 
\usepackage{placeins}
\usepackage{mymacros}

\usepackage{xcolor}

\definecolor{verifier}{rgb}{0,0,0}

\newcommand{\Xdiff}{\ensuremath W} 
\newcommand{\Xdae}{\ensuremath Z} 
\newcommand{\Xdiffrksub}{\ensuremath w} 
\newcommand{\Xdaerksub}{\ensuremath z} 
\newcommand{\Pe}{\ensuremath \mathrm{Pe}}

\newcommand{\initial}{{0}}
\newcommand{\final}{{f}}
\newcommand{\indexcell}{{\mathfrak{i}}}
\newcommand{\indexspecies}{{\mathfrak{k}}}

\newcommand{\speciessubscript}{{sp}}
\newcommand{\dx}{\Delta x}

\newcommand{\indexstage}{i}
\newcommand{\indexstagebis}{j}
\newcommand{\indexstep}{n}
\newcommand{\totalstep}{s}

\newcommand{\amplifactor}{b}
\newcommand{\rZN}{\mathrm{r}}
\newcommand{\kZN}{\mathrm{k}}

\begin{document}

\begin{frontmatter}

\title{High-order adaptive time discretisation of one-dimensional low-Mach reacting flows:
	a case study of solid propellant combustion}

\author[1,2]{Laurent François}
\cortext[mycorrespondingauthor]{Corresponding author}
\ead{laurent.francois@onera.fr}

\author[1]{Joël Dupays}
\author[1]{Dmitry Davidenko}
\author[2]{Marc Massot}

\address[1]{ONERA, DMPE, 6 Chemin de la Vauve aux Granges, 91120 Palaiseau, France}
\address[2]{CMAP, CNRS, École polytechnique, Institut Polytechnique de Paris, Route de Saclay, 91120 Palaiseau Cedex, France}

\begin{abstract}
Solving the reactive low-Mach Navier-Stokes equations with high-order adaptive methods in time is still a challenging problem, in particular due to the handling of the algebraic variables involved in the mass constraint. We focus on the one-dimensional configuration, where this challenge has long existed in the combustion community. We consider a model of solid propellant combustion, which possesses the characteristic difficulties encountered in the homogeneous or spray combustion cases, with the added complication of an active interface. The system obtained after semi-discretisation in space is shown to be differential-algebraic of index 1. A numerical strategy relying on stiffly accurate Runge-Kutta methods is introduced, with a specific discretisation of the algebraic constraints and time adaptation.
High order is shown to be reached on all variables, while handling the constraints properly. Three challenging test cases are investigated:  ignition, limit cycle, and unsteady response with detailed gas-phase kinetics. We show that the time integration method can greatly affect the ability to predict the dynamics of the system. The proposed numerical strategy exhibits high efficiency and accuracy for all cases compared to traditional schemes used in the combustion literature.
\end{abstract}

\end{frontmatter}

\section{Introduction}

Solving for steady and unsteady homogeneous and spray or solid-propellant combustion in one-dimensional flame simulations has attracted enormous attention in the combustion community starting with the seminal work conducted by M.D. Smooke and collaborators between Sandia and Yale University
\cite{docPREMIX,Smooke_Gio,darabiha1992,massot1998,massot2014}. 
To our knowledge, in these studies and subsequent papers, most of the one-dimensional unsteady CFD codes for such applications use a time integration based on splitting and/or implicit methods, which are usually limited to first-order accuracy in time. No mathematical analysis has been reported regarding the nature of the system of coupled equations obtained after semi-discretisation in space, where the handling of some variables (e.g. surface temperature for the solid propellant case, eigenvalues such as strain rate or mass fluxes) requires careful examination in connection with boundary conditions. Besides, relying on low-order integration methods may prove disadvantageous in terms of accuracy, performance and ability to resolve fine details of the dynamics. 
Indeed high-order methods are especially important  when investigating instabilities and nonlinear behaviours, e.g. limit cycles, where growth of some modes can only be captured by high-fidelity numerical methods \cite{cycle_limite_tube_rijke_high_order}.
One exception is the work on the dynamics of non-premixed counterflow flames by Im et al. \cite{petzold_dae_counterflow_2000}, where the high-order time integrator DASSL for differential-algebraic equations (DAEs) was used, but the constraint formulation was rather involved (introduction of compressibility effects to reduce the index) and details on the convergence and efficiency were not the main focus of the paper.

Even if solid propellant combustion brings in additional difficulties and constraints related to the heterogeneous nature of the flow, it involves the same problematic as homogeneous or two-phase flow combustion in one-dimensional low-Mach flows.
Therefore, we choose to focus on the more complete case of solid propellant combustion.
Before getting into the key original findings and contribution of the present paper, let us introduce the specific topic of solid propellant combustion and explain why the one-dimensional configuration is especially important in this case and requires an efficient and accurate numerical resolution.

Solid propellant combustion is a key element in rocket propulsion and has been extensively studied since the 1950s \cite{barrere, deluca1992,novozhilov1992, lengelle}.
It involves a solid phase and a gas phase, separated by an interface.
The solid is heated up by thermal conduction and radiation from the gas phase. At its
surface, the solid propellant decomposes, melts and evaporates through a pyrolysis process as the interface regresses. The resulting gaseous products react and form a flame which heats back the solid, allowing for a sustained combustion.
A key element is the regression speed of the propellant surface and its dependence on the
combustion chamber conditions. This has been extensively studied in a steady-state context through the use of analytical models \cite{DBW_original, BDPmodel, WSBmodel}.
These served as initial tools for the analysis of unsteady combustion dynamics, e.g. linearised frequency response \cite{reponse_culick_1968}, intrinsic combustion instabilities, and combustion chamber stability \cite{culickInstabilities}.

Unfortunately, these simple models are limited in accuracy, and they are not suited to the detailed study of transient dynamics. Since the 1990s, developments have therefore focused on CFD tools, both for steady and unsteady applications.
Recent one-dimensional CFD codes make use of complex chemical mechanisms \cite{meynet2006, theseAllumageSmith2011, Beckstead2007ModelingOC, meredith_HMX_ignition},
however such kinetics is still difficult to evaluate and validate, and the heterogeneity of widely used propellants, such as AP-HTPB, tends to limit the use of the one-dimensional approach to lower pressures \cite{theseAllumageSmith2011}, where three-dimensional effects (such as diffusion flames) can be neglected.
Detailed two- and three-dimensional CFD approaches taking into account the heterogeneity of the propellant have started to emerge \cite{jackson_2D,gallier3D, dmitry_compas}, however they are costly, even with simplified kinetics, and limited to the simulation of a very small burning area.

These models have been developed to simulate laboratory experiments, for example the burning of a small sample in an enclosed vessel.
They can however not be used for larger-scale simulations, such as the computation of the ignition transient in a solid rocket motor.
In such configurations, a three-dimensional CFD tool can hardly be used to solve the flow field inside the chamber and the surface combustion, due to extreme computational requirements imposed by the representation of the propellant flame. 
Indeed, this flame is typically only a few hundred micrometers thick, with strong gradients requiring mesh cells thinner than 1 \textmu m near the surface.
Using such a refined mesh along the complete length of combustion chamber (1-10 meters) would be prohibitive. Therefore, the only viable option is to use a coarse mesh near the surface, and solve all the surface combustion within a simpler submodel, e.g. a one-dimensional CFD tool. Such a coupling has been reported with very simplified propellant models \cite{li2015, johnston1995_couplage3D,bizot_ignition},
most often that of an initially inert solid propellant transitioning instantly to quasi-steady burning once its surface temperature becomes larger than a predefined ignition temperature.
More accurate one-dimensional models such as the ones discussed earlier may lead to a better resolution of the coupled dynamics, however the numerical strategy must be improved to ensure accuracy, efficiency and stability without requiring prohibitively small time steps.

Eventually, there is a convergent need for a detailed one-dimensional combustion model associated with an efficient numerical strategy, either for coupling with a three-dimensional CFD code,
or for detailed parametric studies of flame dynamics in a purely one-dimensional context.
This is specifically what we wish to investigate in this paper. It is instructive to conduct a short overview of the various numerical approaches presented in the literature for the time integration of one-dimensional solid propellant combustion models, which is typical of what is described in the general field of low-Mach one-dimensional combustion.
One of the earliest detailed one-dimensional model is presented by Erikson and Beckstead \cite{Erikson1998}. A splitting method is implemented to integrate the gas phase equation, using the ICE scheme \cite{ICEscheme}  to compute the pressure and velocity fields with an implicit scheme of first-order accuracy in time. The stiff chemical source terms are handled with DVODE \cite{DVODE}. The solid phase energy equation is integrated implicitly. The surface temperature and regression speed are then iterated upon until the interface conditions are met, each time performing the split integration of both phases.
Due to poor computational performance and large splitting errors, they transition in \cite{Erikson1999} to a fully implicit resolution of the gas phase, using the TWOPNT \cite{TWOPNT} algorithm to solve the system discretised in time with the first-order implicit Euler method.
Both phases still need to be iterated upon at each time step. 
The authors mention the attempt to use DASSL \cite{Petzold_DAE} instead, a high-order adaptive multistep method, however, following implementation difficulties, they reverted to the earlier first-order strategy.
Other researchers used a similar approach with an iterative coupling of both phases (see the review \cite{Beckstead2007ModelingOC}), sometimes with dual-time stepping \cite{Liau_RDX} or substepping \cite{meredith_HMX_ignition, theseAllumageSmith2011} to improve convergence of the solid phase.
To our knowledge, most one-dimensional unsteady solid propellant combustion codes use a similar iterative coupling procedure between the different phases, with a time integration that is overall first-order accurate in time.
Furthermore, the fixed-point approach for the coupling of both phases and their various physics only converges if the time step is sufficiently small, potentially hindering the stability and efficiency of the computation, in particular for highly refined meshes. It has indeed been reported \cite{theseAllumageSmith2011} that this iterative approach requires the convective CFL-number be smaller than 1 overall, even though implicit solvers are used to some degree.

The present contribution focuses on the one-dimensional solid propellant combustion in the low-Mach number approximation, and proposes a high-order time integration strategy based on existing singly-diagonally implicit Runge-Kutta methods (up to order 5), with a fully-coupled approach involving an original treatment of the mass conservation constraint, and ensuring unconditional stability of the overall computation.
The proposed approach can be generalized to any one-dimensional low-Mach combustion model for homogeneous or spray flames.

In the first part of the paper, we focus on the mathematical nature of the set of coupled equations obtained after semi-discretisation in space.
We show that it is a differential-algebraic system of index one. 
The knowledge of this particular property is decisive when looking for high-order time integration methods, which is the aim of the second part of this paper.
Many such methods are reported in the literature for this particular class of problems\footnote{It is important to stress the fact that, when the problem is considered in more than one dimension, the system may become differential-algebraic of index two, in particular due to the pressure field \cite{Petzold_DAE}, requiring a partial reformulation of the problem and the use of more advanced integration techniques \cite{nguessan2019}. However, the main results of the one-dimensional case should remain applicable.} \cite{hairer_book2, Petzold_DAE}. Based on a list of computational and mathematical requirements, we choose a family of singly-diagonally implicit Runge-Kutta methods with embedded lower-order solutions. This provide high-order accuracy, unconditional stability and native time adaptation capabilities.
We conduct the numerical analysis of the specific form of the algebraic constraints and introduce an original and efficient way of treating, in particular, the mass conservation constraint on the velocity field. This is a key issue for the efficiency and accuracy of the numerical method.
We then present a one-dimensional CFD code utilising this strategy to compute the time-dependent solution of the solid propellant combustion model. Spatial and temporal discretisations are verified with two test cases, for which quasi-analytical solutions are available: steady-state combustion and response to pressure oscillations.
In order to investigate the potential and efficiency of the method, we tackle three challenging test cases. 
First, we focus on the simulation of ignition transients with a simplified modelling, showing that the time adaptation capability of the proposed strategy is more efficient than traditional CFL- or variation-limited time step evaluation strategies commonly used in such applications.
Second, we study a solid propellant configuration whose steady-state solution is linearly unstable, leading to a initial instability that stabilises nonlinearly on a limit-cycle periodic solution. 
Simulation of this configuration shows that high-order time integration methods offer a major improvement over low-order methods, in terms of quality of the result, robustness, and ability to capture the limit cycle, as well as restitution time and ease of simulation setup.
Eventually, the proposed numerical strategy is tested on the unsteady combustion of an ammonium perchlorate (AP) monopropellant with detailed gas-phase kinetics.
It is shown that the convergence order and
efficiency of the method are not affected by the increased complexity and
stiffness of the detailed modelling level. Important gains in the quality of
the result, ease of simulation setup, and computational times are also observed.
Finally, we conclude on the proposed strategy and discuss possible extensions.

\section{Formulation of a semi-discrete unsteady solid propellant combustion model}
\label{section:formulation}

In this section, we first present a one-dimensional model of solid propellant combustion, very similar to the ones found in the one-dimensional low-Mach combustion literature \cite{docPREMIX,Smooke_Gio,darabiha1992,massot1998,massot2014}. 
The continuous equations are semi-discretised in space via a method of lines based on a finite-volume scheme.

\subsection{General modelling}
\label{section:baseEquations}
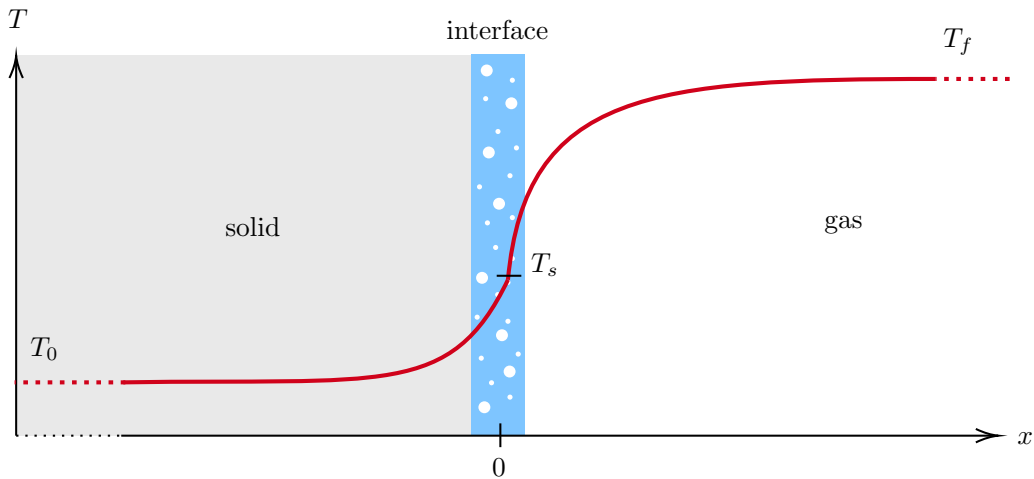
\begin{figure}[hbt!]
\centering

\resizebox{0.7\textwidth}{!}{
\centering

\tikzset{every picture/.style={line width=0.75pt}} 

\begin{tikzpicture}[x=0.75pt,y=0.75pt,yscale=-1,xscale=1]
	
	\draw  [color={rgb, 255:red, 0; green, 0; blue, 0 }  ,draw opacity=0 ][fill={rgb, 255:red, 126; green, 197; blue, 255 }  ,fill opacity=1 ] (349.13,100.2) -- (375.13,100.2) -- (375.13,220.6) -- (349.13,220.6) -- cycle ;
	\draw    (178.67,220.33) -- (604.13,220.33) ;
	\draw [shift={(606.13,220.33)}, rotate = 180] [color={rgb, 255:red, 0; green, 0; blue, 0 }  ][line width=0.75]    (10.93,-3.29) .. controls (6.95,-1.4) and (3.31,-0.3) .. (0,0) .. controls (3.31,0.3) and (6.95,1.4) .. (10.93,3.29)   ;
	\draw  [color={rgb, 255:red, 0; green, 0; blue, 0 }  ,draw opacity=0 ][fill={rgb, 255:red, 255; green, 255; blue, 255 }  ,fill opacity=1 ] (355.01,119.43) .. controls (355.43,117.88) and (357.03,116.97) .. (358.58,117.39) .. controls (360.12,117.81) and (361.04,119.41) .. (360.61,120.96) .. controls (360.19,122.5) and (358.59,123.42) .. (357.05,122.99) .. controls (355.5,122.57) and (354.58,120.97) .. (355.01,119.43) -- cycle ;
	\draw  [color={rgb, 255:red, 0; green, 0; blue, 0 }  ,draw opacity=0 ][fill={rgb, 255:red, 255; green, 255; blue, 255 }  ,fill opacity=1 ] (359.9,145.19) .. controls (359.9,143.59) and (361.21,142.29) .. (362.81,142.29) .. controls (364.42,142.29) and (365.72,143.59) .. (365.72,145.19) .. controls (365.72,146.8) and (364.42,148.1) .. (362.81,148.1) .. controls (361.21,148.1) and (359.9,146.8) .. (359.9,145.19) -- cycle ;
	\draw  [color={rgb, 255:red, 0; green, 0; blue, 0 }  ,draw opacity=0 ][fill={rgb, 255:red, 255; green, 255; blue, 255 }  ,fill opacity=1 ] (351.62,181.33) .. controls (351.62,179.73) and (352.92,178.43) .. (354.52,178.43) .. controls (356.13,178.43) and (357.43,179.73) .. (357.43,181.33) .. controls (357.43,182.94) and (356.13,184.24) .. (354.52,184.24) .. controls (352.92,184.24) and (351.62,182.94) .. (351.62,181.33) -- cycle ;
	\draw  [color={rgb, 255:red, 0; green, 0; blue, 0 }  ,draw opacity=0 ][fill={rgb, 255:red, 255; green, 255; blue, 255 }  ,fill opacity=1 ] (361.33,209.33) .. controls (361.33,207.73) and (362.63,206.43) .. (364.24,206.43) .. controls (365.84,206.43) and (367.14,207.73) .. (367.14,209.33) .. controls (367.14,210.94) and (365.84,212.24) .. (364.24,212.24) .. controls (362.63,212.24) and (361.33,210.94) .. (361.33,209.33) -- cycle ;
	\draw  [color={rgb, 255:red, 0; green, 0; blue, 0 }  ,draw opacity=0 ][fill={rgb, 255:red, 255; green, 255; blue, 255 }  ,fill opacity=1 ] (361.05,109.95) .. controls (361.05,109.27) and (361.6,108.71) .. (362.29,108.71) .. controls (362.97,108.71) and (363.53,109.27) .. (363.53,109.95) .. controls (363.53,110.64) and (362.97,111.2) .. (362.29,111.2) .. controls (361.6,111.2) and (361.05,110.64) .. (361.05,109.95) -- cycle ;
	\draw  [color={rgb, 255:red, 0; green, 0; blue, 0 }  ,draw opacity=0 ][fill={rgb, 255:red, 255; green, 255; blue, 255 }  ,fill opacity=1 ] (366.9,131.53) .. controls (366.9,130.84) and (367.46,130.29) .. (368.15,130.29) .. controls (368.83,130.29) and (369.39,130.84) .. (369.39,131.53) .. controls (369.39,132.21) and (368.83,132.77) .. (368.15,132.77) .. controls (367.46,132.77) and (366.9,132.21) .. (366.9,131.53) -- cycle ;
	\draw  [color={rgb, 255:red, 0; green, 0; blue, 0 }  ,draw opacity=0 ][fill={rgb, 255:red, 255; green, 255; blue, 255 }  ,fill opacity=1 ] (355.9,154.53) .. controls (355.9,153.84) and (356.46,153.29) .. (357.15,153.29) .. controls (357.83,153.29) and (358.39,153.84) .. (358.39,154.53) .. controls (358.39,155.21) and (357.83,155.77) .. (357.15,155.77) .. controls (356.46,155.77) and (355.9,155.21) .. (355.9,154.53) -- cycle ;
	\draw  [color={rgb, 255:red, 0; green, 0; blue, 0 }  ,draw opacity=0 ][fill={rgb, 255:red, 255; green, 255; blue, 255 }  ,fill opacity=1 ] (359.9,166.53) .. controls (359.9,165.84) and (360.46,165.29) .. (361.15,165.29) .. controls (361.83,165.29) and (362.39,165.84) .. (362.39,166.53) .. controls (362.39,167.21) and (361.83,167.77) .. (361.15,167.77) .. controls (360.46,167.77) and (359.9,167.21) .. (359.9,166.53) -- cycle ;
	\draw  [color={rgb, 255:red, 0; green, 0; blue, 0 }  ,draw opacity=0 ][fill={rgb, 255:red, 255; green, 255; blue, 255 }  ,fill opacity=1 ] (365.9,183.53) .. controls (365.9,182.84) and (366.46,182.29) .. (367.15,182.29) .. controls (367.83,182.29) and (368.39,182.84) .. (368.39,183.53) .. controls (368.39,184.21) and (367.83,184.77) .. (367.15,184.77) .. controls (366.46,184.77) and (365.9,184.21) .. (365.9,183.53) -- cycle ;
	\draw  [color={rgb, 255:red, 0; green, 0; blue, 0 }  ,draw opacity=0 ][fill={rgb, 255:red, 255; green, 255; blue, 255 }  ,fill opacity=1 ] (360.9,189.53) .. controls (360.9,188.84) and (361.46,188.29) .. (362.15,188.29) .. controls (362.83,188.29) and (363.39,188.84) .. (363.39,189.53) .. controls (363.39,190.21) and (362.83,190.77) .. (362.15,190.77) .. controls (361.46,190.77) and (360.9,190.21) .. (360.9,189.53) -- cycle ;
	\draw  [color={rgb, 255:red, 0; green, 0; blue, 0 }  ,draw opacity=0 ][fill={rgb, 255:red, 255; green, 255; blue, 255 }  ,fill opacity=1 ] (365.9,202.53) .. controls (365.9,201.84) and (366.46,201.29) .. (367.15,201.29) .. controls (367.83,201.29) and (368.39,201.84) .. (368.39,202.53) .. controls (368.39,203.21) and (367.83,203.77) .. (367.15,203.77) .. controls (366.46,203.77) and (365.9,203.21) .. (365.9,202.53) -- cycle ;
	\draw  [color={rgb, 255:red, 0; green, 0; blue, 0 }  ,draw opacity=0 ][fill={rgb, 255:red, 255; green, 255; blue, 255 }  ,fill opacity=1 ] (370.9,218.53) .. controls (370.9,217.84) and (371.46,217.29) .. (372.15,217.29) .. controls (372.83,217.29) and (373.39,217.84) .. (373.39,218.53) .. controls (373.39,219.21) and (372.83,219.77) .. (372.15,219.77) .. controls (371.46,219.77) and (370.9,219.21) .. (370.9,218.53) -- cycle ;
	\draw  [color={rgb, 255:red, 0; green, 0; blue, 0 }  ,draw opacity=0 ][fill={rgb, 255:red, 255; green, 255; blue, 255 }  ,fill opacity=1 ] (350.9,200.53) .. controls (350.9,199.84) and (351.46,199.29) .. (352.15,199.29) .. controls (352.83,199.29) and (353.39,199.84) .. (353.39,200.53) .. controls (353.39,201.21) and (352.83,201.77) .. (352.15,201.77) .. controls (351.46,201.77) and (350.9,201.21) .. (350.9,200.53) -- cycle ;
	\draw  [color={rgb, 255:red, 0; green, 0; blue, 0 }  ,draw opacity=0 ][fill={rgb, 255:red, 255; green, 255; blue, 255 }  ,fill opacity=1 ] (352.05,136.95) .. controls (352.05,136.27) and (352.6,135.71) .. (353.29,135.71) .. controls (353.97,135.71) and (354.53,136.27) .. (354.53,136.95) .. controls (354.53,137.64) and (353.97,138.2) .. (353.29,138.2) .. controls (352.6,138.2) and (352.05,137.64) .. (352.05,136.95) -- cycle ;
	\draw  [color={rgb, 255:red, 0; green, 0; blue, 0 }  ,draw opacity=0 ][fill={rgb, 255:red, 255; green, 255; blue, 255 }  ,fill opacity=1 ] (370.05,117.95) .. controls (370.05,117.27) and (370.6,116.71) .. (371.29,116.71) .. controls (371.97,116.71) and (372.53,117.27) .. (372.53,117.95) .. controls (372.53,118.64) and (371.97,119.2) .. (371.29,119.2) .. controls (370.6,119.2) and (370.05,118.64) .. (370.05,117.95) -- cycle ;
	\draw  [color={rgb, 255:red, 0; green, 0; blue, 0 }  ,draw opacity=0 ][fill={rgb, 255:red, 255; green, 255; blue, 255 }  ,fill opacity=1 ] (368.05,171.95) .. controls (368.05,171.27) and (368.6,170.71) .. (369.29,170.71) .. controls (369.97,170.71) and (370.53,171.27) .. (370.53,171.95) .. controls (370.53,172.64) and (369.97,173.2) .. (369.29,173.2) .. controls (368.6,173.2) and (368.05,172.64) .. (368.05,171.95) -- cycle ;
	\draw  [color={rgb, 255:red, 0; green, 0; blue, 0 }  ,draw opacity=0 ][fill={rgb, 255:red, 255; green, 255; blue, 255 }  ,fill opacity=1 ] (368.05,151.95) .. controls (368.05,151.27) and (368.6,150.71) .. (369.29,150.71) .. controls (369.97,150.71) and (370.53,151.27) .. (370.53,151.95) .. controls (370.53,152.64) and (369.97,153.2) .. (369.29,153.2) .. controls (368.6,153.2) and (368.05,152.64) .. (368.05,151.95) -- cycle ;
	\draw  [draw opacity=0][fill={rgb, 255:red, 0; green, 0; blue, 0 }  ,fill opacity=0.09 ] (127.45,100.2) -- (349.13,100.2) -- (349.13,220.6) -- (127.45,220.6) -- cycle ;
	\draw  [dash pattern={on 0.84pt off 2.51pt}]  (127.45,220.33) -- (178.67,220.33) ;
	\draw [color={rgb, 255:red, 208; green, 2; blue, 27 }  ,draw opacity=1 ][line width=1.5]    (185.5,208) .. controls (317.5,210) and (339.35,211.31) .. (364,182.6) ;
	\draw [color={rgb, 255:red, 208; green, 2; blue, 27 }  ,draw opacity=1 ][line width=1.5]    (364,182.6) .. controls (387.35,115.31) and (424.5,115) .. (563.5,115) ;
	\draw    (127.45,221) -- (127.45,99.8) ;
	\draw [shift={(127.45,97.8)}, rotate = 90] [color={rgb, 255:red, 0; green, 0; blue, 0 }  ][line width=0.75]    (10.93,-3.29) .. controls (6.95,-1.4) and (3.31,-0.3) .. (0,0) .. controls (3.31,0.3) and (6.95,1.4) .. (10.93,3.29)   ;
	\draw    (363.45,226.93) -- (363.45,215.33) ;
	\draw    (358.67,182.33) -- (370.67,182.33) ;
	\draw [color={rgb, 255:red, 208; green, 2; blue, 27 }  ,draw opacity=1 ][line width=1.5]  [dash pattern={on 1.69pt off 2.76pt}]  (599.5,115) -- (563.5,115) ;
	\draw [color={rgb, 255:red, 208; green, 2; blue, 27 }  ,draw opacity=1 ][line width=1.5]  [dash pattern={on 1.69pt off 2.76pt}]  (185.5,208) -- (127.5,208) ;
	
	\draw (605.2,206.2) node [anchor=north west][inner sep=0.75pt]   [align=left] {$\displaystyle x$};
	\draw (120,74.4) node [anchor=north west][inner sep=0.75pt]   [align=left] {$\displaystyle T$};
	\draw (557.6,89.8) node [anchor=north west][inner sep=0.75pt]   [align=left] {$\displaystyle T_{\final}$};
	\draw (131,181) node [anchor=north west][inner sep=0.75pt]   [align=left] {$\displaystyle T_{\initial}$};
	\draw (375,165) node [anchor=north west][inner sep=0.75pt]   [align=left] {$\displaystyle T_{\surface}$};
	\draw (333,79.6) node [anchor=north west][inner sep=0.75pt]   [align=left] {interface};
	\draw (468.8,146.6) node [anchor=north west][inner sep=0.75pt]   [align=left] {gas};
	\draw (226,147) node [anchor=north west][inner sep=0.75pt]   [align=left] {solid};

\end{tikzpicture}

}

\caption{One-dimensional model of solid propellant combustion}
\label{fig:modele1D}
\end{figure}

The temperature profile of the one-dimensional combustion of a solid propellant is presented schematically in Figure \ref{fig:modele1D}. The space variable is $x$.
The solid phase represents the propellant and is semi-infinite towards $x=-\infty$. As the surface heats up, thermal decomposition occurs inside the solid. However, due to the high activation energy of this process, the pyrolysis is often assumed to be concentrated in an infinitely thin zone at the surface, while the propellant remains inert.
A liquid layer is usually observed at the surface, however due to its low thickness and its difficult experimental characterisation,
most models assume this layer to be infinitely thin \cite{theseShihab, meynet2005simulation, jackson_2D}.
The modelled interface is the location of a singularity, where all decomposition and gasification phenomena occur.
We adopt the same assumption, however the addition of a liquid phase and other effects (e.g in-depth decomposition reactions) would not affect the considerations discussed in Section \ref{section:natureDAE}.
Flows near the surface of a solid propellant are generally at high temperatures but low velocity (1-10 m/s), therefore the gas phase is modelled as a low-Mach one-dimensional reactive flow with $n_e$ species and uniform thermodynamic pressure $P$, which is an input of the model.
Radiative emission from the gas phase and surface are neglected. The inert solid phase is simply modelled with the heat equation.

We recall the equations for both phases and their coupling at the interface. They are originally expressed in a Galilean reference frame, where the interface position $\sigma$ varies in time: $\derivshort{\sigma}{t}=-\regSpeed$ with $\regSpeed \geq 0$ the surface regression rate. We perform the simple variable change $\xcomplet = x_{galilean} - \int \derivshort{\sigma}{t} dt$, so as to keep the interface at $x=0$. This introduces a convective term in the solid phase heat equation, which represents the interface regression.
\noindent The temperature field $T$ in the solid phase at $x<0$ is subject to:
\begin{equation}
    \rhocomplet_c c_c \partialdershort{\Tcomplet}{t} + \rhocomplet_c c_c \regSpeed \partialdershort{\Tcomplet}{\xcomplet} - \partialdershort{\left(\lambda_c \partialdershort{\Tcomplet}{\xcomplet}\right)}{\xcomplet} = 0 \label{eq:base:solid}
\end{equation}
\noindent with $\rho_c$ the propellant density, $c_c$ its heat capacity, $\lambda_c$ the thermal conductivity.
Far below the surface, the solid is at its resting temperature
    $\Tcomplet({-\infty}) = \Tcomplet_\initial$. 
\renewcommand{\ucomplet}{v}
\newcommand{\hcomplet}{h}
\noindent The gas phase at $x>0$ is subject to the following partial differential equations:
\begin{curlyeqset}{1}{1pt}
    & \partialdershort{\rhocomplet}{t} + \partialdershort{\rhocomplet (\ucomplet+\regSpeed)}{\xcomplet} = 0  \label{eq:base:continuity}\\
    & \partialdershort{\rhocomplet \Ycomplet_\indexspecies}{t} + \partialdershort{\left(\rhocomplet (\ucomplet + \regSpeed) \Ycomplet_\indexspecies\right)}{\xcomplet} = - \partialdershort{J_\indexspecies}{\xcomplet} + 
	      \tauxproduction{\indexspecies}
	      \qquad \forall \indexspecies \in \llbracket 1,n_e \rrbracket
	      \label{eq:base:species}\\
& \partialdershort{\rhocomplet \hcomplet}{t} + \partialdershort{\left(\rhocomplet (\ucomplet+\regSpeed) \hcomplet\right)}{\xcomplet} = -\derivshort{P}{t} + \partialdershort{(\lambda \partialdershort{T}{x})}{x} - \partialdershort{(\Sigma_1^{n_e} \hcomplet_\indexspecies J_\indexspecies)}{x}
\label{eq:base:enthalpy}
\end{curlyeqset}

\noindent Here $Y_\indexspecies$ is the mass fraction of $\indexspecies$-th species,
and $J_{\indexspecies}=-\rho D^{eq}_{\indexspecies} \partialdershort{Y_\indexspecies}{x}$ its diffusion flux, where the equivalent diffusion coefficient $D^{eq}_\indexspecies$ is evaluated following the Hirschfelder-Curtiss approximation \cite{hirschfelder_1954}.
The volumetric production rate of the $\indexspecies$-th species is $\tauxproduction{\indexspecies}$.
The enthalpy $\hcomplet$ is the sum of the chemical and sensible enthalpies:
$\hcomplet = \Sigma_{\indexspecies=1}^{n_e} Y_\indexspecies h_\indexspecies$, where $h_\indexspecies = (\Delta h_{f,\indexspecies}^0 + \int_{T_0}^{T} c_{p,\indexspecies}(a) da)$,
with $c_{p,\indexspecies}$ the heat capacity of the $\indexspecies$-th species, and $\Delta h_{f,\indexspecies}$ its formation enthalpy at $T_0$.
The thermal conductivity is $\lambda$.
Soret 
and Dufour effects 
are neglected.
Zero-flux boundary conditions are imposed at $x=\pm \infty$: 
\begin{equation}
  \partialdershort{\Tcomplet}{\xcomplet}(-\infty) = 0, \qquad 
  \partialdershort{\Tcomplet}{\xcomplet}(+\infty) = 0, \qquad 
  \partialdershort{\Ycomplet_\indexspecies}{\xcomplet}(+\infty) = 0 \qquad \forall \indexspecies \in \llbracket 1,n_e \rrbracket 
\end{equation}

\noindent Both phases are coupled at the interface by the following conditions, expressing the continuity of the mass flow rate and temperature, as well as the thermal and species flux balance around the interface:
\begin{curlyeqset}{1}{1pt}
  &\rhocomplet_c \regSpeed = \rhocomplet(0^+) (\ucomplet(0^+)+\regSpeed) \label{eq:base:BCs:interfaceMassFlow}\\
  &\Tcomplet(0^-) = \Tcomplet(0^+) = \Tcomplet_\surface \label{eq:base:BCs:interfaceContinuityT}\\
  &\left( \massflux \hcomplet - \lambda_c \partialdershort{T}{x} 
   \right)_{0^-} =
  \left( \massflux \hcomplet
  + q_{r}
  - \lambda \partialdershort{T}{x}
  + \Sigma_1^{n_e} \hcomplet_\indexspecies J_\indexspecies  \right)_{0^+}
  \label{eq:base:BCs:interfaceBilanH}\\
  &\left(\massflux \Ycomplet_{inj,\indexspecies}\right)_{0^-} =  \left( \massflux \Ycomplet_\indexspecies + J_\indexspecies \right)_{0^+} ~~~\forall~\indexspecies ~\in~ \llbracket 1, n_e \rrbracket\label{eq:base:BCs:interfaceBilanY}
\end{curlyeqset}
\noindent with $\massflux$ the mass flow rate ($\rho_c \regSpeed$ in the solid, $\rho(\ucomplet+\regSpeed)$ in the gas), $q_{r}$ an optional external heat flux (e.g. laser) and $Y_{inj}$ the product mass fractions generated by decomposition and gasification processes, which can be constant or functions of the surface temperature as in \cite{meynet2006}.

The ideal gas law relates the various state variables in the gas phase:
\begin{equation}
\label{eq:base:idealgaslaw}
\rho = {P}/{\left( RT \sum\limits_{\indexspecies=1}^{n_e} \dfrac{Y_\indexspecies}{\molarmass_\indexspecies} \right)}
\end{equation}

\noindent with $\molarmass_\indexspecies$ the molar mass of the $\indexspecies$-th species.
The surface pyrolysis mass flow rate is given by the pyrolysis law:
\begin{equation}
\label{eq:base:pyrolysislaw}
 \massflux(0) = \rho_c \regSpeed = f(T_\surface, P)
\end{equation}
 
In the following, we denote by $u$ the gas velocity relative to the interface $(\ucomplet+\regSpeed)$, and we simplify our notations by using the mass flow rate $\massflux = \rho u$.

Note that, in this one-dimensional framework, the momentum equation is redundant, as the mass flow rate spatial variation is already determined by the continuity equation from the temporal evolution of $\rho$, which itself is known from the temporal variations of $T$ and $Y_{\indexspecies}$ and the thermodynamic pressure $P$ through the equation of state \cref{eq:base:idealgaslaw}.
Hence $\rho$ is not a true variable, although its time derivative is specified by the continuity equation. As we will see further on in Section \ref{section:natureDAE}, the continuity equation only acts as a constraint that determines the velocity field.
Including the momentum equation would only be necessary for two- or three-dimensional flow models, requiring the hydrodynamic pressure field to be computed for the solution of the velocity field.

\subsection{Discretisation with a finite volume approach}
\label{section:schemaVF}

For the numerical implementation of a solid propellant combustion model, we apply the method of lines to obtain a set of discrete evolution equations. The semi-discretisation in space is obtained with a finite volume method, however other approaches could be applied without affecting the conclusions drawn in this paper.

\subsubsection{Gas phase}

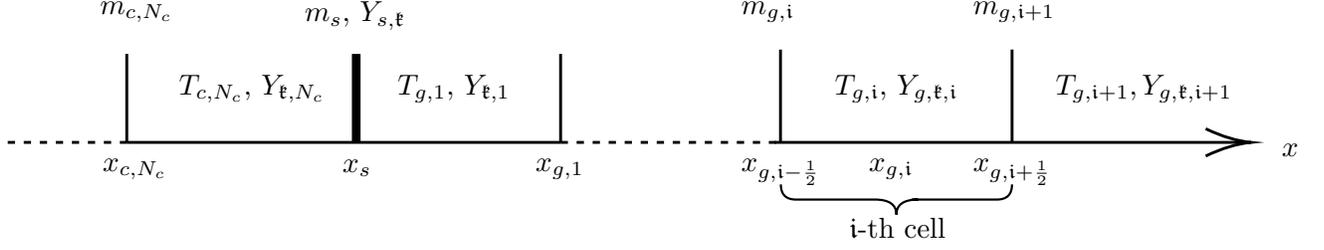
\begin{figure}
	\centering
	\newcommand{\dispcom}{\textstyle}
	
	\resizebox{\textwidth}{!}{	
		
		\tikzset{every picture/.style={line width=0.75pt}} 
		
		\begin{tikzpicture}[x=0.75pt,y=0.75pt,yscale=-1,xscale=1]
			\pgfmathsetmacro{\thicknessNormal}{1.0}
			\pgfmathsetmacro{\thicknessSurface}{3.0}
			\draw [color={rgb, 255:red, 208; green, 2; blue, 27 }  ,draw opacity=1 ][line width=3]    (204.84,70.83) -- (204.84,121.67) ;
			\draw [line width=\thicknessNormal]    (398.16,121.57) -- (614.26,121.57) ;
			\draw [shift={(618.26,121.57)}, rotate = 180] [color={rgb, 255:red, 0; green, 0; blue, 0 }  ][line width=\thicknessNormal]    (20.99,-6.32) .. controls (13.34,-2.68) and (6.35,-0.57) .. (0,0) .. controls (6.35,0.57) and (13.34,2.68) .. (20.99,6.32)   ;
			\draw [line width=\thicknessNormal]    (400.84,69.83) -- (400.84,121.57) ;
			\draw [line width=\thicknessNormal]    (507.82,69.83) -- (507.82,121.57) ;
			\draw   (401.01,147) .. controls (401.01,151.67) and (403.34,154) .. (408.01,154) -- (444.4,154) .. controls (451.07,154) and (454.4,156.33) .. (454.4,161) .. controls (454.4,156.33) and (457.73,154) .. (464.4,154)(461.4,154) -- (500.79,154) .. controls (505.46,154) and (507.79,151.67) .. (507.79,147) ;
			\draw [line width=\thicknessNormal]    (299.26,70.83) -- (299.26,121.57) ;
			\draw [line width=\thicknessNormal]  [dash pattern={on 2.53pt off 3.02pt}]  (299.26,121.57) -- (399.16,121.57) ;
			\draw [line width=\thicknessNormal]    (98.84,121.57) -- (299.26,121.57) ;
			\draw [line width=\thicknessNormal]  [dash pattern={on 2.53pt off 3.02pt}]  (43.86,121.57) -- (98.84,121.57) ;
			\draw [line width=\thicknessNormal]    (98.84,69.83) -- (98.84,121.57) ;
			\draw    (148,116) -- (148,128) ;
			\draw    (251,116) -- (251,128) ;
			\draw    (453,116) -- (453,128) ;
			
			\draw (632,117.74) node [anchor=north west][inner sep=0.75pt]   [align=left] {$\displaystyle x$};
			\draw (385.03,50) node [anchor=north west][inner sep=0.75pt]   [align=left] {$\displaystyle m_{g,\indexcell}$};
			\draw (485.17,50) node [anchor=north west][inner sep=0.75pt]   [align=left] {$\displaystyle m_{g,\indexcell+1}$};
			\draw (437,78) node [anchor=north west][inner sep=0.75pt]   [align=left] {$\displaystyle T_{g,\indexcell}$\\$\displaystyle Y_{g,\indexspecies,\indexcell}$};
			\draw (529,78) node [anchor=north west][inner sep=0.75pt]   [align=left] {$\displaystyle T_{g,\indexcell+1}$\\$\displaystyle Y_{g,\indexspecies,\indexcell+1}$};
			\draw (438.32,128) node [anchor=north west][inner sep=0.75pt]   [align=left] {$\displaystyle x_{g,i}$};
			\draw (429.33,160) node [anchor=north west][inner sep=0.75pt]   [align=left] {$\displaystyle \indexcell$-th cell};
			\draw (238,78) node [anchor=north west][inner sep=0.75pt]   [align=left] {$\displaystyle T_{g,1}$\\$\displaystyle Y_{\indexspecies,1}$};
			\draw (134,78) node [anchor=north west][inner sep=0.75pt]   [align=left] {$\displaystyle T_{c,N_{c}}$\\$\displaystyle Y_{\indexspecies,N_{c}}$};
			\draw (190.52,40) node [anchor=north west][inner sep=0.75pt]  [color={rgb, 255:red, 208; green, 2; blue, 27 }  ,opacity=1 ] [align=left] {$\displaystyle m_{s}$\\$\displaystyle Y_{s,\indexspecies}$};
			\draw (78.03,50) node [anchor=north west][inner sep=0.75pt]   [align=left] {$\displaystyle m_{c,N_{c}}$};
			\draw (195.04,128) node [anchor=north west][inner sep=0.75pt]  [color={rgb, 255:red, 208; green, 2; blue, 27 }  ,opacity=1 ] [align=left] {$\displaystyle x_{s}$};
			\draw (486.32,128) node [anchor=north west][inner sep=0.75pt]   [align=left] {$\displaystyle x_{g,\indexcell+\frac{1}{2}}$};
			\draw (379.32,128) node [anchor=north west][inner sep=0.75pt]   [align=left] {$\displaystyle x_{g,\indexcell-\frac{1}{2}}$};
			\draw (238.32,128) node [anchor=north west][inner sep=0.75pt]   [align=left] {$\displaystyle x_{g,1}$};
			\draw (132.32,128) node [anchor=north west][inner sep=0.75pt]   [align=left] {$\displaystyle x_{c,N_{c}}$};
			\draw (284.03,50) node [anchor=north west][inner sep=0.75pt]   [align=left] {$\displaystyle m_{g,2}$};

		\end{tikzpicture}
	}
	\caption{Localisation of the discretised variables in the finite volume mesh. The thick red line represents the surface.}
	\label{fig:schemaVF:mesh}
\end{figure}

The set of conservative equations for the gas phase is semi-discretised in space with a finite volume approach: the domain is split in $N_g$ cells (control volumes). The discretised variables are the temperature $T$, mass fractions $Y_\indexspecies$, mass flow rate $\massflux$.
The temperature and mass fractions are taken at the centers of each cell, while the mass flow rate is taken at the left face of each cell, which is convenient for the free boundary problem with a flux defined at the surface.
The surface temperature $T_\surface$, the surface mass fraction of the $\indexspecies$-th species $Y_{\surface,\indexspecies}$ are taken at $x=0$ (rightmost face of the solid domain, leftmost face of the gas domain). The localisation of each variable  is sketched in Figure \ref{fig:schemaVF:mesh}.

Using the notation $q$ to identify any of the conservative variables $\rho Y_\indexspecies$ and $\rho h$, or $\rho$, and the subscript $\indexcell$ as the index of the mesh cell considered, the conservative equations
\cref{eq:base:continuity,eq:base:species,eq:base:enthalpy}
become:
\begin{equation}
\label{eq:VF:dqdt}
\dx_\indexcell \deriv{q_\indexcell}{t} = -\left[ F_{d,q} + F_{c,q} \right]_{\indexcell-\frac{1}{2}}^{\indexcell+\frac{1}{2}} + \dx_\indexcell s_{q,\indexcell}
\end{equation}
\noindent with $\dx_\indexcell=x_{\indexcell+\frac{1}{2}} - x_{\indexcell-\frac{1}{2}}$ the size of the $\indexcell$-th cell, $F_{d,q}$ the diffusive fluxes, $F_{c,q}$ the convective fluxes and $s_{q,\indexcell}$ the chemical source term.
Source terms are evaluated at the cell centers.
Thermodynamic and transport properties are evaluated at the cell centers, and their values at the interfaces are taken as averages of the adjacent cells values. The gradient at the $(\indexcell-\frac{1}{2})$-th interface of a variable $q$ discretised at the cell centers is computed as:
\begin{equation}
\nabla q_{\indexcell-\frac{1}{2}} = \frac{q_{\indexcell}-q_{\indexcell-1}}{x_{\indexcell}-x_{\indexcell-1}}
\end{equation}
\noindent resulting in a second-order approximation of the diffusive fluxes.
The transported interface values at $x_{\indexcell-\frac{1}{2}}$ are defined as
$q_{\indexcell-\frac{1}{2}} = \Phi^{-}_{\indexcell-\frac{1}{2}} q_{\indexcell-1} + \Phi^{+}_{\indexcell-\frac{1}{2}} q_{\indexcell}$,
where $\Phi^-_{\indexcell-\frac{1}{2}}$ and $\Phi^+_{\indexcell-\frac{1}{2}}$ are coefficients which sum up to one.
Convective fluxes are then computed as $F_{c,q,\indexcell-\frac{_1}{^2}} =  \massflux_\indexcell  ~ q_{\indexcell-\frac{1}{2}}$.
The maximal order of accuracy is obtained with both coefficients set to 0.5, resulting in a centered convection scheme.
Although unconditionally stable with an appropriate implicit temporal discretisation, the centered scheme may result in a lack of diagonal dominance for the Jacobian used in the Newton algorithm at each step,
essentially rendering the linear system solution impossible to perform \cite{jameson1981implicit,khosla1979filtering} unless the time step is dramatically reduced.
To circumvent this issue, the centered scheme needs to be locally upwinded if the flow is convection-dominated, ensuring diagonal dominance at the expense of falling back to first-order accuracy.
This is done dynamically via a Péclet-weighted average of the first-order upwind and second-order centered schemes, similarly to what is done in \cite[page 448]{book_anderson_tannehill_pletcher} with the concept of ``mesh Reynolds number''.
The Péclet number $c_{p}m/\lambda$ measures the relative importance of diffusive phenomena compared to convective transport in the energy equation.
We define a numerical Péclet number at the $(\indexcell-\frac{1}{2})$-th interface as
$\Pe_{\indexcell-\frac{1}{2}} = \left( \Pe_{\indexcell} + \Pe_{\indexcell-1} \right) \left(x_{\indexcell}-x_{\indexcell-1}\right)/2$, 
i.e. it is an average Péclet number with a reference length taken as the distance between the centers of the neighbouring cells.
If $|\Pe_{\indexcell-\frac{1}{2}}|<0.5$, we use the centered scheme: $\Phi^+=\Phi^-=0.5$, i.e. the centered scheme is used as the thermal diffusion dominates energy convection.
If $|\Pe_{\indexcell-\frac{1}{2}}|>1$, we use the upwind scheme: $\Phi^+ = 0$ if $\Pe>0$, else $\Phi^+=1$.
The transition between these two cases is smooth with respect to $\Pe_{\indexcell-\frac{1}{2}}$, so as not to cause numerical issues later on.
Overall, the second-order scheme is used where the mesh is sufficiently refined, while the first-order scheme is only used in poorly resolved areas, typically far away from the surface, where precise representation of the flow is not required.
The final semi-discrete conservative equations in the $\indexcell$-th cell read:

\vspace{-1em}
\begin{align}
	\phantom{i + j + k}
	&\begin{aligned}
		\label{eq:VF:continuity}
		&\mathllap{ \dx_\indexcell\deriv{\rho_\indexcell}{t} }~=~& &\massflux_\indexcell - \massflux_{\indexcell+1}
	\end{aligned}\\
	&\begin{aligned}
		\label{eq:VF:species}
		&\mathllap{\dx_\indexcell\deriv{\rho_\indexcell Y_{\indexspecies,\indexcell}}{t}}~=~&
		&\massflux_{\indexcell}     \left( \Phi^+_{\indexcell-\frac{1}{2}} Y_{\indexspecies,\indexcell} + \Phi^-_{\indexcell-\frac{1}{2}} Y_{\indexspecies,\indexcell-1} \right)
		- \massflux_{\indexcell+1} \left( \Phi^+_{\indexcell+\frac{1}{2}} Y_{\indexspecies,\indexcell+1} + \Phi^-_{\indexcell+\frac{1}{2}} Y_{\indexspecies,\indexcell} \right)
		\\
		&\qquad +& &J_{\indexspecies,\indexcell-\frac{1}{2}} - J_{\indexspecies,\indexcell+\frac{1}{2}}
		+\dx_\indexcell\omega_{\indexspecies,\indexcell}
		\quad \forall \indexspecies \in \llbracket 1, n_e \rrbracket
	\end{aligned}\\
	&\begin{aligned}
	\label{eq:VF:enthalpy}
	&\mathllap{ \dx_\indexcell \deriv{\rho_\indexcell h_\indexcell}{t} }~=~&
		~&\dx_\indexcell \derivshort{P}{t}
		+ \massflux_{\indexcell}   \left( \Phi^+_{\indexcell-\frac{1}{2}} h_{\indexcell} + \Phi^-_{\indexcell-\frac{1}{2}} h_{\indexcell-1} \right)
		- \massflux_{\indexcell+1} \left( \Phi^+_{\indexcell+\frac{1}{2}} h_{\indexcell+1} + \Phi^-_{\indexcell+\frac{1}{2}} h_{\indexcell} \right)
		\\
		&\qquad -& &\frac{\lambda_{\indexcell-1} + \lambda_{\indexcell}}{2} \frac{T_{\indexcell} - T_{\indexcell-1}}{x_{\indexcell}-x_{\indexcell-1}}
		+\frac{\lambda_{\indexcell} + \lambda_{\indexcell+1}}{2} \frac{T_{\indexcell+1} - T_{\indexcell}}{x_{\indexcell+1}-x_{\indexcell}}
		\\
		&\qquad +& &\sum\limits_{\indexspecies=1}^{n_e} \left(
		\frac{h_{\indexspecies,\indexcell}+h_{\indexspecies,\indexcell-1}}{2} J_{\indexspecies,\indexcell-\frac{1}{2}}
		-\frac{h_{\indexspecies,\indexcell+1}+h_{\indexspecies,\indexcell}}{2} J_{\indexspecies,\indexcell+\frac{1}{2}}
		\right)
	\end{aligned}
\end{align}

\subsubsection{Solid Phase}
 The solid phase energy equation \cref{eq:base:solid} is replaced by a conservative equation for the enthalpy and discretised in the same way as the conservative equations in the gas phase. The solid mesh contains $N_c$ cells.
 As explained further in Section \ref{section:cost}, having a block-tridiagonal Jacobian for the complete system is computationally efficient. 
 In order to keep a consistent Jacobian structure between the gas and solid phases, the mass flow rate field $\massflux$ and the species profile $Y_{\indexspecies}$ are kept as dummy variables.
 As the solid propellant is assumed inert and incompressible, the continuity equation is equivalent to $\derivshort{\massflux}{x}=0$ which is discretised as $\massflux_\indexcell=\massflux_{\indexcell+1}$, with the boundary condition $\massflux_{N_{c}+1} = \massflux(T_\surface, P)$.
 Species evolution is not considered in the solid phase, therefore the simple equation $Y_{\indexspecies,\indexcell}=0$ is used.
 Due to the similarity in structure between the gas and solid phases semi-discrete systems, we assume when necessary in the rest of the article that the solid phase equations are contained in the previous gas phase system.
 
\subsubsection{Surface coupling conditions}
The surface variables are the surface temperature $T_\surface$ and the surface mass fractions $Y_{\surface,\indexspecies}$ on the gas side. The surface matching conditions  \cref{eq:base:BCs:interfaceBilanY,eq:base:BCs:interfaceMassFlow,eq:base:BCs:interfaceBilanH} and the continuity equation \cref{eq:base:continuity} are discretised as follows, with a first-order approximation of the gradients:
\begin{curlyeqsetat}{8}{1pt}
 0 &~=~&  &g_{th}                              &&~\coloneqq~ -&&\lambda_c \dfrac{T_\surface-T_{c,-1}}{x_\surface-x_{c,1}} + \lambda \dfrac{T_{g,1}-T_\surface}{x_{g,1} - x_\surface} + \massflux \left( \hcomplet (T_\surface,Y_{\surface,\indexspecies}) - \hcomplet (T_\surface,Y_{inj,\indexspecies}) \right) + \Sigma_1^{n_e} \hcomplet (T_\surface,Y_{inj,\indexspecies}) J_{\surface,\indexspecies} \label{eq:contraintes:energy}\\
 0 &~=~&  &g_{\speciessubscript,\indexspecies} &&~\coloneqq~  &&\massflux(T_\surface) ( Y_{inj,\indexspecies}-Y_{\surface,\indexspecies} )
 + J_{\surface, \indexspecies}
	\qquad \forall \indexspecies \in \llbracket 1, n_e \rrbracket \label{eq:contraintes:species}
\end{curlyeqsetat}

\noindent with $T_{\solid,-1}$ the temperature in the last cell of the solid phase below the surface, and $T_{\gas,1}$ the temperature in the first cell of the gas phase, just above the surface.
The species surface diffusion fluxes are computed as:
\begin{equation}
	J_{\surface,\indexspecies}
	=
	\rho(T_\surface, Y_\surface, P) D^{eq}_\indexspecies (T_\surface, Y_\surface, P)
	\dfrac{Y_{g,\indexspecies,1}-Y_{\surface,\indexspecies}}{x_{g,1} -x_{\surface}}
\end{equation}

\subsection{Final semi-discrete form}

Let $\Xdiff$ be the vector of the discrete variables $T$ (solid and gas) and $Y_{\indexspecies}$ in each cell, and let $\Xdae$ be the vector containing the remaining variables $T_\surface$, $Y_{\surface,\indexspecies}$ and $\massflux$:
$\Xdiff =(Y_{\indexspecies}, T)^t$ and $\Xdae = (T_\surface, Y_{\surface,\indexspecies},\massflux)^t$.
The semi-discretisation in space yields the following system:
\begin{curlyeqsetat}{2}{1pt}
	\derivshort{Q(\Xdiff)}{t} &~=~& f(\Xdiff, \Xdae) \label{eq:semiexplicit:conservative_ode}\\
	0 &~=~& g(\Xdiff, \Xdae) \label{eq:semiexplicit:conservative_contraintes}
\end{curlyeqsetat}
\noindent with $Q(\Xdiff)= (\rho Y_{\indexspecies}, \rho h)^t$ the vector of conserved differential variables, $f$ the corresponding spatial operator, i.e. the right-hand side of Equations \cref{eq:VF:species,eq:VF:enthalpy}, and $g$ the vector obtained by concatenating Equations \eqref{eq:contraintes:energy}, \eqref{eq:contraintes:species}, \eqref{eq:base:BCs:interfaceMassFlow} and \eqref{eq:VF:continuity}.
The continuity equation \eqref{eq:VF:continuity} is not included in Equation \cref{eq:semiexplicit:conservative_ode}, but rather in Equation \eqref{eq:semiexplicit:conservative_contraintes}, as will be explained later in Sections \ref{section:constraints:identification} and \ref{section:RK:contrainte_debit}.
A more classical form can be obtained if we were to reformulate the conservation equations so that they directly give the temporal derivatives of $W$, not of $Q(W)$:

\vspace{-0.7em}
\begin{curlyeqsetat}{2}{1pt}
	\derivshort{\Xdiff}{t} &~=~& f_1(\Xdiff, \Xdae) \label{eq:semiexplicit:ode}\\
	0 &~=~& g(\Xdiff, \Xdae) \label{eq:semiexplicit:contraintes}
\end{curlyeqsetat}
We will use the simpler first form for the numerical implementation, and we will mostly use the second formulation for theoretical discussions (dropping the subscript 1 for $f$).

\section{Differential-algebraic nature of the semi-discretised system}
\label{section:natureDAE}

We aim at developing a one-dimensional unsteady CFD tool for high-fidelity simulations of transient phenomena. Relying on a finite-volume space discretisation, we have obtained a system of semi-discrete evolution equations.
In this section, we show that a difficulty arises from the nature of this system: some variables are not defined by differential equations but by algebraic ones. The system thus belongs to the class of Differential-Algebraic Equations (DAE). We refer the reader to \cite{hairer_book2, Petzold_DAE} for details on DAEs and only the necessary aspects of this class of problem will be discussed hereafter.

\subsection{Identification of the constraints}
\label{section:constraints:identification}
The discretised surface variables $Y_{\surface,\indexspecies}$ and $T_\surface$ only appear in  equations \cref{eq:contraintes:energy,eq:contraintes:species}, however no time derivative appear. 
Such variables are called algebraic and will ``instantly'' adapt to the variations of the other variables in the cells adjacent to the surface, i.e. they are not directly affected by their time histories.
Note that this characteristic is linked to the assumption that no accumulation of energy or mass takes place at the interface, hence the latter has no inertia and its associated variables must always be such that the interface conditions are met.

The remaining algebraic equations come from the discrete gas continuity equation \cref{eq:VF:continuity}, which we recall here:
\begin{curlyeqsetat}{6}{1pt}
 0 &~=~& g_{\massflux,1} 		  &~\coloneqq~&& \massflux_1 - \massflux(T_\surface) \label{eq:contraintes:rhou_interface}\\
 0 &~=~& g_{\massflux,\indexcell} &~\coloneqq~&& \deriv{\rho_{\indexcell-1}}{t} + \dfrac{\massflux_{\indexcell} - \massflux_{\indexcell-1}}{\dx_\indexcell} \qquad \forall \indexcell \in \llbracket 2,N_{g} \rrbracket \label{eq:contraintes:rhou}
\end{curlyeqsetat}
Equation \cref{eq:contraintes:rhou_interface} is the boundary condition for the gaseous mass flow rate field and is equivalent to Equation \cref{eq:base:BCs:interfaceMassFlow}.
As underlined in Section \ref{section:baseEquations},
the density $\rho_\indexcell$ in the $\indexcell$-th cell is a function of $T_\indexcell$ and $Y_{\indexspecies,\indexcell}$ via the ideal gas law \cref{eq:base:idealgaslaw}. 
These variables are governed by
the discrete equations \cref{eq:VF:species,eq:VF:enthalpy}
and by the pressure $P$.
Therefore $\rho_\indexcell$ is not a true differential variable. Its time derivative $\derivshort{\rho_\indexcell}{t}$ appearing in Equation \cref{eq:contraintes:rhou} is entirely determined from the variations of the other gas-phase variables $T_\indexcell$ and $Y_{\indexspecies,\indexcell}$. 
The continuity equation is solely used to constrain the flow rate field in the gas phase $\massflux$, the time derivative of which does not appear in this one-dimensional low-Mach model. 
As a consequence, the discrete values of the mass flow rate are also algebraic variables, which adapt instantly to variations of the other variables in coherence with  the parabolic nature of the low-Mach number limit, 
where all pressure waves propagate at infinite speed and are relaxed instantly. This situation is generic in low-Mach number combustion modelling.

\subsection{Index of the algebraic equations}
\label{section:index}
The presence of algebraic equations usually brings in various types of additional difficulties in the numerical resolution of the system of equations. A useful mathematical criterion, the ``index'', is used to classify the different types of algebraic equations
and is essential in order to chose the proper numerical method.
It is defined as the number of times a given algebraic equation must be differentiated with respect to time to obtain a differential equation for the corresponding algebraic variable \cite{Petzold_DAE}.

Deriving our discretised Equations \cref{eq:contraintes:energy,eq:contraintes:species,eq:contraintes:rhou_interface,eq:contraintes:rhou} with respect to time, terms in $\derivshort{T_\surface}{t}$, $\derivshort{Y_{\surface,\indexspecies}}{t}$ and 
$\derivshort{\massflux_\indexcell}{t}$ appear, and manipulations lead to explicit expressions for each of them.
Apart from exceptional cases which we never encountered in practice,
the denominators in these expressions will never be zero, thus the obtained ODEs are well-posed and the index of the original DAE system is 1.
Note that when differentiating Equation \cref{eq:contraintes:rhou} with respect to time, the second-derivative $\derivshort{\rho_\indexcell}{tt}$ appears. 
It can however be expressed by differentiating the ideal gas law and the other conservation equations.

\begin{figure}
 \centering
 \includegraphics[width=0.35\textwidth]{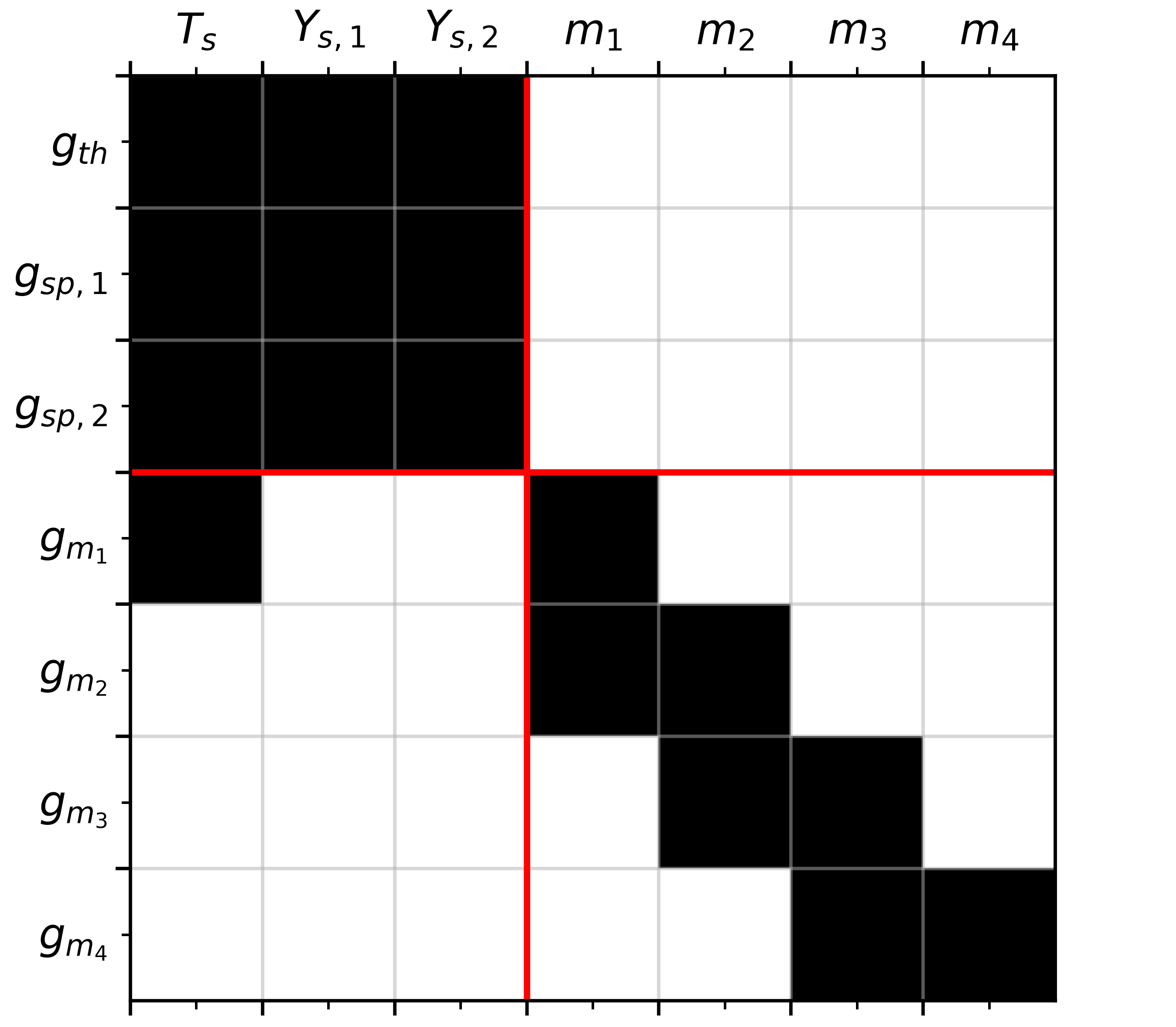}
 \caption{Sparsity pattern of the Jacobian of $g$, the vector of constraints}
 \label{fig:sparsity_contraintes}
\end{figure}

Another approach is to directly use System \cref{eq:semiexplicit:ode,eq:semiexplicit:contraintes}, also known as {\it semi-explicit form} for DAE problems, where the separation between the differential and algebraic variables $\Xdiff$ and $\Xdae$, respectively, leads to a direct identification of the constraints $g = (g_{th},g_{\speciessubscript,1} \ldots g_{\speciessubscript,n_e}, g_{\massflux_0} \ldots g_{\massflux_{N_{g}}})^t$.
In this form, the index of the DAE is 1 if the Jacobian $\partialdershort{g}{\Xdae}$ is non-singular \cite{hairer_book2}. 
We can verify this by forming this Jacobian. For the sake of simplicity, we show in Figure \ref{fig:sparsity_contraintes} its sparsity pattern when considering only 2 species. 
The labels on the vertical axis describe the constraint being derived, and the labels on the horizontal axis denote the differentiation variable. 
We can decompose the matrix into smaller blocks following the red lines.
The first block on the diagonal corresponds to the Jacobian of the nonlinear system $(g_{th}, g_{\speciessubscript,1} \ldots g_{\speciessubscript,ne})^t = 0$.
The second block on the diagonal is lower triangular with non-zero elements on the diagonal, therefore it is invertible. 

In the simplified case where only 2 species of equal properties are considered with a unitary Lewis number, and assuming $c_p=c_s$, we have been able to show that $g_{th}$ and $g_{sp}$ can be formulated such that the first block is lower-triangular with diagonal terms which are similar to the denominators previously discussed, and thus only zero in specific cases which will never be encountered in practice.
For instance, under the previous assumptions, the terms in involving species diffusion
and enthalpy in equation \eqref{eq:base:BCs:interfaceBilanH} for $g_{th}$ can be gathered under the
form $\massflux Q_p$, with $Q_p$ the heat associated to the pyrolysis and gasification processes at the surface,
which is a constant under the previous assumptions.
The diagonal term $\partialdershort{g_{th}}{T_s}$ is
$1/ ( \lambda_c/\Delta x_{c,N_c} + \lambda/\Delta x_{g,1} - Q_p \partialdershort{m}{T_s} )$, and similar expressions can be obtained for the other diagonal terms.
The term $\partialdershort{m}{T_s}$ being positive and $Q_p$ potentially being positive as well, there is no guarantee that the denominator will not cancel. However, this term can be bounded in magnitude so that, for asymptotically fine meshes, it becomes negligible, and therefore the block is ensured to invertible.

A thorough demonstration of the invertibility of this block in the general case seems out of reach, owing to the large number of nonlinear terms that appear via the diffusion fluxes, enthalpy expressions, and the pyrolysis law mostly.
However, our numerical experience shows that the Jacobian $\partialdershort{g}{\Xdae}$ is indeed invertible in all studied cases.
Overall we claim that the Jacobian is invertible, hence the index of the corresponding algebraic equations is 1.

\section{Requirements for the time integration method}
\label{section:strategy}
Now that we have determined that our semi-discretised problem is an index-1 differential-algebraic system of equations, we look for a high-order time integration method to maximise the performance of our code, while ensuring high-accuracy. Adaptive time stepping is an additional capability that we seek, so as to minimise the number of time steps and thus the computational time, while also guaranteeing a prescribed solution accuracy.
We focus on Runge-Kutta methods and shortly discuss the applicability of other methods.

\subsection{Ensuring high-order convergence for DAEs}
\subsubsection{Formulation of Runge-Kutta methods for DAEs}
A generic \totalstep-stage Runge-Kutta integration method applied to Equations \cref{eq:semiexplicit:ode,eq:semiexplicit:contraintes} yields \cite{hairer_book2}:
\begin{curlyeqset}{1}{1pt}
\Xdiffrksub_{\indexstep \indexstage} ~=~& \Xdiff_\indexstep + \Delta t \sum\limits_{\indexstagebis=1}^\totalstep a_{\indexstage \indexstagebis} f(\Xdiffrksub_{\indexstep \indexstagebis},\Xdaerksub_{\indexstep \indexstagebis}) \label{eq:RK:souspas_yni} \\
0 ~=~& g(\Xdiffrksub_{\indexstep \indexstage}, \Xdaerksub_{\indexstep \indexstage})   \label{eq:RK:souspas_zni} \\
\Xdiff_{\indexstep+1} ~=~& \Xdiffrksub_\indexstep + \Delta t\sum\limits_{\indexstage=1}^\totalstep b_\indexstage f(\Xdiffrksub_{\indexstep \indexstage}, \Xdaerksub_{\indexstep \indexstage})\label{eq:RK:pas_y} \\
\Xdae_{\indexstep+1} ~=~& (1-\sum\limits_{\indexstage,\indexstagebis=1} b_\indexstage \omega_{\indexstage \indexstagebis}) \Xdae_{\indexstep} + \sum\limits_{\indexstage,\indexstagebis=1}^\totalstep b_\indexstage \omega_{\indexstage \indexstagebis} \Xdaerksub_{\indexstep \indexstagebis} \label{eq:RK:pas_z}
\end{curlyeqset}

\noindent where $\Xdiffrksub_{\indexstep \indexstage}$ and $\Xdaerksub_{\indexstep \indexstage}$ represent the values of $\Xdiff$ and $\Xdae$ at stage $\indexstage$ (time $t_{\indexstep \indexstage} = t_n + c_i \Delta t$),
and $a_{\indexstage \indexstagebis}$, $b_\indexstage$ and $c_\indexstage$ are the coefficients of the Runge-Kutta method, and $\omega_{\indexstage \indexstagebis}$ are the coefficients of the inverse of $A=(a_{\indexstage \indexstagebis})$, which we assume to exist, i.e. we temporarily restrict ourselves to implicit methods.
At each stage, the algebraic variables $\Xdae$ are determined via Equation \cref{eq:RK:souspas_zni} such that the constraints are satisfied.
This can be interpreted as a systematic projection of the algebraic variables onto the set $\{ z \mid g(\Xdiff, z)=0\}$.
After all stages are computed, the advancement to the next time step is performed via Equations \cref{eq:RK:pas_y,eq:RK:pas_z}.

Equation \cref{eq:RK:pas_z} which does not necessarily ensure that $g(\Xdiff_{\indexstep+1}, \Xdae_{\indexstep+1})=0$, hence a deviation from the correct solution may occur. Runge-Kutta methods that are {\it stiffly accurate} satisfy the following conditions:
$\Xdiff_{\indexstep+1} = \Xdiffrksub_{\indexstep \totalstep}$, $\Xdae_{\indexstep+1} = \Xdaerksub_{\indexstep \totalstep}$, i.e. the last stage is the solution at the next time step. With such methods, we directly obtain $0 = g(\Xdiff_{\indexstep+1}, \Xdae_{\indexstep+1})$, and Equation \cref{eq:RK:pas_z} becomes redundant.
Non-stiffly accurate methods have their order of convergence severely reduced or might even become unstable when applied to index-1 DAEs \cite{hairer_book2}, whereas stiffly accurate methods retain the same order of convergence as for ODEs.
Therefore such methods are of particular interest for our problem.
Implicit Runge-Kutta methods also allow for a natural treatment of the mass conservation constraint, preserving the order of convergence on both differential and algebraic variables, as will be discussed in Section \ref{section:RK:contrainte_debit}.
Finally, fine meshes and chemical kinetics often produce stiff systems, for which L-stability is a very advantageous property. It ensures a stable solution by instantly relaxing modes with time scales much shorter than the time step.

\subsubsection{Examples}
\label{section:strategy:examplesRK}
The most widely used implicit method is implicit Euler (or backward Euler), which is single-stage, first-order, stiffly accurate and L-stable.
Applied to our system, it yields:
\begin{curlyeqset}{1}{1pt}
\Xdiff_{\indexstep+1} ~=~& \Xdiff_\indexstep + \Delta t f(\Xdiff_{\indexstep+1}, \Xdae_{\indexstep+1}) \label{eq:RK:pas_y_IE} \\
0 ~=~& g(\Xdiff_{\indexstep+1}, \Xdae_{\indexstep+1}) \label{eq:RK:pas_z_IE}
\end{curlyeqset}
\noindent A classical second-order scheme is the Crank-Nicolson method, which is not L-stable:
\begin{curlyeqset}{1}{1pt}
\Xdiff_{\indexstep+1} ~=~& \Xdiff_\indexstep + \dfrac{\Delta t}{2}\left( f(\Xdiff_{\indexstep+1}, \Xdae_{\indexstep+1}) + f(\Xdiff_{\indexstep}, \Xdae_{\indexstep})\right)  \label{eq:RK:pas_y_CKN} \\
g(\Xdiff_{\indexstep+1}, \Xdae_{\indexstep+1}) ~=~& -g(\Xdiff_{\indexstep}, \Xdae_{\indexstep}) \label{eq:RK:pas_z_CKN}
\end{curlyeqset}
Equation \cref{eq:RK:pas_z_CKN} shows that this method is sensitive to error accumulation on the algebraic variables. In particular, it is crucial that the initial condition satisfies $g(Y_0,Z_0)=0$.

\subsection{Optimising the computational cost}
\label{section:cost}
When advancing forward in time, Equations \cref{eq:RK:souspas_yni,eq:RK:souspas_zni} must be solved, usually via a Newton algorithm, which iterates on the values $\Xdiffrksub_{\indexstep \indexstage}$ and $\Xdaerksub_{\indexstep \indexstage}$ for $\indexstage \in [1,\totalstep]$. Fully implicit Runge-Kutta methods are such that all the stages must be solved simultaneously. The very popular stiffly accurate method Radau5 \cite{hairer_book2}, a 3-stage fifth-order fully implicit method based on Gauss-Radau quadrature points, is one such method. It possesses very interesting properties, however if the problem has $N$ unknowns, each time step requires solving a $3N \times 3N$ system, which can be rather costly.
An appealing subclass of Runge-Kutta methods is that of diagonally-implicit Runge-Kutta methods (DIRK) \cite{alexander1977}. These methods are such that the summations in Equations \cref{eq:RK:souspas_yni,eq:RK:souspas_zni} for the $\indexstage$-th stage only go up to $\indexstage$ instead of $\totalstep$, i.e. the stages can be solved sequentially.
Such methods require more stages (typically twice as many) to reach the same order of convergence as fully implicit methods, however a complete time step only requires the resolution of $\totalstep$ systems of size $N\times N$, which is usually more computationally efficient. Therefore we narrow down our choices to DIRK methods.

For each stage we need to solve Equations \cref{eq:RK:souspas_yni,eq:RK:souspas_zni} for the unknowns $\Xdiffrksub_{\indexstep \indexstage}$ and $\Xdaerksub_{\indexstep \indexstage}$. These equations can be combined to form the nonlinear problem $F(X) = 0$, 
with $X=(\Xdiffrksub_{\indexstep \indexstage}, \Xdaerksub_{\indexstep \indexstage})^t$ the vector of unknowns. This problem is solved iteratively using a quasi-Newton method.
\modif{
At each Newton step, the Newton increment $\Delta X$ is obtained by solving the linear system $J \Delta X = \modif{-}F(X)$, with $J=\partialdershort{F}{X}$ the jacobian of the residual vector.
A more detailed description of the structure of $J$ is given in \ref{appendix:jacobienne}.
It is computed by finite differences and is only updated when convergence is poor. We did not implement a way to evaluate and store the separate Jacobians of $f$, $Q$ and $g$ which appear in $J$, hence the latter cannot be simply refactored upon time step changes. In the case of nearly linear dynamics, this can be disadvantageous. This matter will be discussed in Section \ref{sec:limitcycle:comparison_fixed_dt}.}

\modif{
Following Section \ref{section:schemaVF}, in particular the three-point stencil of the spatial discretisation and the ordering of our state vectors (variables sorted by cell)}, $J$ is block-tridiagonal and a Thomas algorithm can be used, after having performed a  LU-decomposition.
If the diagonal coefficients $a_{\indexstage \indexstage}$ of the method are not equal, $J$ needs to be updated at each stage.
Therefore, we focus on singly-diagonally implicit Runge-Kutta methods (SDIRK), which satisfy the property $a_{\indexstage \indexstage} = a_{\indexstagebis \indexstagebis} ~\forall~ (\indexstage,\indexstagebis) \in \llbracket 1, \totalstep \rrbracket$.

\subsection{Time adaptation}
\label{section:RK:adaptation}
A typical ignition transient is shown in Figure \ref{fig:allumage:Ts_esdirk} for a solid propellant ignited by a laser source. The evolution of the surface temperature is very rapid at the beginning of the propellant heating and around the time of ignition. The ignition transient can clearly be split into successive phases with very different time scales:
inert heating, ignition {\it per se} and stabilisation onto steady-state.
Consequently, ensuring the time step is dynamically adapted is important to improve accuracy and lower the computational time.

\newcommand{\ylow}{\ensuremath\hat{X}}
\newcommand{\yhigh}{\ensuremath X}
Some simulations presented in the literature use a constant time step, taken as sufficiently low compared to the ignition time \cite{Erikson1998, Liau_RDX_ignition}. When time adaptation is used, it usually relies on a CFL criterion \cite{theseAllumageSmith2011}. CFL limitation may be irrelevant before ignition, as the gas-phase flow velocity is negligible, therefore it may be supplemented with an additional requirement that the relative solution variation between two successive time steps is sufficiently small (e.g. 1\%). However such an approach requires fine-tuning and does not provide any guarantee regarding the accuracy of the solution.
So-called {\it embedded Runge-Kutta methods} provide a local error estimate by comparing two solutions $\yhigh_{\indexstep+1}$ and $\ylow_{\indexstep+1}$ of different orders $p$ and $\hat{p}<p$, which are constructed to share as many internal stages as possible.
This allows to cheaply obtain an objective estimate $\epsilon(\Delta t)$ of the integration error that scales as $\mathcal{O}(\Delta t^{\hat{p}+1})$.
We can the compute a normalised error $err(\Delta t) = || \epsilon(\Delta t) / (atol + rtol |X_\indexstep|) ||$, where the division is performed element by element, and $atol$, $rtol$ are user-specified absolute and relative error tolerances. In the present article, we use the 2-norm.
Requiring $err(\Delta t) < 1$ yields an explicit way to determine an optimal time step $\Delta t_{opt}$ that satisfies this requirement \cite{hairer_book1}:
$\Delta t_{opt} = \Delta t . err(\Delta t)^{-1/(q+1)}$.

If $\Delta t$ is such that the asymptotic regime of convergence for the integration method is already reached, then the estimated local integration error for $\Delta t_{opt}$ will be close to the prescribed tolerance.
If the estimated error is larger than the tolerance, the current step is restarted with the new (smaller) time step $\Delta t_{opt}$. Otherwise, it is accepted and the next step is computed with the new optimal time step length $\Delta t_{opt}$.
This allows the time step to be dynamically reduced or increased, ensuring that the solution is solved at least as precisely as specified, while minimising computational cost.
In practice, to avoid over-correcting the time step, we do not allow it to change by a factor lower than 0.2 or higher than 5 between two successive error estimations. A safety factor of $0.9$ is also applied to $\Delta t_{opt}$ to ensure the tolerance is strictly satisfied.
If a time step fails due to floating arithmetic errors or results in non-convergence of the Newton iterates, the current step is started over again with a decreased step length.
Application of a dead zone, i.e. not changing the time step if its required relative variation is below a certain threshold, helps lowering the number of time step changes and Jacobian updates, however it generally leads to worse performance overall in our test cases.

\subsection{Final choice of the method}
\label{section:choice}
Considerations on the accuracy of the method for index-1 DAEs have led us to consider stiffly accurate Runge-Kutta methods. Minimisation of the computational cost due to Newton iterations is ensured by using singly-diagonally implicit methods, and embedded methods allowing for robust time adaptation are favoured.
Additionally, improvements in the error estimation for stiff system and DAEs can be obtained if the lower-order embedded method is stiffly accurate as well, as discussed in \cite{Kvaerno2004}. This reference introduces several such schemes with an additional interesting property, which is that the first stage is the last stage of the previous step. This allows to have a ``free'' stage to improve the accuracy of the method without any additional cost. We retain three schemes from this reference: ESDIRK-32A, ESDIRK-43B and ESDIRK-54A.
For comparison with more classical schemes, we also include backward Euler and Crank-Nicolson. All methods and their properties are listed in Table \ref{table:rk_integrators}.

\begin{table}
	\begin{center}
		\begin{tabular}{c|c|c|c}
			{\bf name} & {\bf number of stages} & {\bf L-stable} & {\bf order }\\
			\hline
			ESDIRK-54A & 7 (6) & yes & 5\\
			\hline
			ESDIRK-43B & 5 (4) & yes & 4\\
			\hline
			ESDIRK-32A & 4 (3) & yes & 3\\
			\hline
			CKN (Crank-Nicolson) & 2 (1) & no & 2\\
			\hline
			IE (Implicit Euler) & 1 & yes & 1
		\end{tabular}
	\end{center}
	\caption{\label{table:rk_integrators}Selected Runge-Kutta methods. The numbers in brackets correspond to the number of stages actually solved.}
\end{table}

\subsection{Alternative choices}

From a practical point of view, if an existing code uses an implicit Euler or Crank-Nicolson scheme, as is often the case in the literature, the implementation of ESDIRK methods is straightforward, as it only requires solving more stages.
This is also the case if the original code uses a multistep BDF method. Once the ESDIRK schemes are properly implemented, adding time adaptation is also straightforward.

Other methods may however be of interest. 
Multistep methods can be applied directly to the semi-explicit form of the system. In particular the DASSL algorithm \cite{Petzold_DAE} has been extensively applied to index-1 DAEs. However, it is well known that their stability properties degrade \modif{for orders higher than 2}.
Rosenbrock methods \cite{hairer_book2} are similar to SDIRK methods, but are formulated such that only linear systems need to be solved at each stage. It is however possible that strong system nonlinearities cause severe time step restrictions \cite{ostermann_1986}, and these methods may also be less precise with regards to the satisfaction of the algebraic constraints \cite{hairer_book2}.

Another possibility is to take advantage of the fact that the DAE system is of index 1, thus the algebraic variables can be uniquely determined from the differential variables, i.e. we can consider that there exists a function $p$ such that $\Xdae=p(\Xdiff)$ satisfies $g(\Xdiff,\Xdae)=0$. Usually this relation is not explicitly known and $\Xdae$ may instead be iteratively determined via a Newton method.
The differential variables $\Xdiff$ are then governed by the ODE $\derivshort{\Xdiff}{t}=f(\Xdiff,p(\Xdiff))$, which can be integrated with any ODE solver, in particular explicit ones.
This approach is also referred to as \textit{state-space form} \cite{hairer_book2}.
In practice, we have observed that our fully implicit approach is able to produce accurate results while having CFL numbers much higher than 1, therefore explicit integration algorithms would be relatively inefficient due to their stability limitations.
This can be overcome by using an implicit algorithm to integrate $f$, however the cost of solving for the algebraic variables at each evaluation of $f$ may become prohibitive.
All operators (convection, diffusion, reaction) are stiff, thus partially implicit algorithms, e.g. IMEX methods \cite{IMEX_ascher}, are not suited.

Splitting methods could be used so that each operator is integrated with an adequate efficient method, however the order of accuracy in time would generally not exceed 2. Additional difficulties may appear when handling the interface conditions, and time step adaptation is more involved compared to embedded methods \cite{descombes_2016}.
As already mentioned, fully implicit methods like Radau5 \cite{hairer_book2} could be used, however the linear algebra becomes more complex as larger systems need to be solved.
An interesting alternative could be the parallel DIRK iterations of a fully implicit Runge-Kutta method \cite{van1991iterated}, where multiple stages can be solved in parallel to gradually approach the solution of the fully implicit method.

Finally, it should be pointed out that it is difficult to use existing solvers directly.
For the sake of clarity we used the classical form \eqref{eq:semiexplicit:ode} for the differential part in this section,
but Equation \cref{eq:semiexplicit:conservative_ode} actually prescribes the time derivatives of the conserved variables $Q(W)=(\rho Y_\indexspecies, \rho h)$.
Thus, the differential part is of the form $d_t Q(\Xdiff) = f(\Xdiff, \Xdae)$, which is not directly compatible with existing ODE/DAE solvers.
A first remedy would be to rewrite the conservation equations so that they directly give $d_t Y_\indexspecies$ and $d_t T$, in a similar manner to \ref{appendix:contrainte:instant}.
This would yield Equations \cref{eq:semiexplicit:ode,eq:semiexplicit:contraintes}.
Another solution is to discretise the conserved variables instead, and use a Newton method at each call of $f$ to compute $T, Y_k$ from $Q=(\rho, \rho Y_\indexspecies, \rho h)^t$. Thus, Equation \eqref{eq:semiexplicit:conservative_ode} would implicitly take the form $d_t Q = f(\Xdiff(Q),\Xdae)$.
Both solutions would add complexity to the code and are therefore not used in the rest of the article.
Also, since the continuity equation \cref{eq:VF:continuity} gives the time derivative of $\rho$ which is not one of our discretised variables, this Equation would need to reformulated as in Section \ref{section:RK:contrainte_debit:instant}, so as to transform it into a constraint of the $0=g_\massflux (\Xdiff, \Xdae)$. However, as discussed in the next section, this reduces the computational efficiency.

\subsection{Handling the continuity equation}
\label{section:RK:contrainte_debit}

In our one-dimensional low-Mach approach, the density $\rho$ cannot be considered a true variable of our problem, as it is uniquely determined from the temperature, mass fractions and pressure.
Consequently, the continuity equation \cref{eq:VF:continuity} should not be considered as an ODE on $\rho$, but rather as a constraint on the mass flow rate field.
We now focus on the handling of this specific problem.
In particular, we show that a Runge-Kutta formulation can be used in order to provide an original and natural treatment of the mass flow rate constraint.

\subsubsection{General formulation of our DAE problem}

We previously used Equations \cref{eq:semiexplicit:ode,eq:semiexplicit:contraintes} to represent our complete semi-discrete system.
The continuity equation \eqref{eq:VF:continuity} was included in Equation \eqref{eq:semiexplicit:conservative_contraintes}.
We can instead use the following more detailed formulation to highlight the particular characteristic of this equation:
\begin{curlyeqset}{1}{1pt}
\partialdershort{y}{t} &= f(y,z_1,z_2) \label{eq:dae:generic:ode}\\
0 &= \partialdershort{\phi(y)}{t} + g_1(y,z_1,z_2) \label{eq:dae:generic:odae}\\
0 &= g_2(y,z_2) \label{eq:dae:generic:dae}
\end{curlyeqset}
with index-1 constraints $g_2$ and, in our case, $y=(T,Y_k)$, $z_1=(\rho u)$, $z_2=(T_s, Y_{s,k})$, and $\phi(\modif{y})=\rho$ as per the equation of state.
The constraint \cref{eq:dae:generic:odae} \modif{corresponds to} the semi-discrete continuity equation \cref{eq:contraintes:rhou}, while the constraint \cref{eq:dae:generic:dae} \modif{gathers} the interface
coupling conditions \modif{\cref{eq:contraintes:species}, \eqref{eq:contraintes:energy} and \cref{eq:contraintes:rhou_interface}}. Equation \cref{eq:dae:generic:ode} contains the remaining differential equations on $\rho h$, $\rho Y_k$\modif{, i.e. the semi-discrete equations \cref{eq:VF:species,eq:VF:enthalpy}}.
This form clearly shows that the time derivative appearing in \eqref{eq:dae:generic:odae} is not that \modif{of any of the discretised variables $y$, $z_1$ or $z_2$. Hence, it is impossible} to directly apply existing DAE solvers to this system, as discussed in the previous section.
Multiple time discretisation approaches for this particular equation can be envisioned and are described next.

\subsubsection{Instantaneous reformulation of the first constraint}
\label{section:RK:contrainte_debit:instant}
Equation \cref{eq:dae:generic:odae} can be transformed so that the
time derivative of $\phi$, which is not a true variable of our system,
does not appear. Writing:
$\partialdershort{\phi(y)}{t} = (\partialdershort{\phi}{y}) \derivshort{y}{t}$,
\noindent
Equation \cref{eq:dae:generic:odae} can be replaced by:
\begin{equation}
\label{eq:dae:generic:instantaneous_constraint}
0 = (\partialdershort{\phi}{y}) f(y,z_1,z_2) + g_1(y,z_1,z_2)
\end{equation}

First, we may apply a Runge-Kutta scheme on
\cref{eq:dae:generic:ode} to compute the evolution of $y$.
At each evaluation of $f$, $z_2$ is obtained
by solving Equation \cref{eq:dae:generic:dae}, then $z_1$ is computed by solving
Equation \cref{eq:dae:generic:instantaneous_constraint}. Thus, $f$ can be evaluated.
The detailed discrete equations for our solid propellant system are derived in \ref{appendix:contrainte:instant}.

This state-space form is often used in the literature for fractional-step approaches, in particular in low-Mach or incompressible flow models to compute the perturbed pressure field by solving a Poisson equation to ensure the velocity field remains divergence-free \cite{MOTHEAU2016430,knio99,najm05}.
However in our 1D case, as discussed in Section \ref{section:choice}, we believe it is best to solve for the algebraic and differential variables at the same time, i.e. use a fully-coupled implicit approach, \modif{since the cost of solving the velocity field is small compared to that of the species and energy}.

\subsubsection{Use of the Runge-Kutta temporal quadrature}
Equation \cref{eq:dae:generic:odae} constitutes an ODE on $\phi(y_1)$,
which is not a true variable of our problem,
but is computed from the other ones via the equation of state.
We can nonetheless apply the Runge-Kutta scheme to this equation. We
then obtain, for the $i$-th stage of a DIRK method:

\begin{curlyeqset}{1}{1pt}
Y_{ni} &= y_{n} + \Delta t \sum\limits_{j=1}^{i} a_{ij} f(Y_{nj}, Z_{1,nj}, Z_{2,nj})\\
\phi(Y_{ni}) &= \phi(y_{n}) + \Delta t \sum\limits_{j=1}^{i} a_{ij}
g_1(Y_{nj}, Z_{1,nj}, Z_{2,nj}) \label{eq:dae:generic:rk_const}\\
0 &= g_2(Y_{ni},Z_{2,ni})
\end{curlyeqset}
Equation \cref{eq:dae:generic:rk_const} constitutes an algebraic constraint on $Z_{1,ni}$.
The detailed equations are derived in \ref{section:appendix:drhodt_ordre_eleve}.

\subsubsection{Comparing both approaches}
We can interpret the previous formulation as performing a quadrature on the
continuity equation, i.e. approximate the integral in the exact
solution:
\begin{equation}
	\phi(Y_{ni}) = \phi(y_n) + \int_{t_n}^{t_{n}+c_i\Delta t} g_1(y(t), z_1(t), z_2(t)) dt
\end{equation}
The quadrature is the same as the one used for the ODE
\cref{eq:dae:generic:ode} and has the same order of accuracy in time.
On the other hand, the instantaneous reformulation yields:
\begin{equation}
\label{eq:dae:generic:instantaneous_constraint_stage}
0 = (\partialdershort{\phi}{y})(Y_{nj}) f(Y_{nj},Z_{1,nj},Z_{2,nj})
+ g_1(Y_{nj},Z_{1,nj},Z_{2,nj})
\end{equation}

The quadrature performed in Equation \eqref{eq:dae:generic:rk_const} 
introduces an error $O(\Delta t^p)$ with $p$ the accuracy order of the stage considered (2 for internal stages of ESDIRK methods, and the overall order of the method for the last stage of a stiffly accurate method). An error of the same order is introduced in the instantaneous formulation \cref{eq:dae:generic:instantaneous_constraint_stage} through $Y_{nj}$ and $Z_{2,nj}$.

It was assessed on the test cases of this article that both approaches yield virtually identical results, in terms of dynamics, time step evolution, orders of convergence in time and space.
Still, the instantaneous formulation yields a more nonlinear system
because of the term $\partialdershort{\phi}{y}$.
Therefore, convergence of the Newton algorithm is more challenging and we have found that, in many cases,
this results in a dramatic increase in the number of Jacobian updates, roughly $30$\% more in highly dynamic simulations (ignition), and up to a
factor 10 in certain situations (near steady-state AP combustion from Section \ref{section:chimie_complexe}), substantially slowing down the computation.
Thus our original approach is more robust and performant, while also being easier to derive and implement.
To our knowledge, the Runge-Kutta quadrature approach has not yet been presented in the context of one-dimensional low-Mach flows,
and we claim that it is both simpler and more efficient than the traditional instantaneous reformulation.
All the results presented in this paper were obtained with this approach.

\section{Numerical verification}
\label{section:verification}
The previous numerical strategy has been implemented in the \textsc{Vulc1D} Fortran code.
In this section, we verify the spatial and temporal discretisations by comparing the results to those obtained with other approaches.

\subsection{Steady-state solution}
\label{section:simplemodel}
In the case of temperature-independent properties in both phases and unitary Lewis number in the gas phase,
we have developed a semi-analytical tool \cite{articleCTM}, which solves the complete problem in steady-state with arbitrary precision. It can be used as reference to assess the steady-state solution produced by our CFD code.
This comparison has been presented in details in \cite{articleCTM} and we only summarise the most important aspects here for the sake of completeness.

\subsubsection{Model parameters}
\label{section:simplemodel:parameters}
We consider a simple one-dimensional model of an AP-HTPB-Al propellant, whose main characteristics are summarised hereafter.

\paragraph{Solid phase}
The solid phase is composed of the solid species $P$ and has the following properties:
$\rho_c = 1806$ kg.m$^2$, $\Delta h_f^o(P) = 0$ J/kg at $T=0$ K, $c_c = 1253$ J/kg/K, $\lambda_c = 0.65$ W/m/K.
The initial temperature is $T_\initial = 300$ K.

\paragraph{Surface}
The pyrolysis mass flow rate is computed as:
$\massflux = A_p \exp(-T_{ap}/T_\surface)$, with $A_p=6.07 \times 10^7$ kg/s/m$^2$ and $T_{ap}=15082$ K.
The pyrolysis process converts the solid phase into the gaseous species $G_1$.

\paragraph{Gas phase}
Two global species are considered: the reactant $G_1$ and product $G_2$, which have the same properties except standard enthalpies. Their molar mass is $\molarmass = 74$ g/mol, and their heat capacities are
$c_p=c_c$. The standard enthalpies at $T=0$ K are $\Delta h_f^o (G_1) = -1.80 \times 10^{5}$ J/kg and $\Delta h_f^o (G_2) = -4.06 \times 10^{6}$ J/kg.
The unique global reaction is $G_1 \rightarrow G_2$ and is irreversible. The reaction rate is computed as:
$\omega = A [G_1] \exp(-T_a/T)$, with $A=435.5$ s$^{-1}$, $T_a=7216$ K, and $[G_1]$ the concentration of $G_1$.
The diffusion coefficients are equal for both species and taken as a linear function of $T$ such that the Lewis number is one throughout the gas phase. The thermal conductivity is $\lambda = 0.464$ W/m/K.

\subsubsection{Verification process and results}
First, the complete steady-state problem is solved with the semi-analytical tool.
Multiple meshes are then generated for the one-dimensional CFD code: knowing the temperature profile from the semi-analytical solution and starting from $x=0$ (surface), grid points are placed such that the difference in interpolated temperature between two successive grid points is
equal to or below a given threshold. By varying this threshold (from 0.05K to 50 K), grids with
varying levels of refinement are obtained, whose point distribution is relatively well
adapted to the problem. The finite volume mesh is then generated by taking these
grid points as the positions of the cell faces.
The generated meshes are then extended with a geometric progression of cell sizes so that their outer limits are far beyond the characteristic thermal length scales ($\approx 10^{-4}$ m here).
The semi-analytical solution is taken as the initial state and advanced forward in time with backward Euler and large time steps, until stabilisation of the CFD solution.

We show in Figure \ref{fig:validation:steady_convergence} the convergence of the CFD result towards the semi-analytical
solution for the reference case, as a function of the number of cells.
Second-order accuracy in space is attained, and the relative errors reach $10^{-8}$ on $T_\surface$ at around 4000 adapted mesh cells, and similar results are obtained on temperature profiles.
Pushing the error below $10^{-8}$ on $T_\surface$ is difficult since it is very close the maximum accuracy achievable by the Newton algorithm in double-precision arithmetic.
Still, the error is sufficiently small so that we can consider the CFD solution spatially converged.
The steady-state solutions are coherent between both approaches, thus verifying our spatial semi-discretisation.

\begin{center}
\begin{minipage}{0.45\textwidth}
	\centering
	\includegraphics[width=\textwidth]{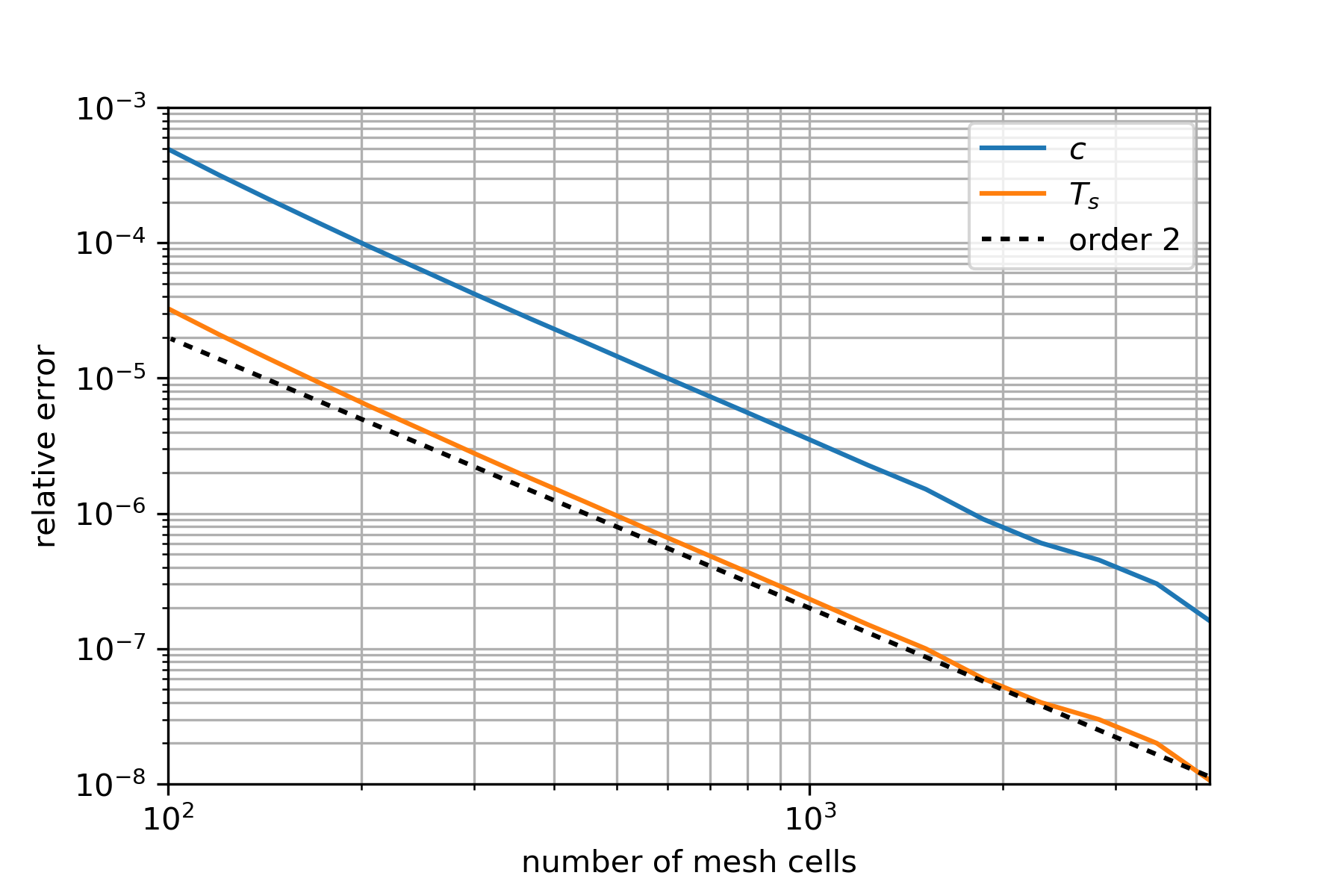}
	\captionof{figure}{Convergence of the steady-state CFD solution towards the semi-analytical solution \cite{articleCTM}}
	\label{fig:validation:steady_convergence}
\end{minipage}
~
\begin{minipage}{0.05\textwidth}
~
\end{minipage}
~
\begin{minipage}{0.45\textwidth}
\centering
\includegraphics[width=\textwidth]{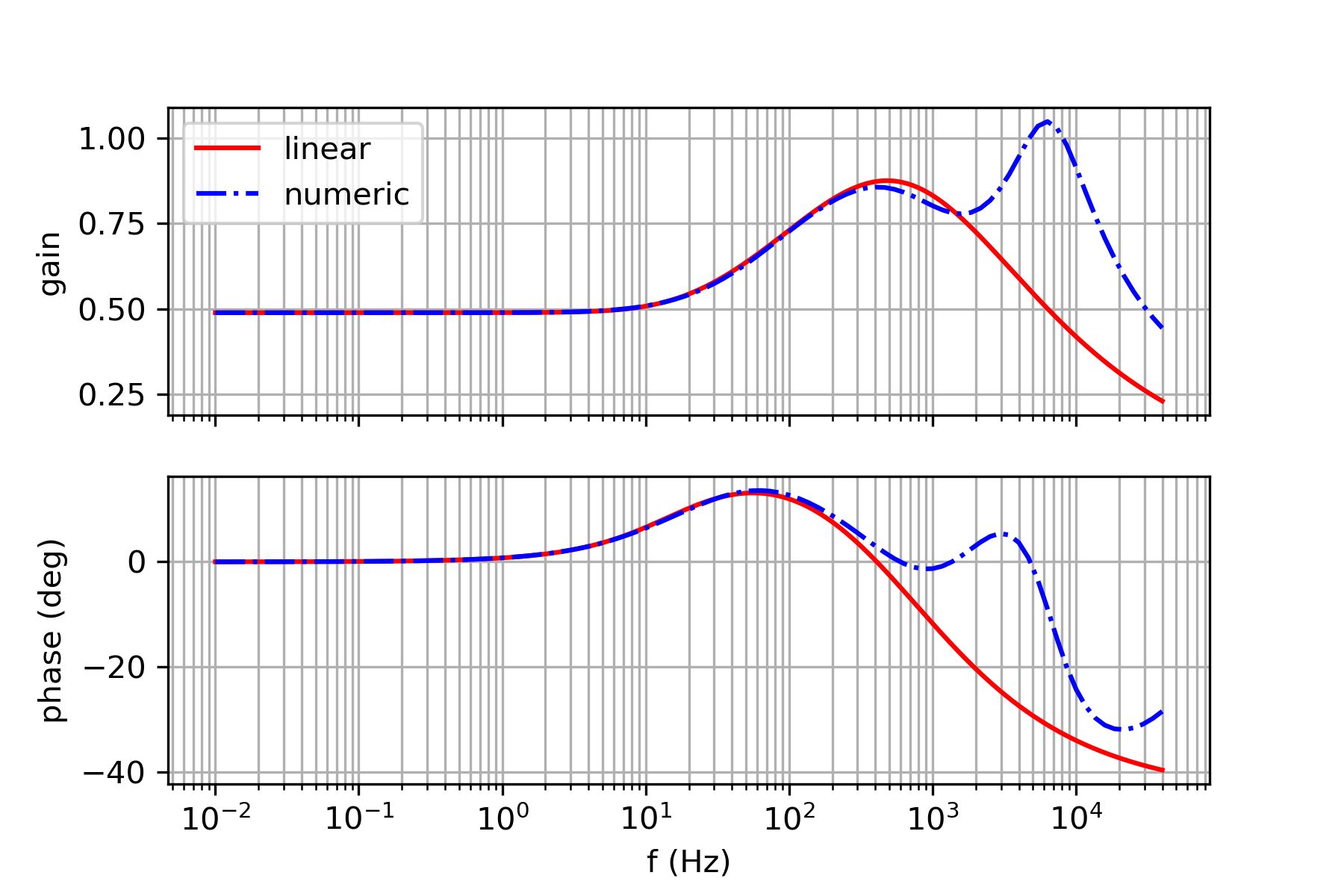}
\captionof{figure}{Bode diagram of the response function $R_{mp}$ to pressure fluctuations}
\label{fig:reponse_pression}
\end{minipage}
\end{center}

\subsection{Verification in unsteady regime}
We now wish to verify the time integration algorithm implemented in the CFD code.
No analytical solution is available in the general unsteady regime, however linear frequency responses to pressure oscillations are available \cite{reponse_culick_1968}, assuming a quasi-steady gas phase.
The linearised response can be constructed analytically using coefficients that represent the sensitivity of the steady-state solution to certain parameters.
These sensitivities can be computed via finite differences.
The 1D CFD code is initialised at the steady-state profile and a sinusoidal pressure oscillation $P(t) = \overbar{P} + p' \sin(2\pi f t)$ is enforced in the gas phase, with $f$ a given frequency.
The mean pressure $\overbar{P}$ is set to $=50\times 10^5$ Pa, and the amplitude is set to $p'=0.001 \overbar{P}$.
In the linear regime ($p'$ small), this leads to sinusoidal fluctuations of amplitude $m'$ in the pyrolysis mass flow rate $\massflux(t) = \overbar{\massflux} + \massflux'\sin(2\pi f t + \phi)$, with $\phi$ the phase shift.
After a few periods, these oscillations stabilise and we can determine the response function $R_{mp} = ({\massflux'/\overbar{\massflux})}/({p'/\overbar{P}})$ at the corresponding frequency.

Figure \ref{fig:reponse_pression} shows the comparison of the linearised and numerical frequency responses. We see that the agreement between both methods is excellent up to approximately 500 Hz, where the gas phase no longer has a quasi-steady behaviour, thus introducing a larger error in the linear response function. This serves as a global verification of our unsteady model. The secondary peak at high frequencies in the response can be obtained analytically if the unsteady gas phase equations are also linearised, as in \cite{clavin_reponse_1992}, however this is a much more involved process. It can also be removed by enforcing the quasi-steadiness of our gas-phase model.

An additional verification of the orders of convergence in time for unsteady simulations is presented in \ref{appendix:convergence}, thus completing the verification process in terms of both spatial and temporal discretisations. In particular, algebraic constraints do not hinder the high-order convergence. We now wish to tackle three much more challenging test cases and investigate the behaviour of the proposed strategy in terms of accuracy and computational efficiency.

\section{Simulation of ignition transients}
\label{section:allumage:simplifie}

\subsection{Setup}
We use the same simplified model as previously described in Section \ref{section:simplemodel:parameters}. The initial solution is a uniform temperature field at 300 K, with only combustion products in the gas phase, as a simpler alternative to adding nitrogen as initial gas, without much effect on the ignition process itself.
We add a laser heat flux $q_{r} = 1$ MW.m\textsuperscript{-2} absorbed at the surface via Equation \cref{eq:base:BCs:interfaceBilanH}.
The mesh has 99 cells in the solid phase and 291 in the gas phase. The cells are distributed such that the steady-state temperature profile is well resolved.

We compute the ignition transient with the methods listed in Table \ref{table:rk_integrators}.
The ESDIRK methods use the time adaptation strategy presented in Section \ref{section:RK:adaptation}, while the classic schemes implicit Euler (IE) and Crank-Nicolson (CKN), with a time adaptation based on the requirement that the solution has a relative variation that is below a certain value between two consecutive time steps. A discussion on a CFL-based time adaptation is presented at the end of this section.
We use the abbreviation $rtol$ to refer to the relative integration error tolerance for ESDIRK methods, and to the allowed relative variation of the solution between consecutive time steps for IE and CKN.
The maximum time step allowed is 0.1 s.

\subsection{Physical interpretation of the ignition}
Figure \ref{fig:allumage:Ts_esdirk} shows the evolution of the surface temperature during ignition as computed by ESDIRK-54 with $rtol=10^{-6}$. The first phase is the inert heating of the solid propellant. The constant laser heat flux produces an evolution of $T_\surface$ which is proportional to $\sqrt{t}$.
When $T_\surface$ is sufficiently high, the pyrolysis mass flow rate given by \cref{eq:base:pyrolysislaw} increases rapidly, causing the release of gaseous pyrolysis products in the gas phase, which chemically react and form a flame that heats up the solid even more. Typically at this point, more thermal energy is stored in the solid as compared to steady-state. This results in a momentarily higher regression rate at ignition, seen here near $t = 0.35$ s, which evacuates this excess of solid phase thermal energy. The temperature profile then converges the steady-state solution.

\subsection{Result}

\begin{center}
	\begin{minipage}{0.45\textwidth}
		\centering
		\includegraphics[width=\textwidth]{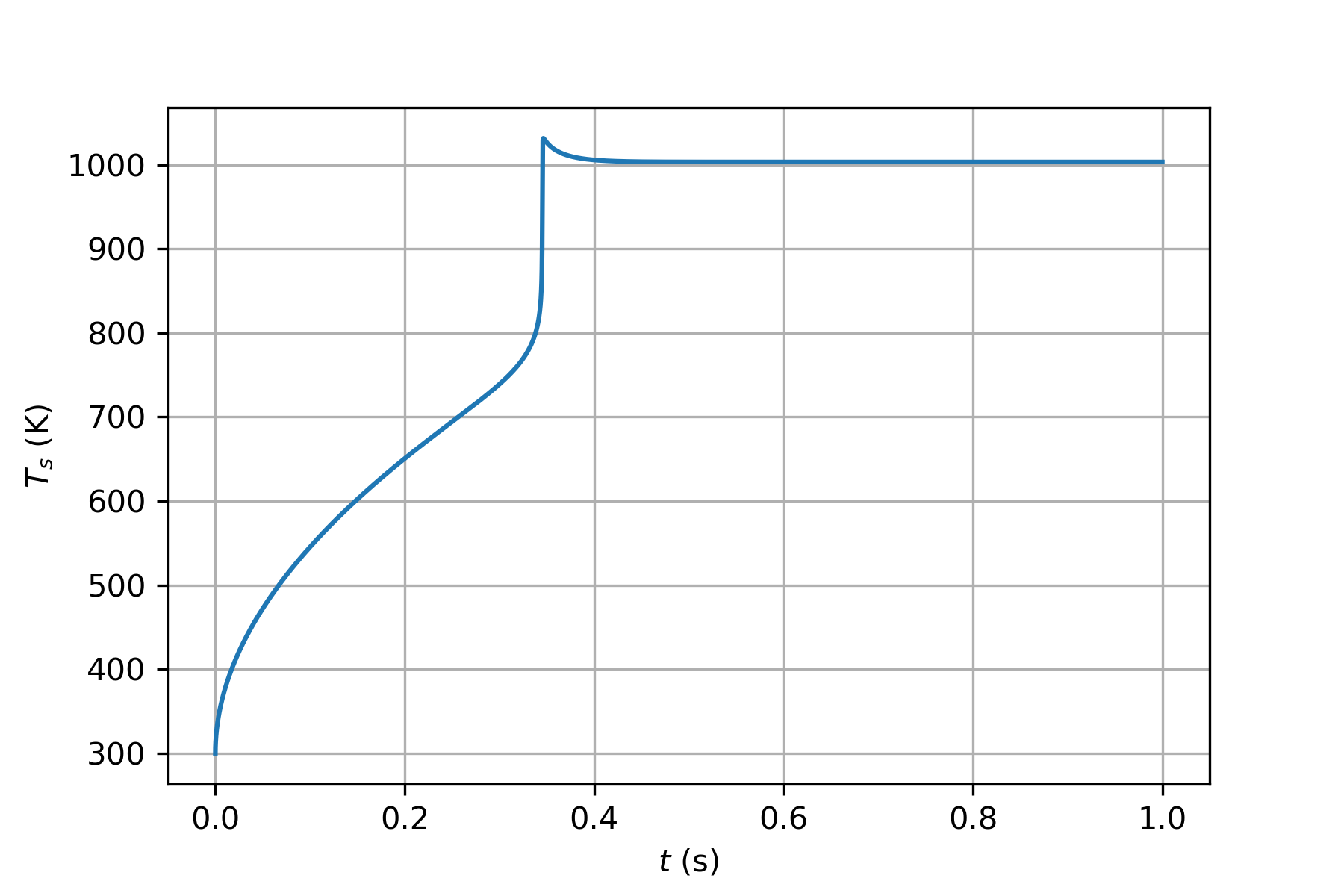}
	\end{minipage}
	~
	\begin{minipage}{0.05\textwidth}
		~
	\end{minipage}
	~
	\begin{minipage}{0.45\textwidth}
		\centering
		\includegraphics[width=\textwidth]{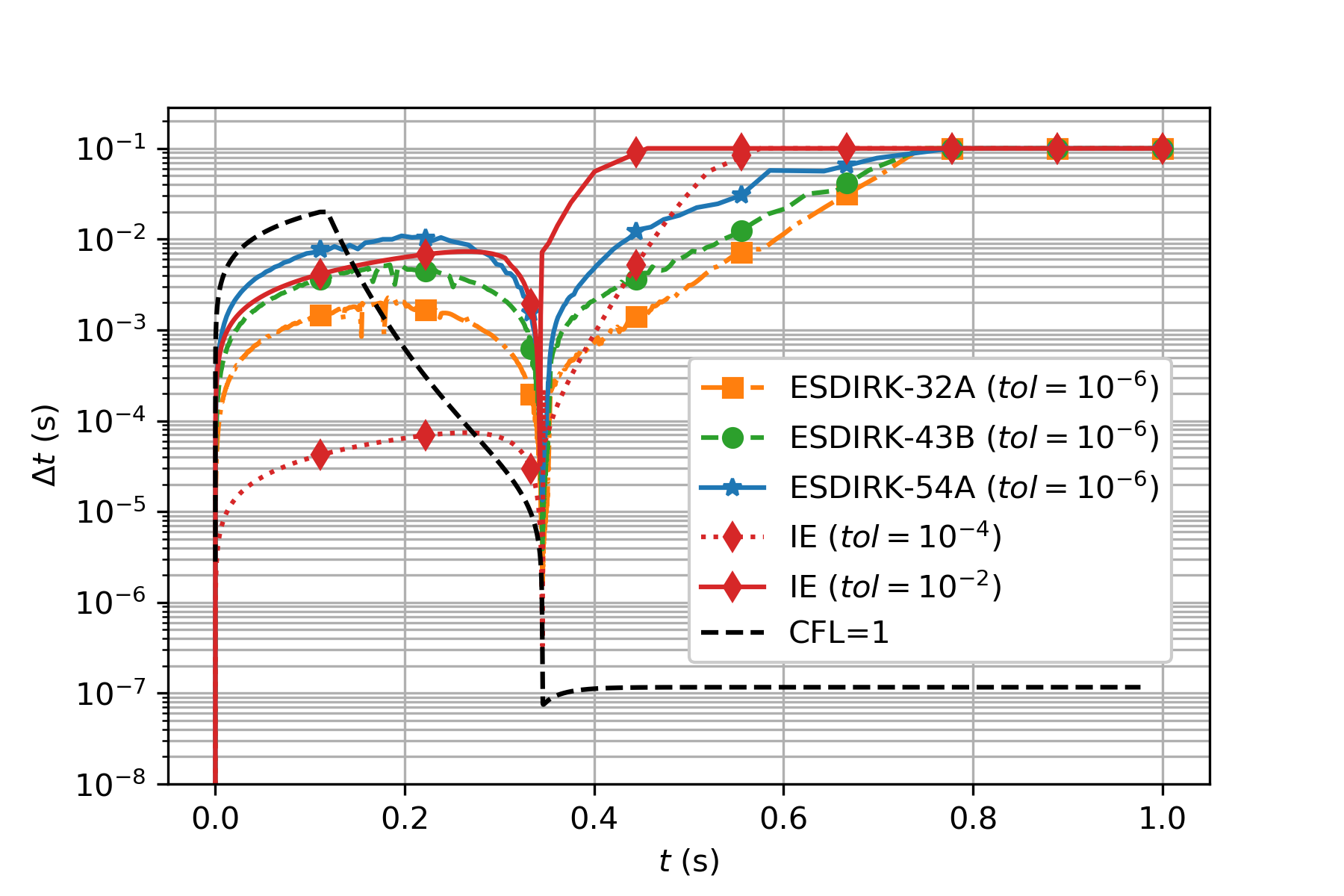}
	\end{minipage}
\end{center}
\begin{center}
	\begin{minipage}{0.45\textwidth}
		\captionof{figure}{Surface temperature evolution during laser ignition}
		\label{fig:allumage:Ts_esdirk}
	\end{minipage}
	~
	\begin{minipage}{0.05\textwidth}
		~
	\end{minipage}
	~
	\begin{minipage}{0.45\textwidth}
		\captionof{figure}{Time steps used by each method}
		\label{fig:allumage:dt_esdirk}
	\end{minipage}
\end{center}

We mainly compare the evolution of $T_\surface$, the computational time, and the ignition time $t_{ign}$. \modif{The latter is defined as the time at which the surface temperature first exceeds 1000 K and is obtained by high-order interpolation}.
Although not shown here, the curves of $T_\surface$ for each method are very similar, except for IE simulations with large tolerances that deviate slightly during the inert heating and ignite a few milliseconds earlier.
Figure \ref{fig:allumage:dt_esdirk} shows the evolution of the time step for some of the simulations. We observe that, for ESDIRK embedded methods, increasing the order of the method allows for larger time steps to be used throughout the integration while maintaining the same accuracy.
For example, the fifth-order method ESDIRK-54A is capable of taking steps 5 times larger in average than the third-order method ESDIRK-32A.
Finally, it is clear that IE needs many more steps to achieve a similar result as the ESDIRK methods. 

\begin{figure}[htpb]
\centering
\begin{subfigure}[t]{0.45\textwidth}
\includegraphics[width=\textwidth]{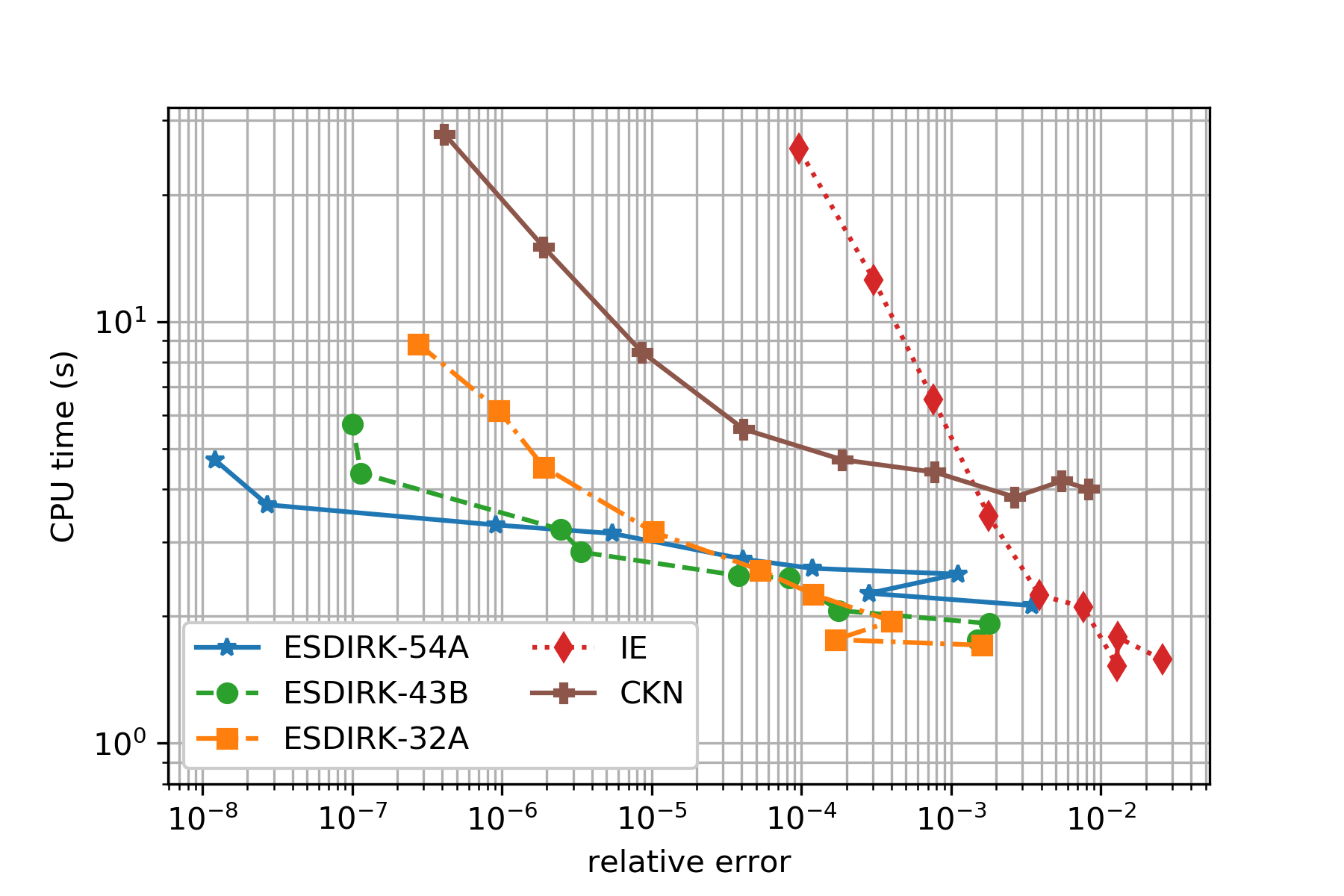}
\caption{coarse mesh (400 cells)}
\label{fig:ignition:workprecision:coarse}
\end{subfigure}
~
\begin{subfigure}[t]{0.45\textwidth}
\includegraphics[width=\textwidth]{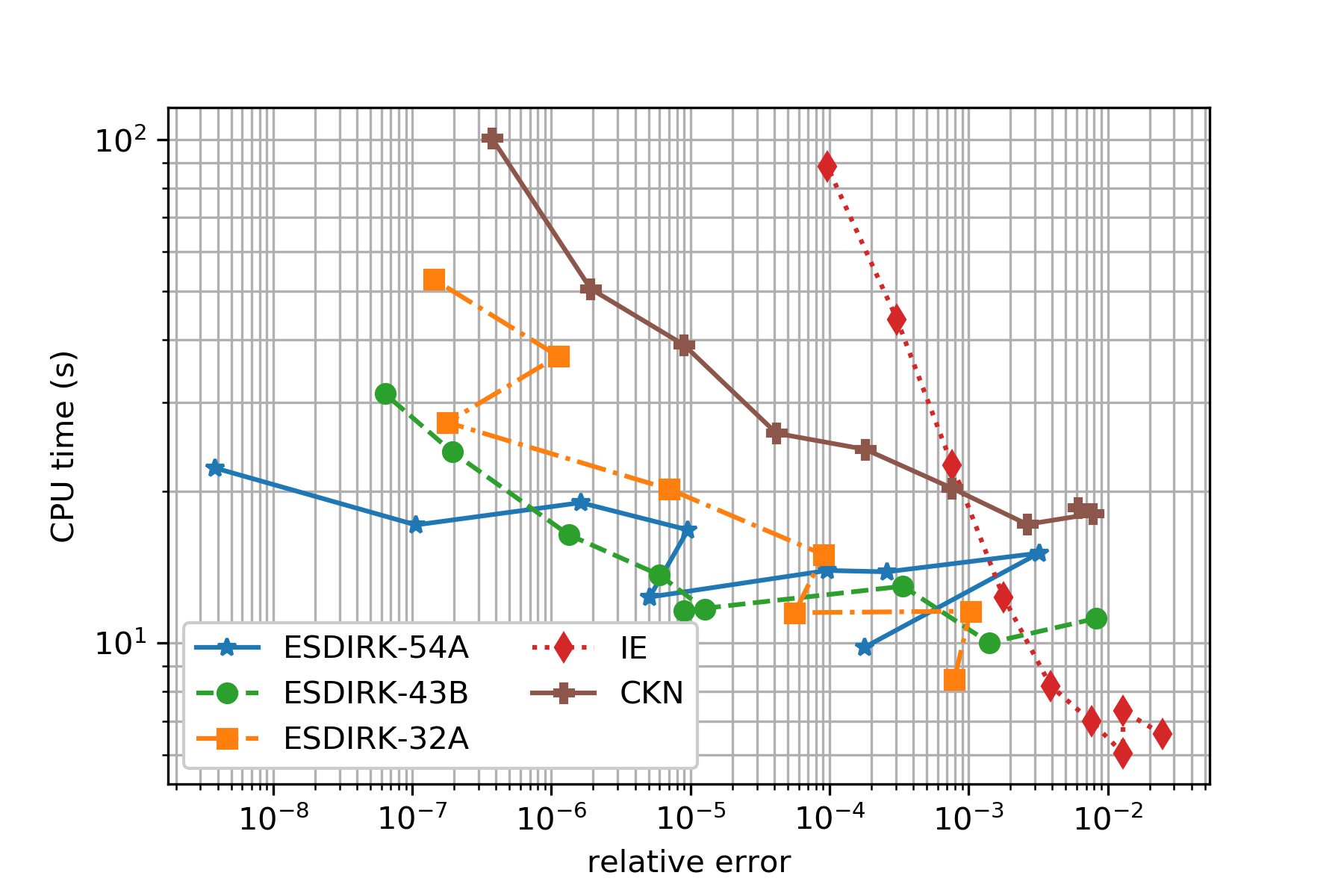}
\caption{fine mesh (2049 cells)}
\label{fig:ignition:workprecision:fine}
\end{subfigure}
\caption{Work-precision diagram for the determination of $t_{ign}$}
\end{figure}

To more quantitatively assess the efficiency of each method in precisely determining the ignition time, multiple simulations were run with each scheme. For the ESDIRK methods, the relative integration error tolerance was varied from $10^{-1}$ to $10^{-6}$. For IE and CKN, the relative solution variation allowed between successive time steps was varied from $10^{-1}$ to $10^{-4}$, without any CFL-limitation. For each simulation the value of $t_{ign}$ is evaluated and a relative error on this value can be inferred by comparing it to the ignition time obtained with ESDIRK-54A and the tightest tolerance.
Figure \ref{fig:ignition:workprecision:coarse} shows the computational time required by each method to achieve a given level of relative error on the initial mesh.
Figure \ref{fig:ignition:workprecision:fine} shows the computational time required by each method on a refined mesh with 734 cells in the solid phase and 1315 cells in the gas phase.
Overall in both cases, if a relatively large error of 1\% on the ignition time is deemed sufficient, IE is a relevant choice. If greater precision is however required, ESDIRK schemes with adaptive time stepping as described in Section \ref{section:RK:adaptation} are much more efficient.

The black dashed curve in Figure \ref{fig:allumage:dt_esdirk} shows the time step evolution corresponding to $\mathrm{CFL}=1$ for the coarse mesh. We clearly see that the CFL-limitation is irrelevant during the inert heating phase, as the mass flow rate is very low. If, after the inert heating, the CFL is to be limited to moderate values (1 to 100) as is usually done, many more time steps would be performed. ESDIRK methods are able to give very accurate results without any such limitation.  We have determined that, if a simulation was to be performed with a maximum CFL of 10 with ESDIRK-54A and $rtol=10^{-6}$ on the physical time interval $[0,0.35]$ s, i.e. only up to $T_\surface \approx 820$ K before the ignition, the simulation would take 4 times more steps than required without CFL-limitation. If the final physical time is increased to $0.4$ s to include most of the transient, this ratio would be 160. Knowing that the simulation without CFL-limitation already achieves an error smaller than $10^{-6}$ on the ignition time and on the rest of the evolution, this clearly shows that a CFL-constraint is not a good choice in terms of computational efficiency.
 
Such an ignition simulation is much more time-consuming when no adaptive time stepping is used.
Indeed, the fixed time step used must be such that the fastest phases, i.e. initial heating and transition to ignition, are sufficiently well resolved.
For instance, a simulation with ESDIRK-54A and $\Delta t=10^{-4}$ s took approximately 400 s to complete.
This is 20 times more than an adaptive simulation with $rtol=10^{-6}$, which also produces a more accurate result.

The advantage of embedded methods is that only a relative integration tolerance needs to be specified. No tuning of a CFL-criterion, maximum relative variation or fixed time step is required, hence such methods speed up the engineer task of simulating different scenarios, while still ensuring a controlled error. One interesting observation we made is that the time step values taken by ESDIRK methods were almost identical with both meshes. In all our testing, no correlation was found between  the time step evolution of the embedded methods and the mesh refinement. This would not be the case if a CFL-criterion was used.

\section{Investigation of limit cycles}

\begin{figure}[hbt!]
\centering
\begin{subfigure}[b]{0.45\textwidth}
\centering
\includegraphics[width=\textwidth]{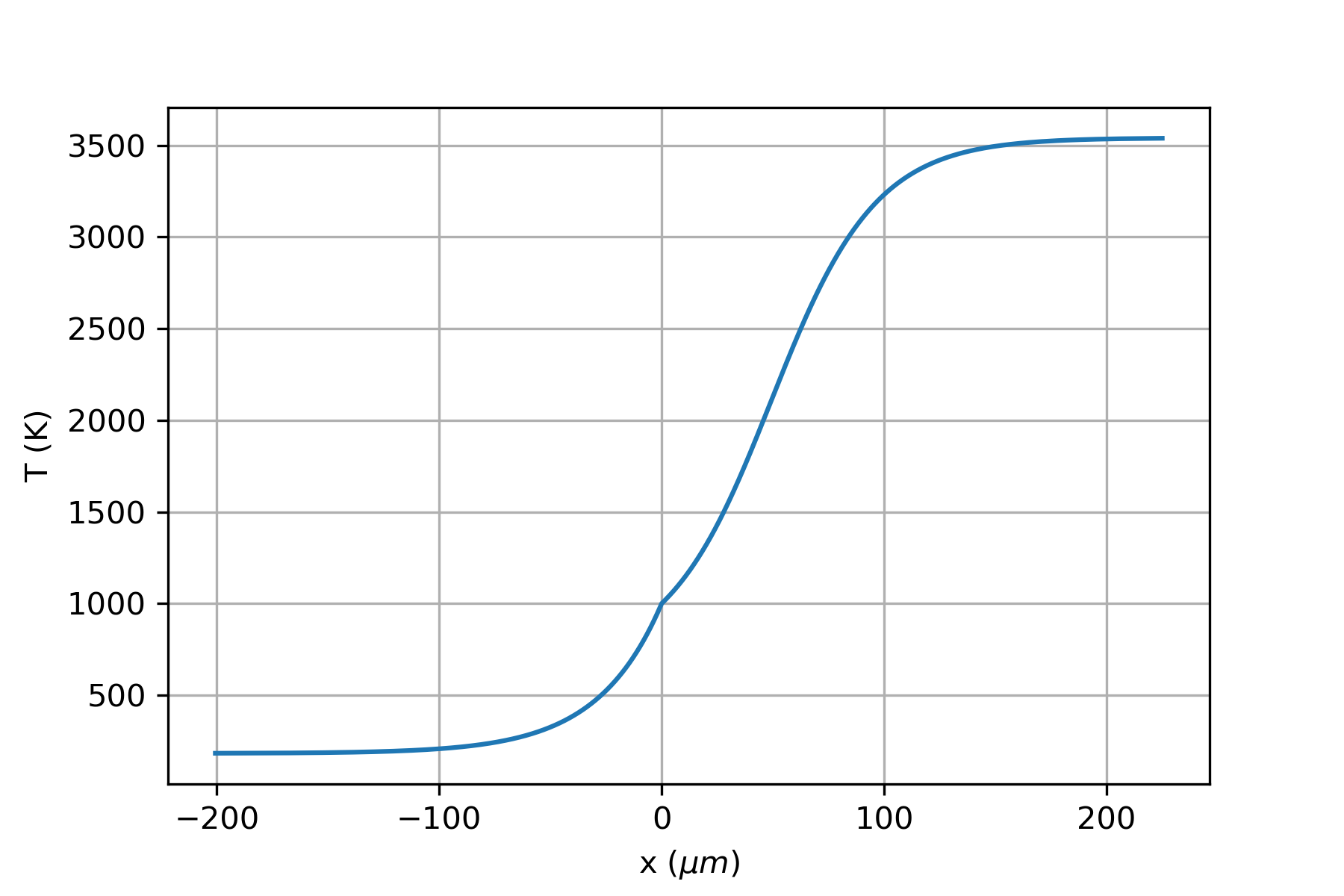}
\caption{Steady-state temperature profile}
\label{fig:optim:steadysol}
\end{subfigure}
~
\begin{subfigure}[b]{0.45\textwidth}
 \centering
\includegraphics[width=\textwidth]{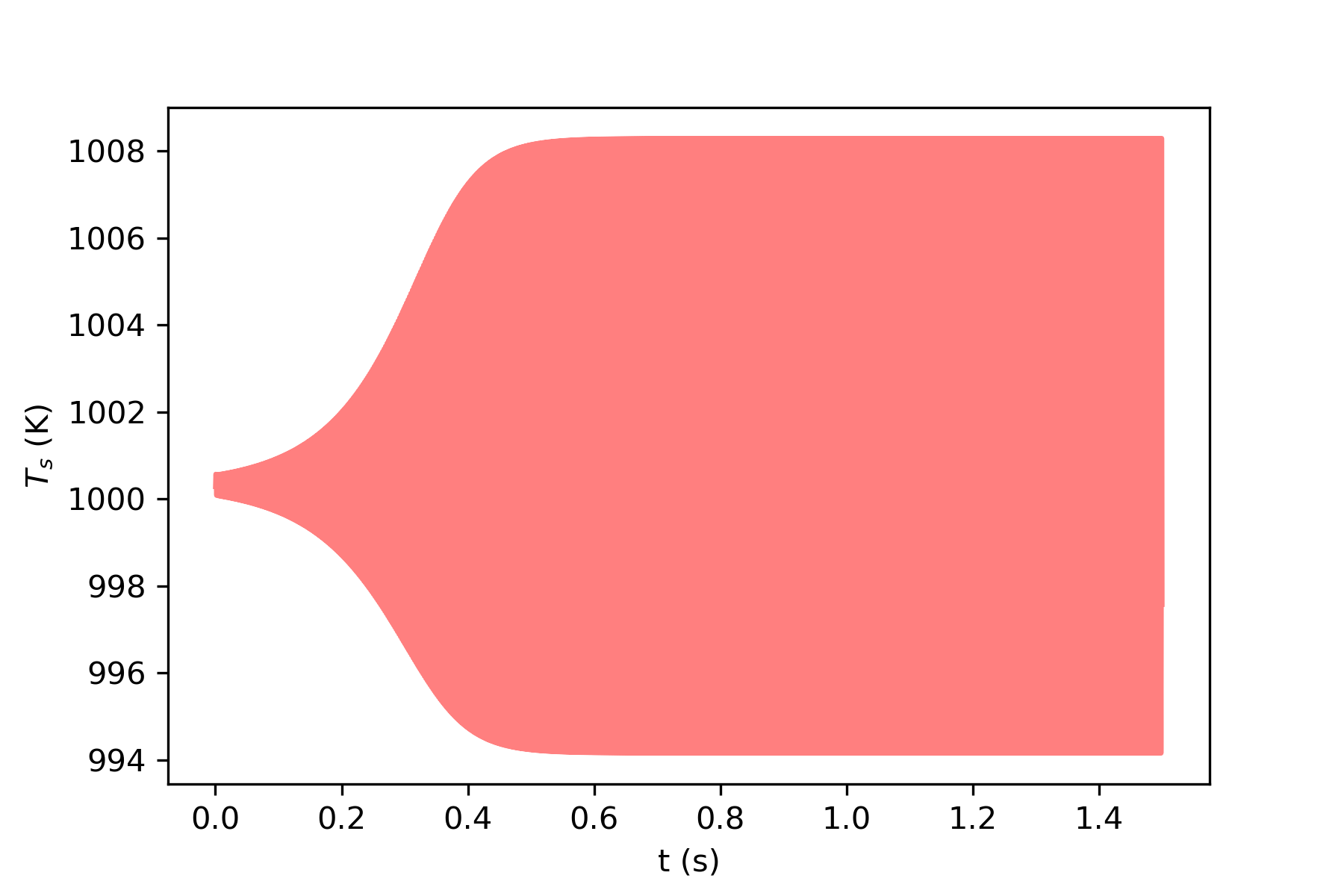}
\caption{Limit cycle}
\label{fig:limit_cycle_reference}
\end{subfigure}

\begin{subfigure}[b]{0.6\textwidth}
\centering
\includegraphics[width=.8\textwidth]{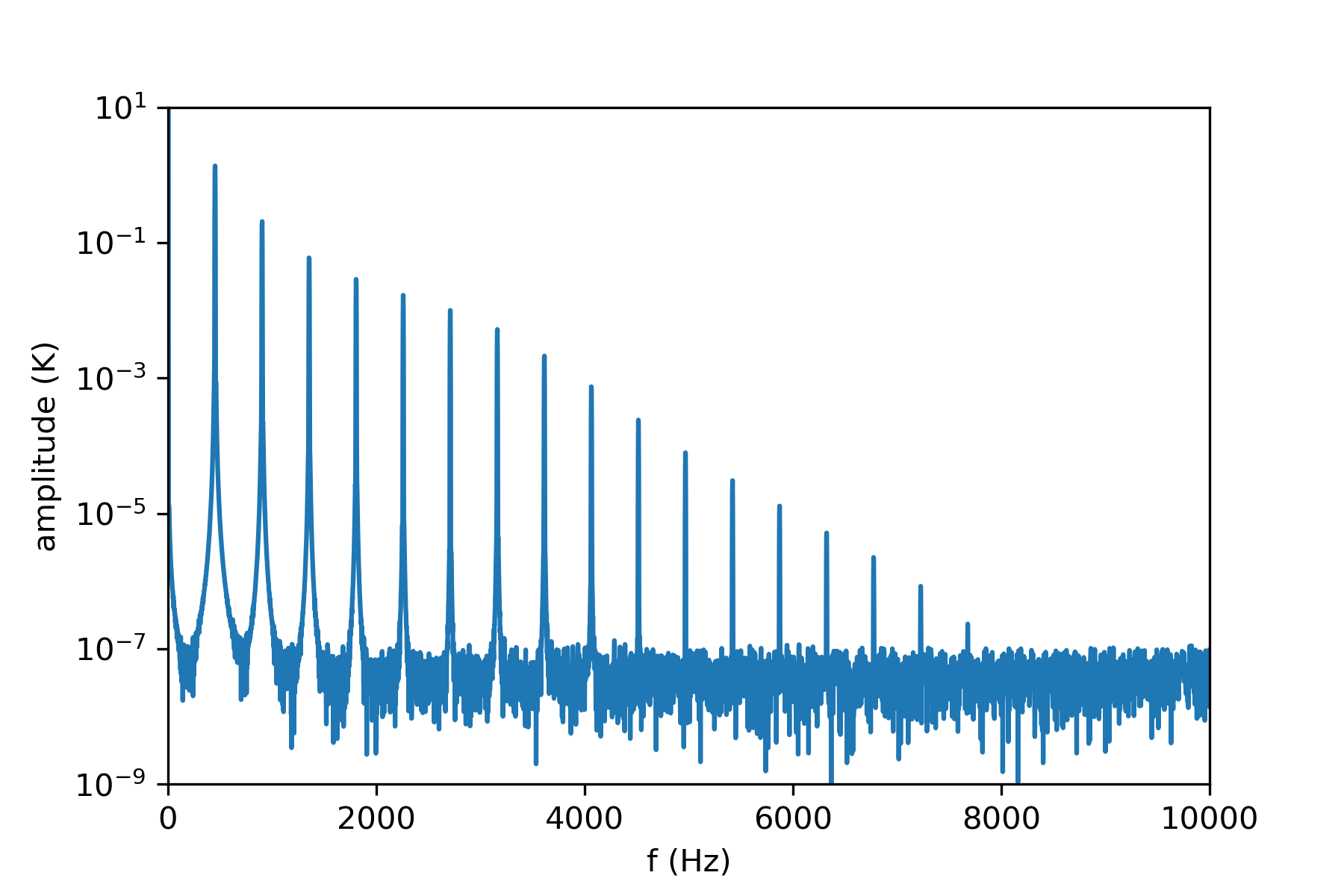}
\caption{Discrete Fourier Transform of the limit cycle}
\label{fig:limitcycle:FFT}
\end{subfigure}
\caption{Main features of the studied configuration}
\end{figure}

The effort made in terms of time integration strategy can be used to accurately study the nonlinear behaviour of the propellant combustion, and in particular potential departures from an unstable steady-state travelling wave solution.
Using the procedure detailed in \ref{appendix:explications_ligneRK}, we have generated a variant of the simplified combustion model used in Section \ref{section:allumage:simplifie}, whose steady state solution is linearly unstable under perturbations but stabilises onto a limit cycle via nonlinear effects.
Using its steady-state temperature profile plotted in Figure \ref{fig:optim:steadysol} as initial solution and applying a small pressure perturbation, the transient destabilisation is computed with ESDIRK54-A and $rtol=10^{-6}$.
The surface temperature history is plotted in Figure \ref{fig:limit_cycle_reference}, clearly exhibiting a limit cycle, the discrete Fourier Transform of which is presented in Figure \ref{fig:limitcycle:FFT}. The limit cycle can be decomposed as a sum of sinusoidal harmonic oscillations, with a fundamental frequency of 452 Hz, close to the propellant natural frequencies defined in \cite{reponse_culick_1968} and \cite{theseShihab}, at 518 Hz and 348 Hz respectively.

As reported in the literature for other applications \cite{cycle_limite_tube_rijke_high_order}, high-order methods are often needed to be able to numerically reproduce such a dynamical and nonlinear behaviour. It is therefore instructive to compare the methods from Table \ref{table:rk_integrators}
to see how they affect the unsteady result. Namely, it is expected that the integration methods will, depending on their order, stability properties and the time step used, dampen the oscillating nature of the system and potentially cause a non-physical stabilisation of the solution.

First, constant time steps simulations are performed. Comparisons are made based on the ability to reproduce the initial amplification of the oscillations, the fundamental frequency of the limit cycle and its amplitude.
Second, simulations with time adaptation are presented.

\subsection{Method}

As described in \ref{appendix:explications_ligneRK}, the simulations are performed on a non-uniform mesh  based on the temperature profile of the steady-state solution with  55 cells for the solid phase and 146 for the gas phase. It has been verified that additional refinement would not affect the solution dynamics.

\label{section:limitcycle:comparisonprocess}
We focus on several aspects. First we analyse visually the envelope curve of the surface temperature time history. For each simulation, this curve is
constructed  as the junction of the successive maxima of $T_\surface$. Unfortunately, for large values of the time step $\Delta t$, the sampling frequency $1/\Delta t$ can be too large compared to the fundamental frequency of the limit cycle, causing artefacts to appear in the form of an oscillation of the envelope. This can be improved by using a cubic interpolant of $T_\surface$ to determine the successive maxima with greater precision, however envelope oscillations are still present, for example in Figure \ref{fig:limitcycle:dt_constant:reveil_esdirk54}.

Second, we analyse the growth of the linear instability near $t=0$. For each simulation, the best exponential fit for the envelope of the evolution of $T_\surface$ is determined, i.e. the curve joining the successive maxima of $T_\surface$, obtained from a cubic interpolation of the temporal evolution of $T_\surface$. Such a fit is of the form
$T_{fit}(t) - T_{\surface}(0) = A \exp(\amplifactor t)$, with $\amplifactor$ the fitted amplification factor.

Third, the established limit cycle is considered, on the time window $1\leq t \leq 1.5$ s.
A discrete Fourier transform of the surface temperature signal is computed via an FFT algorithm, as shown for the reference simulation in Figure \ref{fig:limitcycle:FFT}. Interpolation of the solution on a uniform time grid is performed if the simulation was not conducted with a constant time step. This FFT helps determine the approximate frequencies of the different harmonics with a precision of approximately 1-10 Hz. For each of these peaks in the spectrum, the peak frequency is then precisely computed by maximising the correlation between $T_\surface(t)$ and $\exp(2 \mathrm{i} \pi f t)$, from which we can also determine the precise amplitude of the corresponding peak. These values offer a trustworthy and precise indication of how well the limit cycle is captured.

Finally, work-precision diagrams are given which represent the computational time required to achieve a specific level of relative error on the quantitative results, i.e. fitted amplification factor, fundamental frequency and amplitude. Relative errors are computed relative to the values obtained with the most refined solution.

\subsection{Analysis of schemes efficiency for a constant time step}
\label{section:limitcycle:constant_dt}

\subsubsection{Envelope of the surface temperature history}
\begin{figure}[!htpb]
\centering
\begin{subfigure}[b]{0.45\textwidth}
\centering
\includegraphics[width=\textwidth]{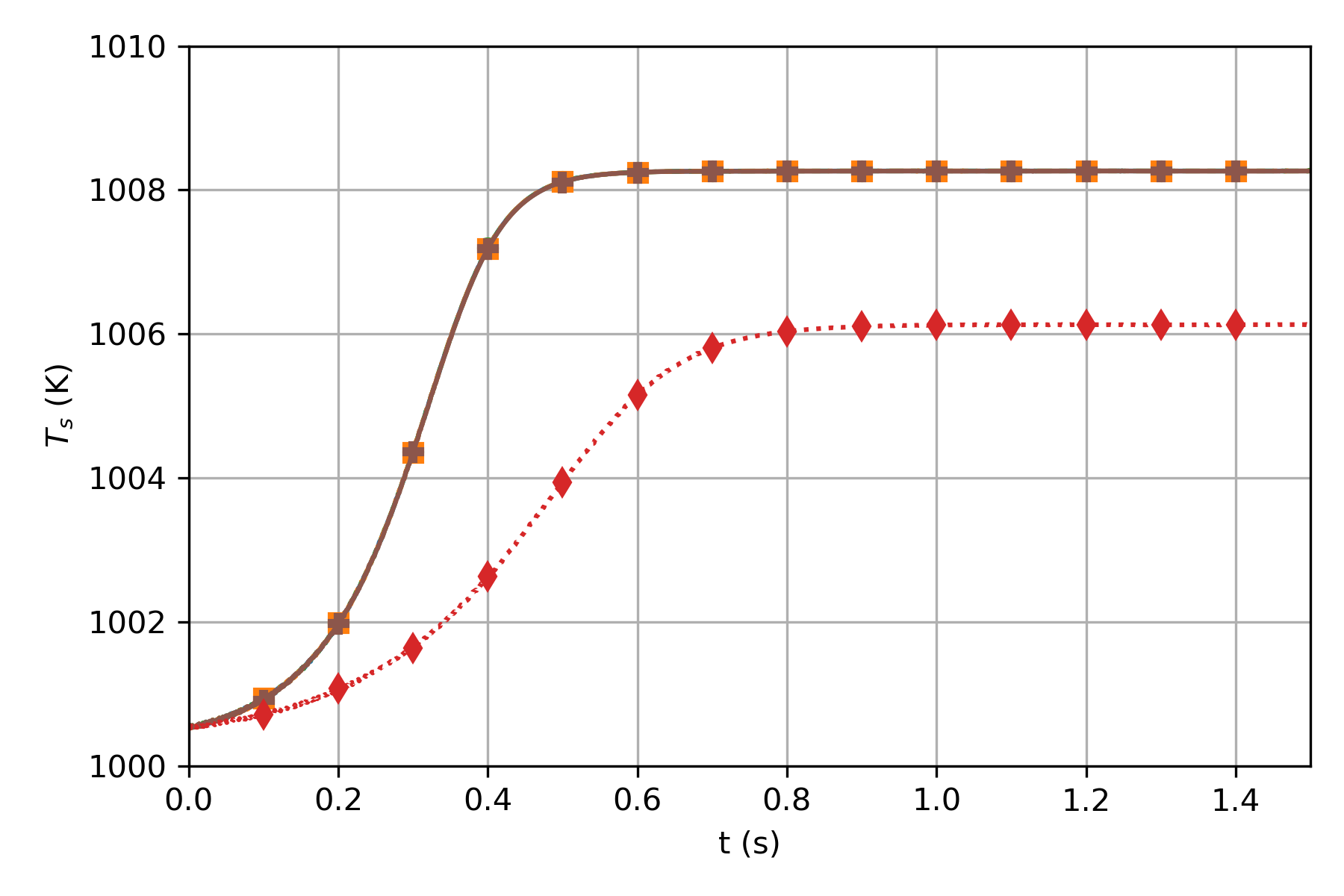}
\caption{$\Delta t=1\times10^{-6}$ s}
\end{subfigure}
~
\begin{subfigure}[b]{0.45\textwidth}
\centering
\includegraphics[width=\textwidth]{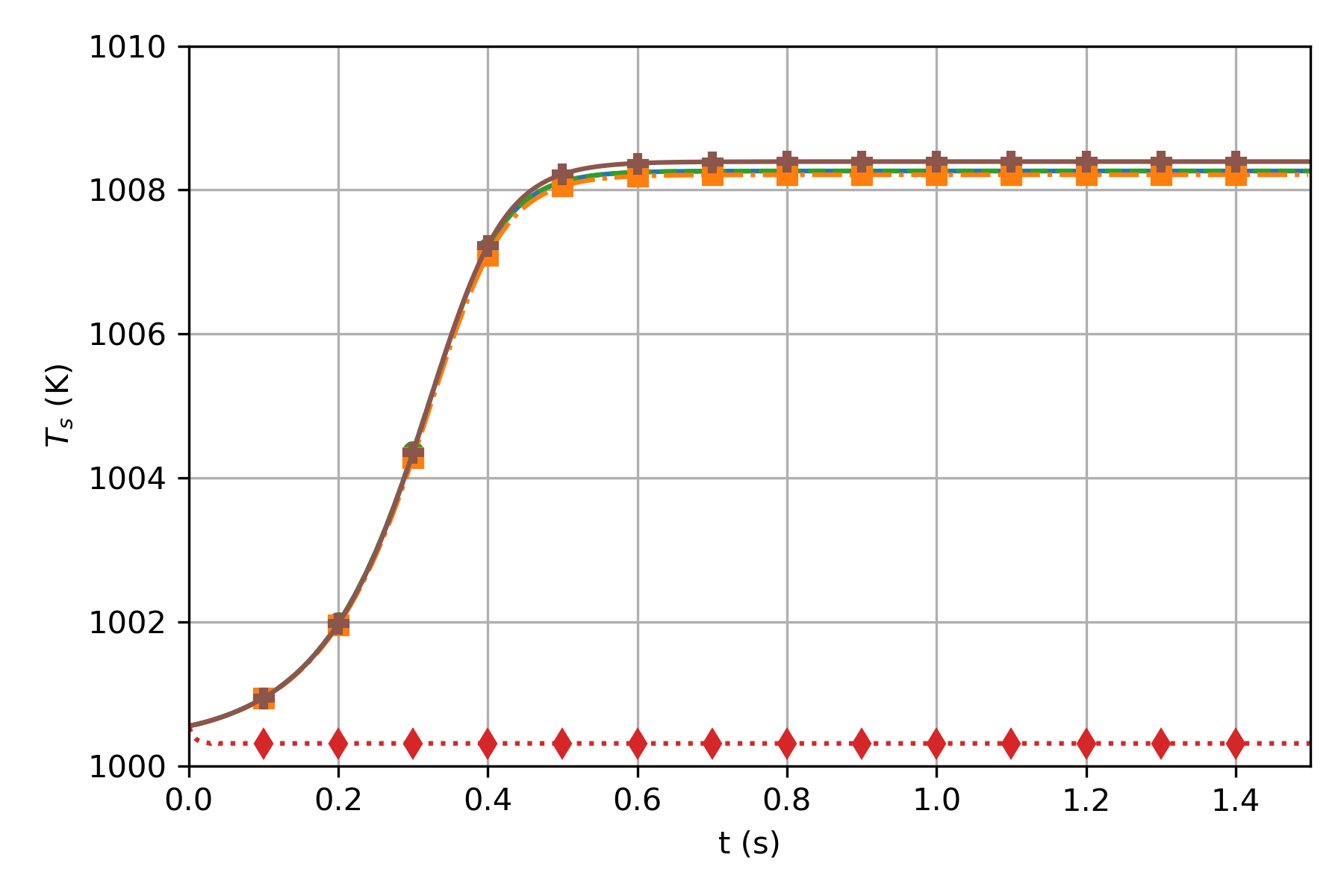}
\caption{$\Delta t=3.79\times10^{-5}$ s}
\end{subfigure}

\begin{subfigure}[b]{0.45\textwidth}
\centering
\includegraphics[width=\textwidth]{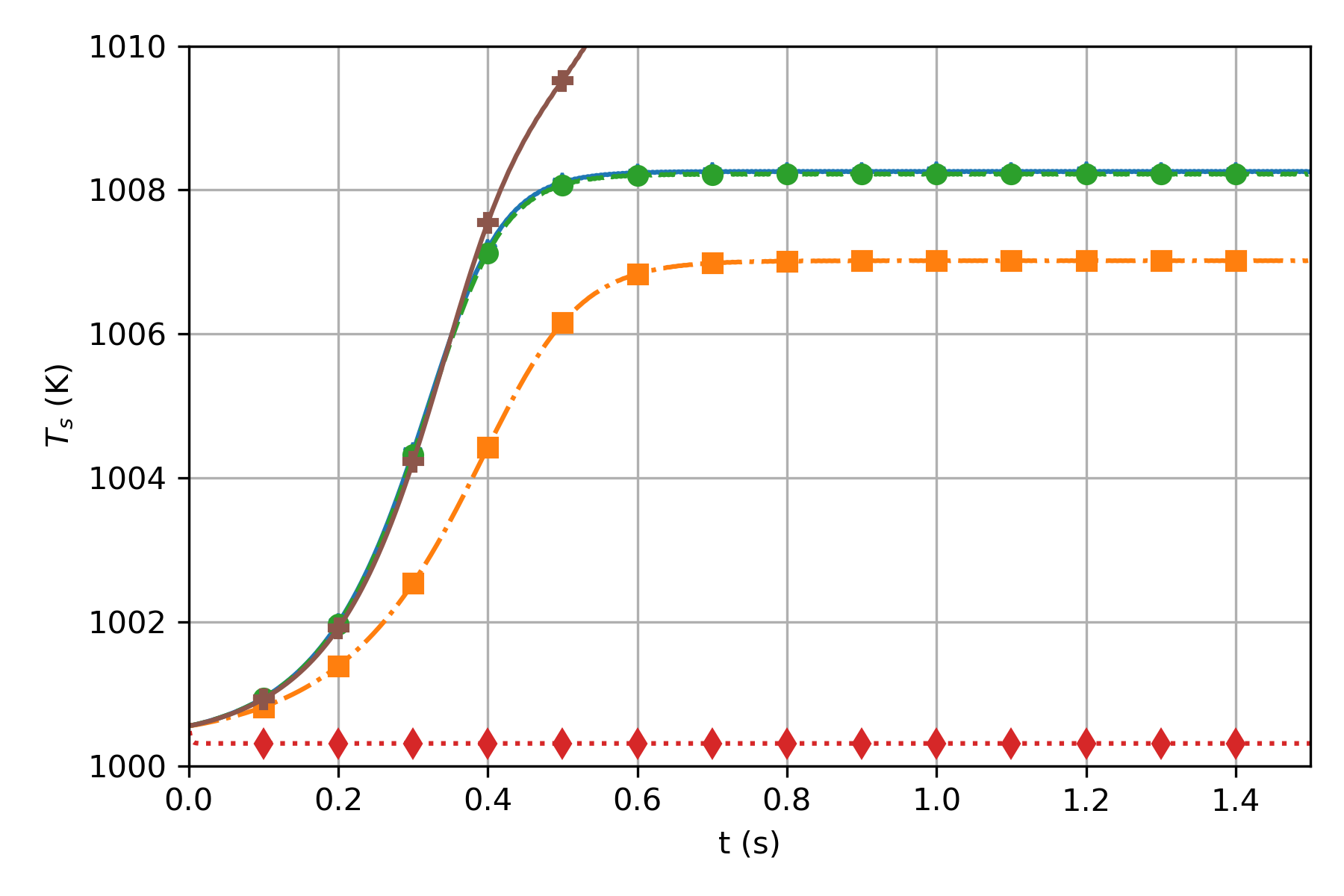}
\caption{$\Delta t=1.13\times10^{-4}$ s}
\end{subfigure}
~
\begin{subfigure}[b]{0.45\textwidth}
\centering
\includegraphics[width=\textwidth]{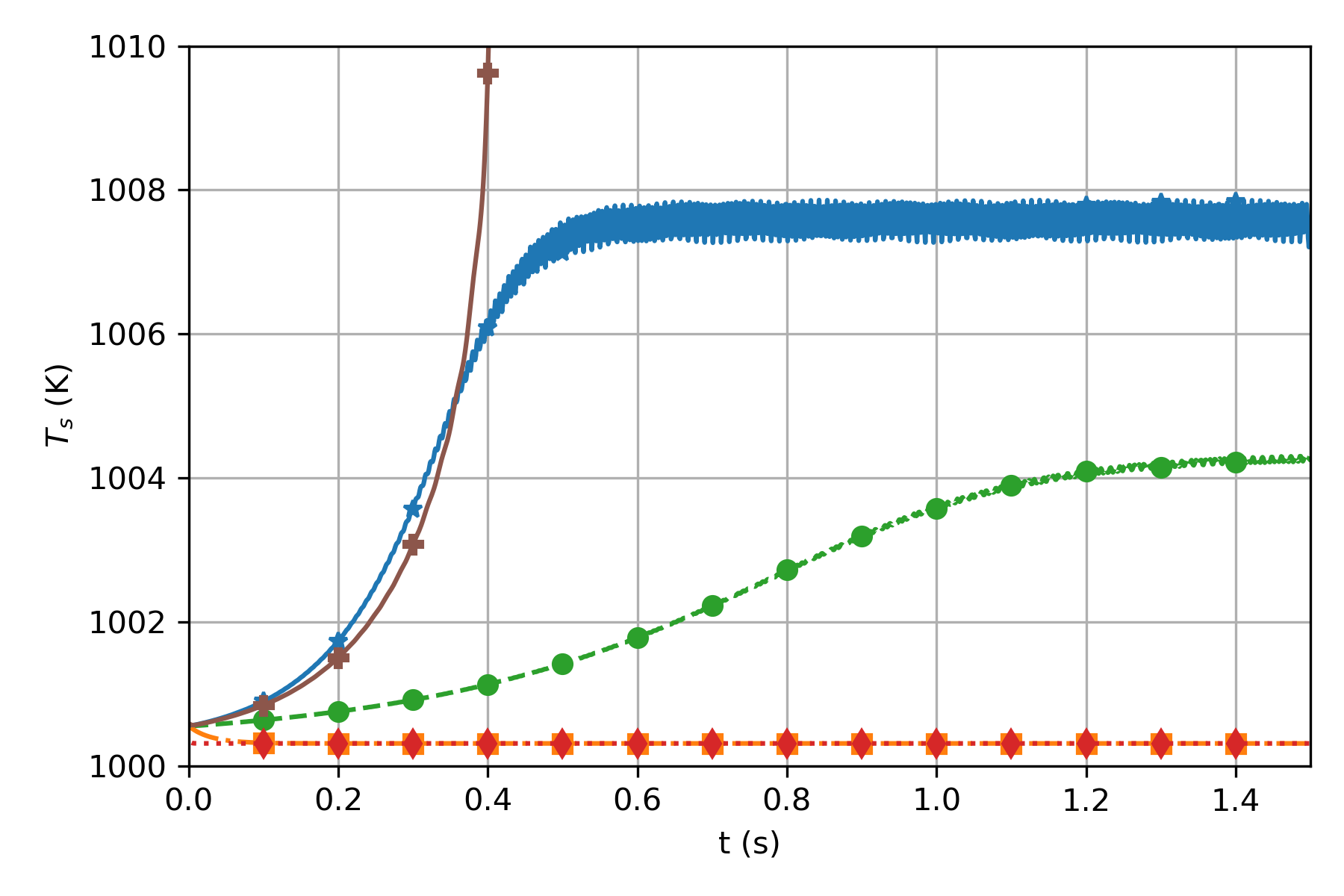}
\caption{$\Delta t=3.36\times10^{-4}$ s}
\end{subfigure}

\begin{subfigure}[b]{0.45\textwidth}
\centering
\includegraphics[width=\textwidth]{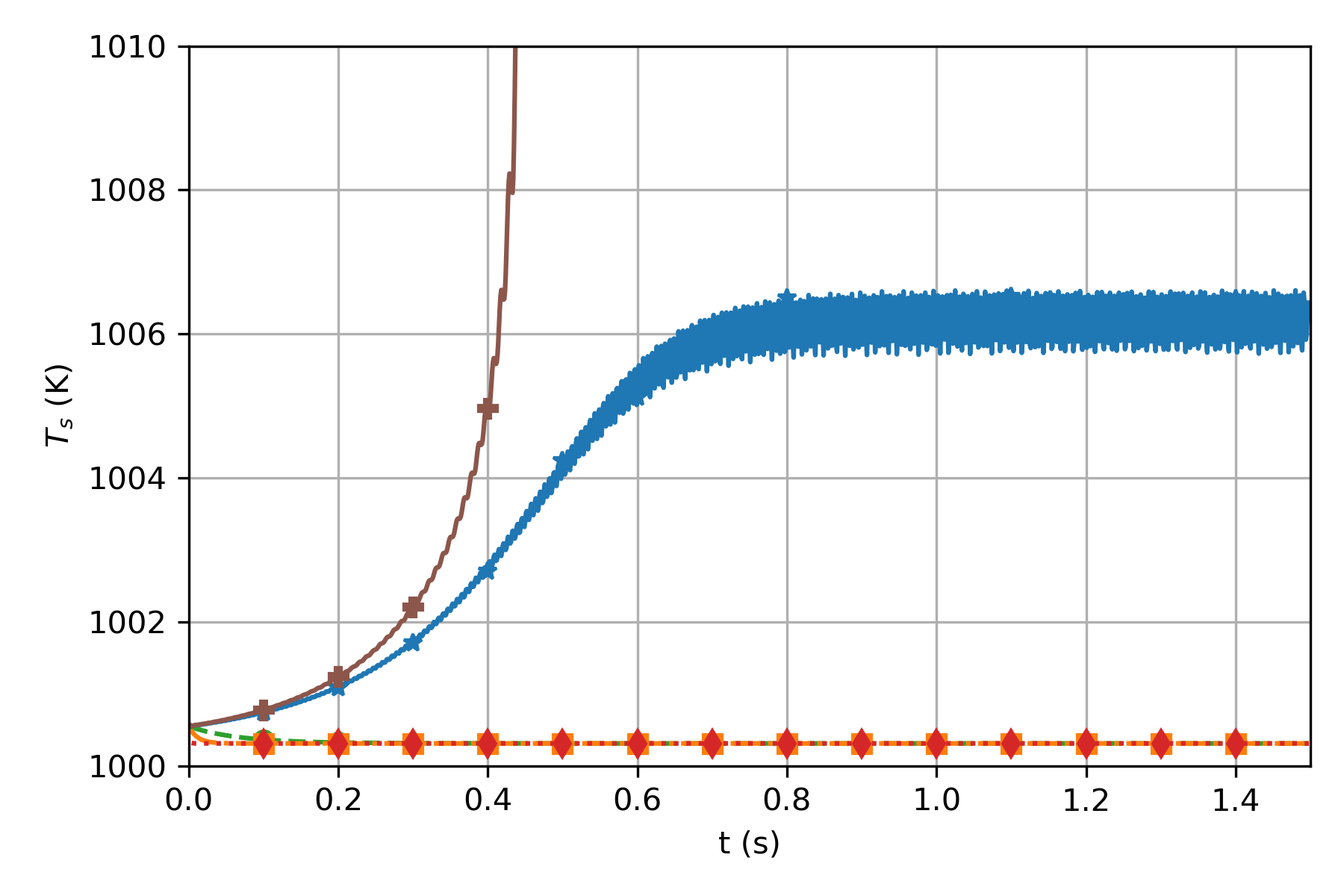}
\caption{$\Delta t=4.83\times10^{-4}$ s}
\label{fig:limitcycle:dt_constant:reveil_esdirk54}
\end{subfigure}
~
\begin{subfigure}[b]{0.45\textwidth}
\centering
\includegraphics[width=\textwidth]{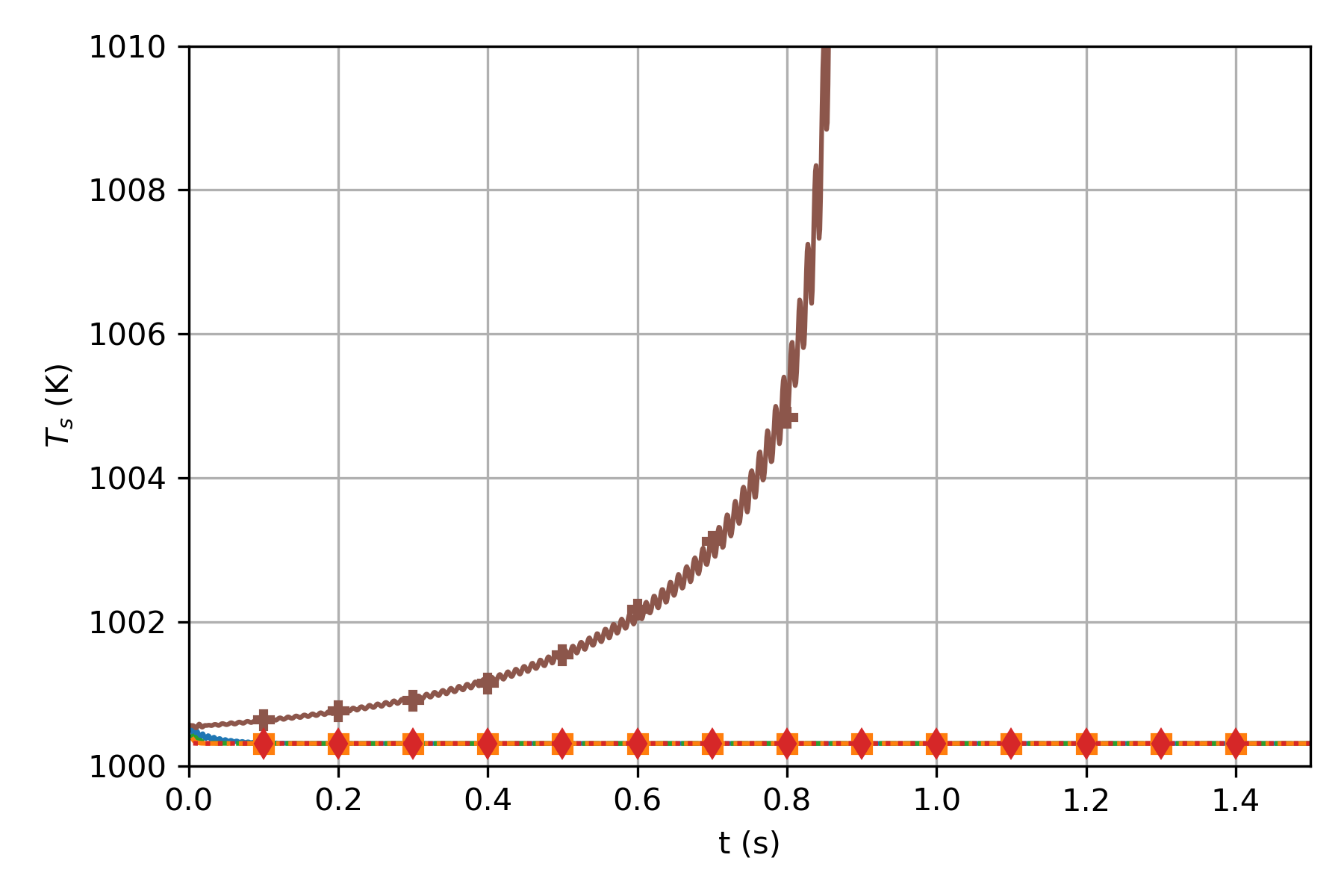}
\caption{$\Delta t=1\times10^{-3}$ s}
\label{fig:limitcycle:enveloppes:dt1em3}
\end{subfigure}

\begin{subfigure}[b]{0.45\textwidth}
\centering
\includegraphics[width=\textwidth]{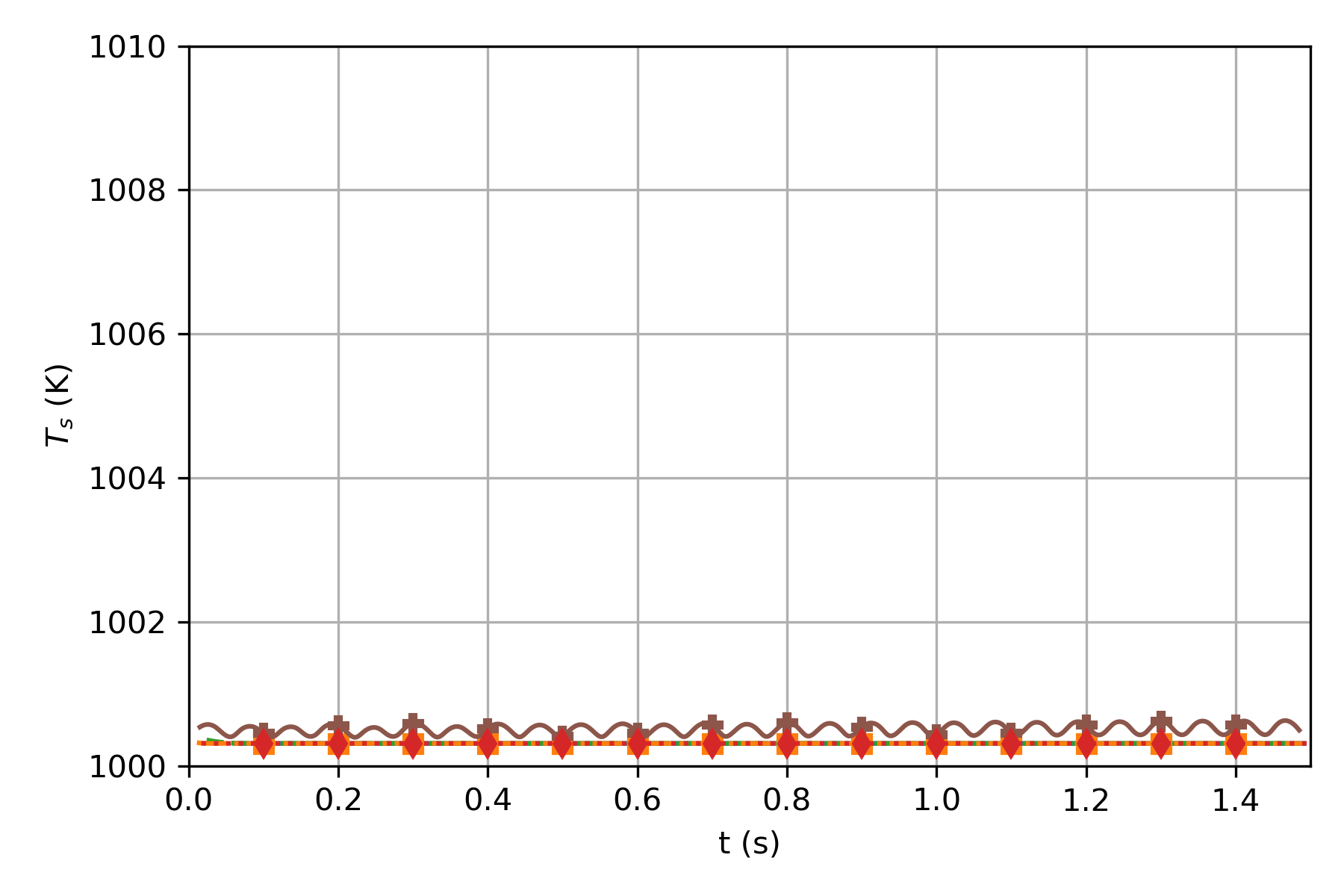}
\caption{$\Delta t=5\times10^{-3}$ s}
\end{subfigure}
~
\begin{subfigure}[b]{0.45\textwidth}
\centering
\includegraphics[width=0.75\textwidth , trim={0cm -10cm 0cm -2cm}, clip]{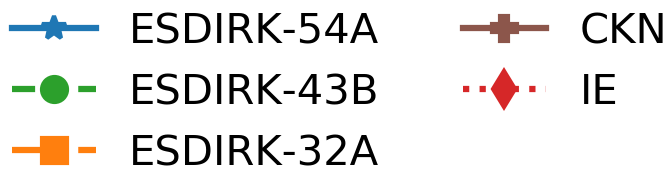}
\end{subfigure}

\caption{Envelopes of the surface temperature histories computed for different time steps and integration methods}
\label{fig:limitcycle:enveloppes}
\end{figure}

We first observe the envelope curves of $T_\surface$ for various time step values in Figure \ref{fig:limitcycle:enveloppes}.
We see that all methods dampen out the oscillations when the time step is too large, except CKN which stabilises at a small oscillating amplitude. However, if we gradually decrease the time step, each method eventually produces a limit cycle. We see that ESDIRK-54A has the best behaviour in terms of reproducing the actual reference limit cycle. It generally seems that the higher the order of the method, the larger the time step can be while still resolving the limit cycle. An interesting behaviour is observed for the CKN method: though it is second-order accurate, it finds a relatively correct initial amplification with larger time steps than required by the fourth- and third-order methods\modif{, as seen in Figures \ref{fig:limitcycle:value_amplification} and \ref{fig:limitcycle:convergence_amplification} where CKN is the first method to obtain a sensible initial amplification of the disturbance for $\Delta t>4\times 10^{-4}$ s, which is also confirmed visually in Figure \ref{fig:limitcycle:enveloppes}}.
\modif{On the opposite, all L-stable methods dampen the oscillations when the time step is too large, and none of them diverges.
	This can be understood by noting that, for the simple Dahlquist test equation $y'=\lambda y$, $\lambda \in \mathbb{C}$, L-stable schemes possess a bounded ``instability" zone (the complement of the domain of absolute stability in $\mathbb{C}$). In our case, the initial linear instability is linked to positive eigenvalues of our system Jacobian. If the time step is too large, the product of the time step times the eigenvalue can be outside of the instability zone, thus the instability is artificially damped by the time scheme.
CKN has an unbounded instability domain (the whole half-plane $Re(\lambda)>0$), therefore this phenomenon does not occur and the instability can freely develop, even though its growth rate is poorly captured when the time step is too large, as we may observe in Figure \ref{fig:limitcycle:enveloppes:dt1em3}.}

Finally, we see that the first-order IE is not able to correctly reproduce the limit cycle, even with the smallest time step of $10^{-6}$ s. Using such a time step is already prohibitive, therefore lowering it further cannot be considered a viable solution to achieve an accurate result. Convergence results presented hereafter are obtained on the range $\Delta t \in [10^{-5}, 10^{-1}]$ s, where IE never produces an initial amplification and established limit cycle, therefore its results are omitted for the sake of  readability.

\subsubsection{Initial growth}

\begin{figure}[hbt!]
\centering
\begin{subfigure}[b]{0.45\textwidth}
\centering
\includegraphics[width=\textwidth]{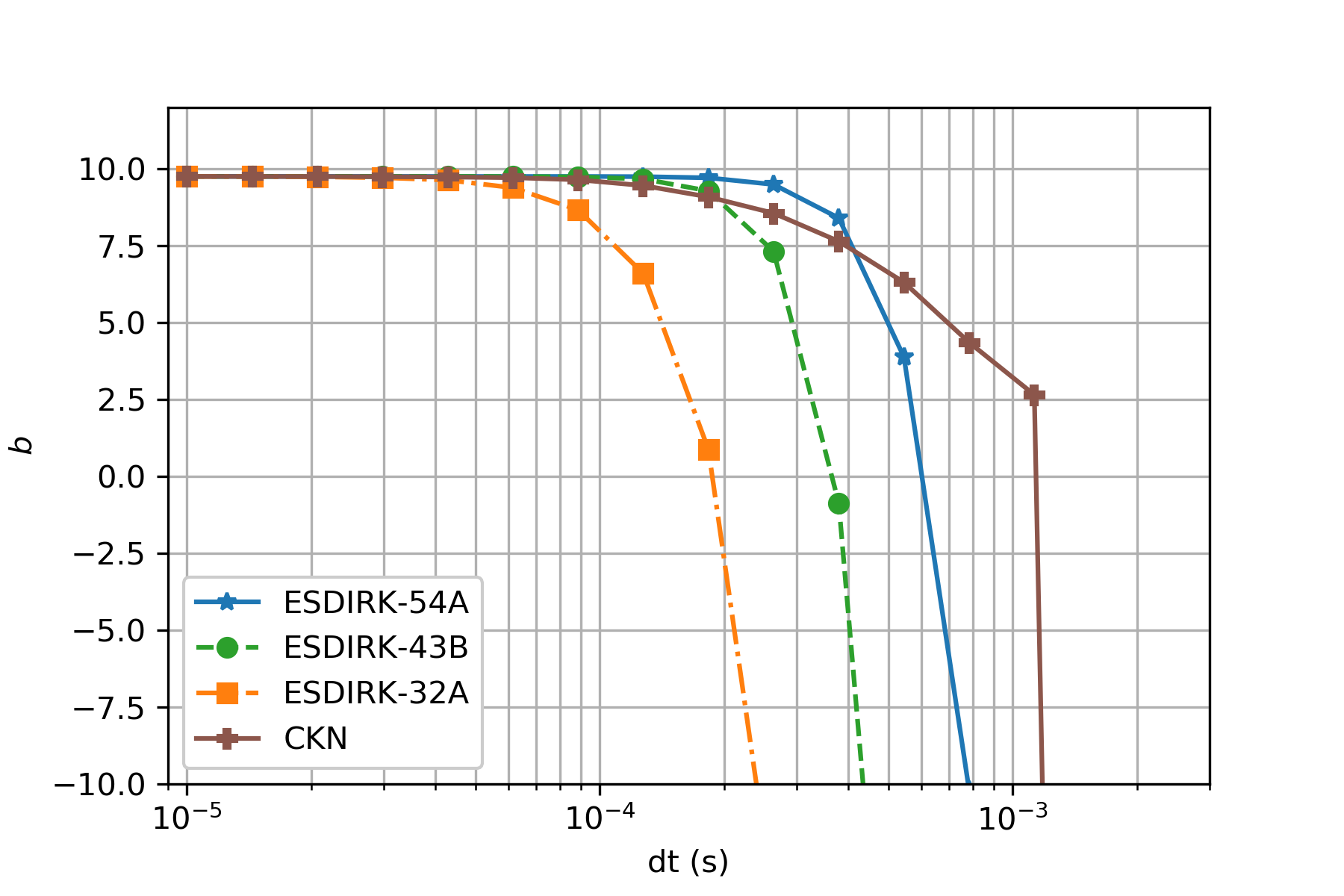}
\caption{Value}
\label{fig:limitcycle:value_amplification}
\end{subfigure}
~
\begin{subfigure}[b]{0.45\textwidth}
\centering
\includegraphics[width=\textwidth]{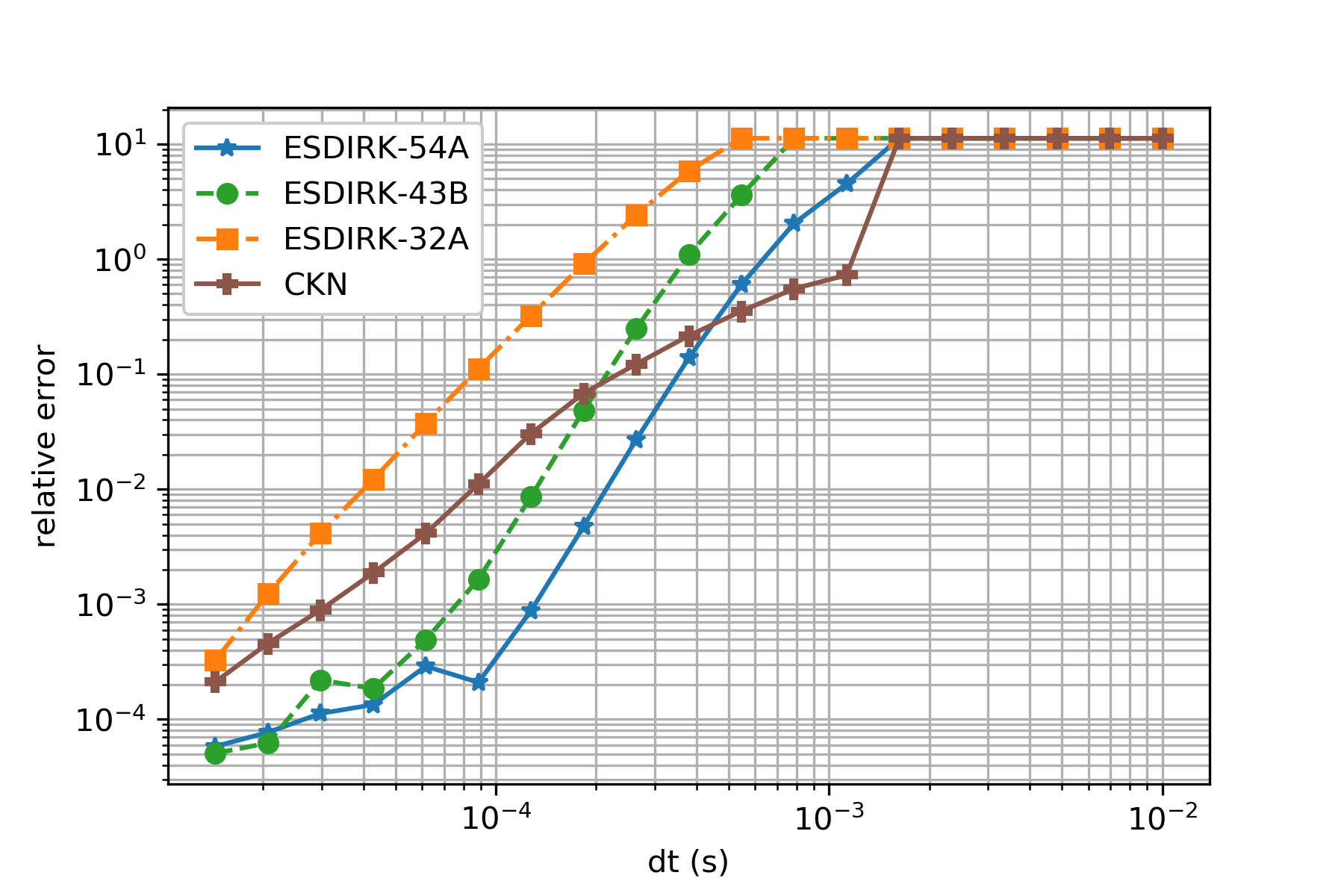}
\caption{Convergence of relative error}
\label{fig:limitcycle:convergence_amplification}
\end{subfigure}
\caption{Fitted amplification factor}
\end{figure}

Figure \ref{fig:limitcycle:value_amplification} shows the evolution of the fitted exponential amplification factor $\amplifactor$ as the time step is lowered. We see that all methods converge to the same value, however ESDIRK-54A and CKN are the first methods that manage to capture a growth (crossing the line $\amplifactor=0$), and also the quickest to converge to the correct value.

Figure \ref{fig:limitcycle:convergence_amplification} shows the relative error of the amplification factor with respect to the reference solution. We see once again the same ranking in terms of ability to find the correct factor. At any time step $\Delta t \leq 5 \times 10^{-4}$ s, ESDIRK-54A yields the best accuracy. Moreover we can observe that each method has an asymptotic convergence region  where the order of convergence is close to the order of the method.
In particular CKN, which initially performs well for moderate time steps, is quickly overtaken by the other methods that possess a higher convergence rate.
However, we notice that both ESDIRK-43B and ESDIRK-54A converge with order 4 for the value of the amplification factor $k$. This can be simply explained by the fact that $k$ is determined via cubic interpolation of the curve of $T_s(t)$ (see Section \ref{section:limitcycle:comparisonprocess}), which introduces an error proportional to $\Delta t^4$ in the determination of $k$, which dominates the true error of the scheme.

For low time step values, the flattening of the convergence curves can be simply explained. The amplification factor is defined as a coefficient from an exponential fit, however this fit is only an approximation, as the nonlinear behaviour will let the unsteady evolution slightly deviate from the theoretical exponential initial growth. Also, the fit is based only on the successive maxima, not on the complete oscillating curve, which induces additional errors, e.g. imprecision in the abscissas of the maxima. Therefore there is a point at which the precision achieved with an exponential approximation cannot be improved further.

\subsubsection{Limit cycle}

\begin{figure}[hbt!]
\centering
\begin{subfigure}[b]{0.45\textwidth}
\centering
\includegraphics[width=\textwidth]{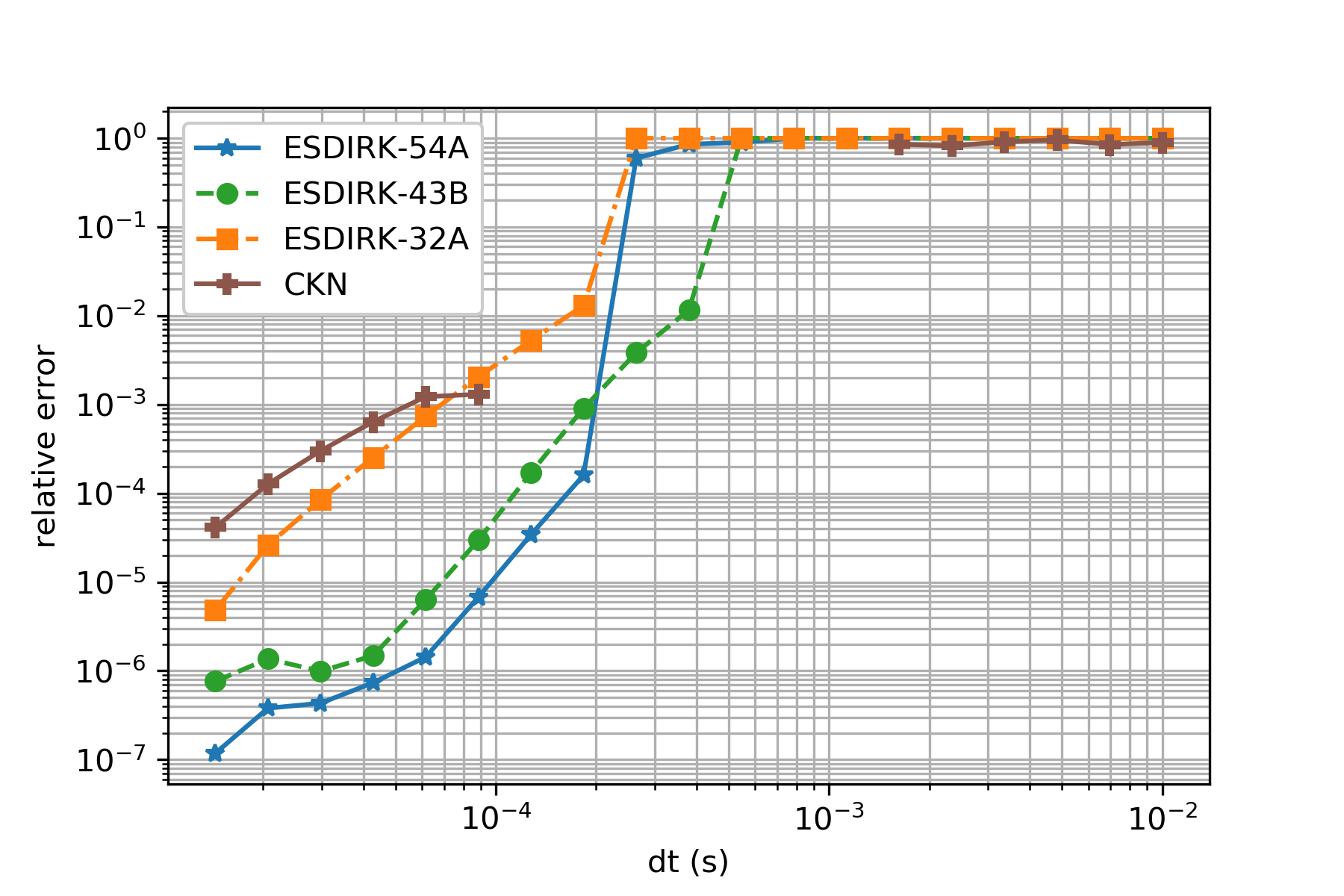}
\caption{Fundamental frequency}
\label{fig:limitcycle:convergence_freq_harm0}
\end{subfigure}
~
\begin{subfigure}[b]{0.45\textwidth}
\centering
\includegraphics[width=\textwidth]{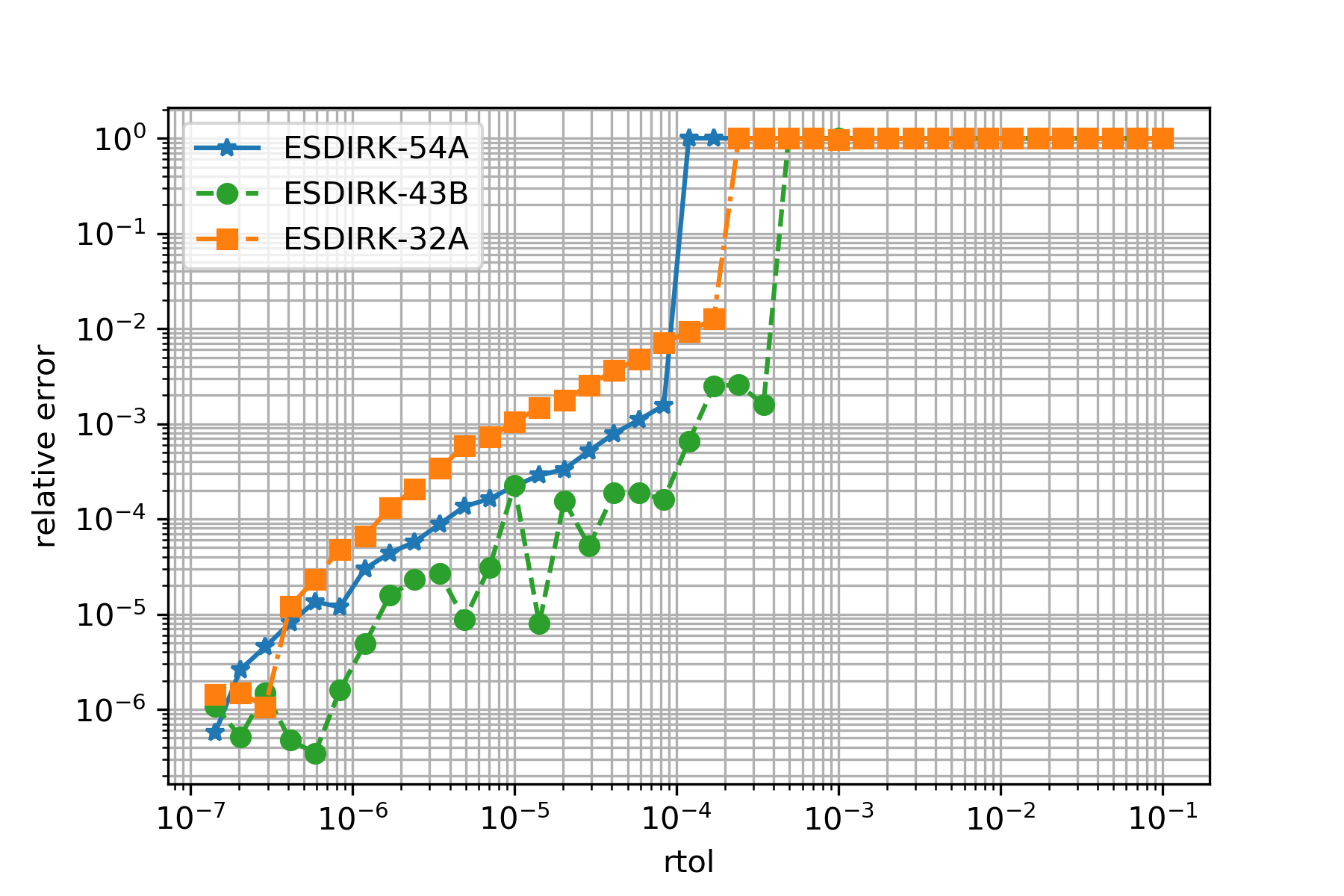}
\caption{Amplitude of the fundamental}
\label{fig:limitcycle:convergence_amplitude_harm0}
\end{subfigure}
\caption{Convergence of the limit cycle properties}
\end{figure}

The fundamental frequency of the established limit cycle is $f\approx 452$ Hz. Figure \ref{fig:limitcycle:convergence_freq_harm0} shows how the relative error on this frequency evolves with the time step, and Figure \ref{fig:limitcycle:convergence_amplitude_harm0} shows the convergence of the amplitude of the fundamental frequency.
The brown crossed curve for CKN is interrupted in the intermediate range of time step values, as the solution diverges, thus not allowing for an established limit cycle to be analysed.
We can observe that ESDIRK-43B and ESDIRK-54A yield the most precise solutions at any given time step. In particular, ESDIRK-43B is able to capture a non-zero oscillation amplitude with larger time steps than required by the other methods.
The frequency-finding process is limited in accuracy for the determination of the fundamental frequency, hence the flattening of the convergence curves when the relative error reaches
$10^{-6}$.

\subsubsection{Computational cost}

\begin{figure}[hbt!]
\centering
\begin{subfigure}[t]{0.45\textwidth}
\centering
\includegraphics[width=\textwidth]{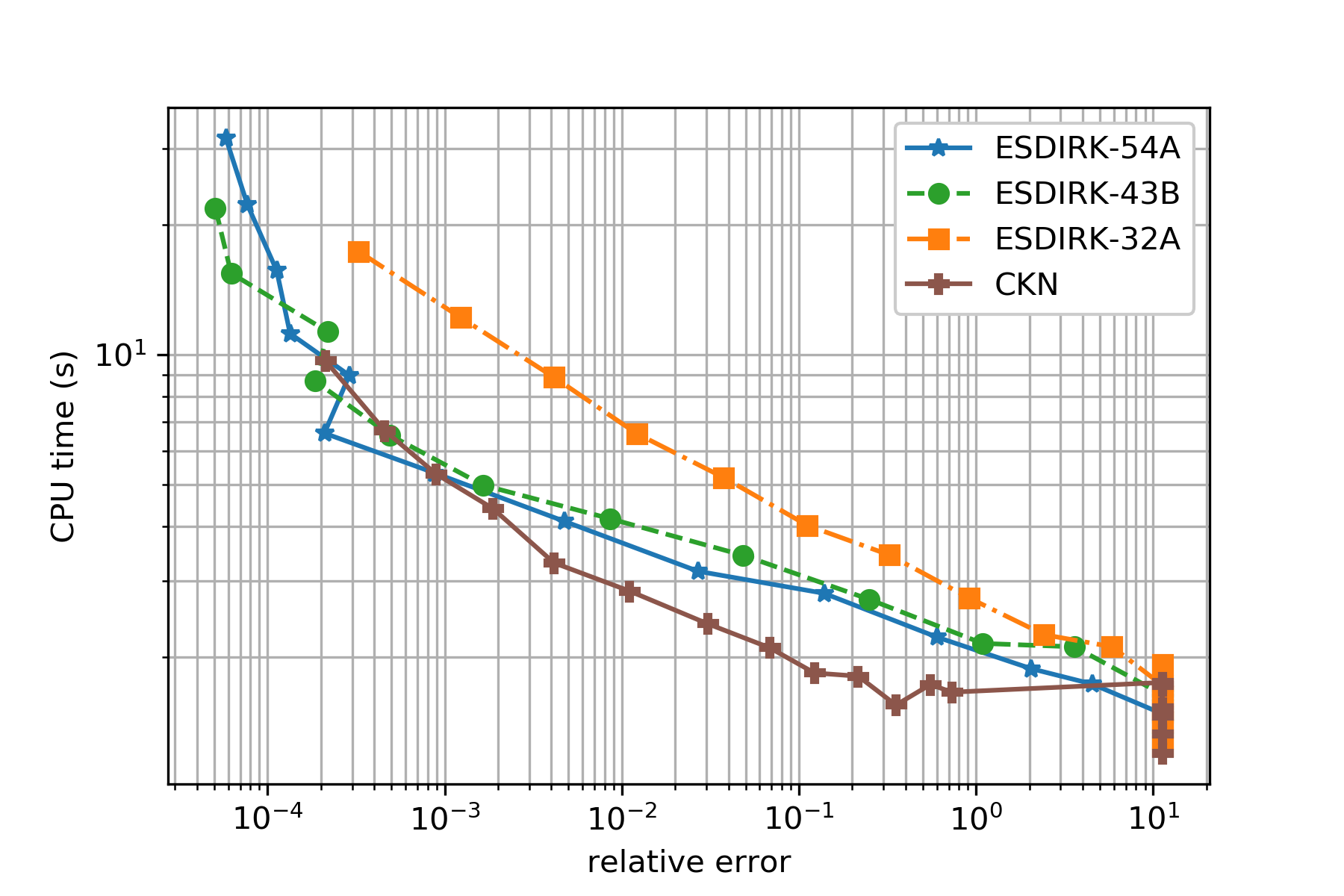}
\caption{Amplification factor (transient-only)\\}
\label{fig:limitcycle:work_ampli}
\end{subfigure}
~
\begin{subfigure}[t]{0.45\textwidth}
\centering
\includegraphics[width=\textwidth]{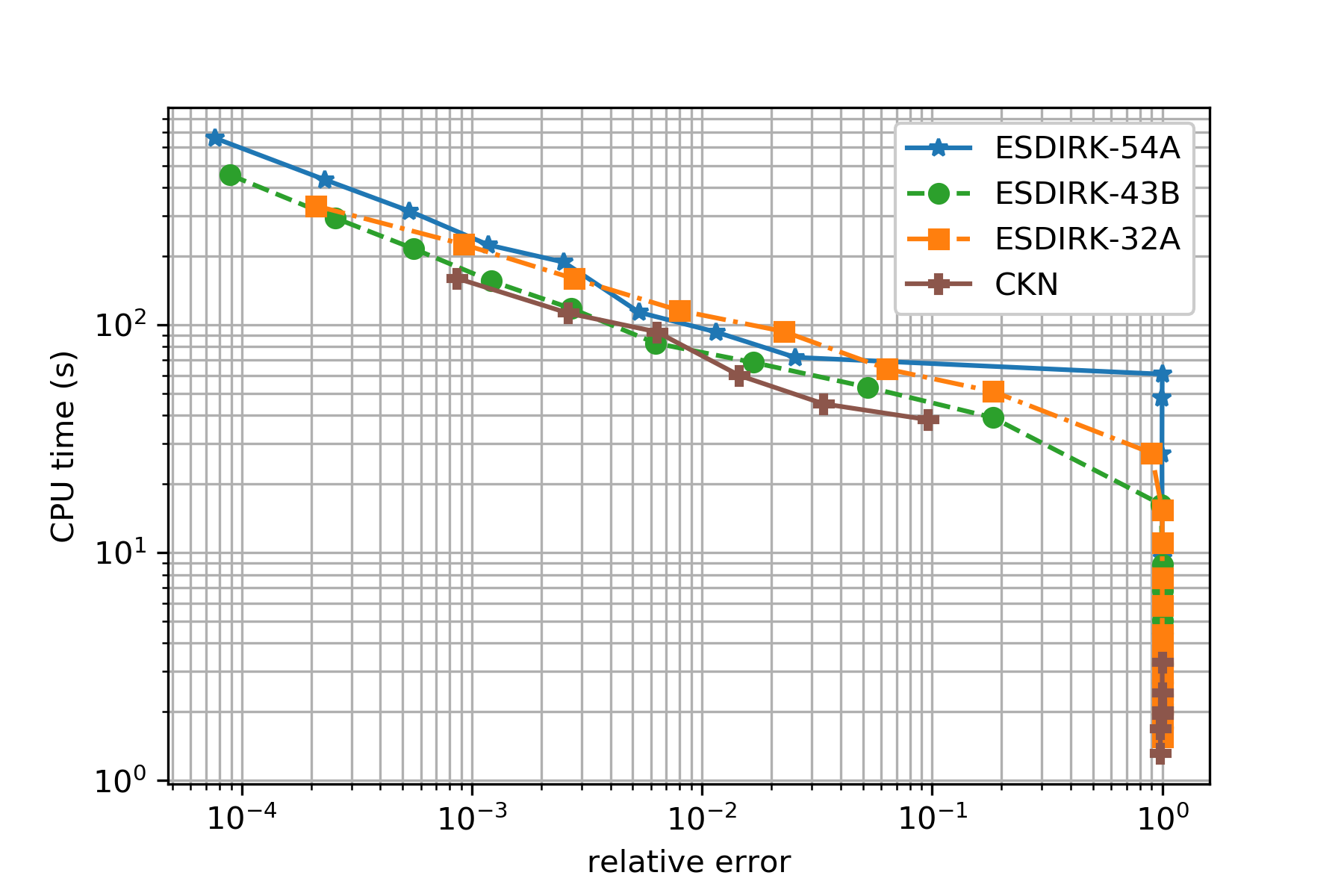}
\caption{Amplitude of the fundamental of the established limit cycle}
\label{fig:limitcycle:work_amplitude_harm0}
\end{subfigure}
\caption{Work-precision diagrams for fixed time step simulations}
\end{figure}

Based on the previous analysis, the high-order methods ESDIRK-54A and ESDIRK-43B seem particularly promising, as they require fewer steps to achieve good results.
However, since these methods possess more stages than lower-order ones, this does not ensure that the actual computational times will be shorter.
Figure \ref{fig:limitcycle:work_ampli} shows the computational time versus the achieved relative error on the amplification factor. Note that these simulations were run only on the physical time interval $t\in[0,0.1]$ s, so that computational times are truly representative. We see that CKN is the fastest method for a relative error higher than $5 \times 10^{-3}$, however this roughly corresponds to the zone were its solution diverges in finite time. 
ESDIRK-32A is not a very good performer, whereas ESDIRK-43B and ESDIRK-54A are performing well and have similar error levels.

Figure \ref{fig:limitcycle:work_amplitude_harm0} shows the computational time required for a given level of relative error on the fundamental amplitude in the established limit cycle. Computational times are those of simulations run on the physical time interval $t\in[0,1.5]$ s. ESDIRK-43B is only marginally better than the other methods, no clear winner is to be picked.

Overall, when using constant steps, the high-order methods ESDIRK-54A and ESDIRK-43B are almost identical and offer overall a very good performance. The Crank-Nicolson method is slightly misleading: its lack of stability when applied to DAEs leads to an easier destabilisation of  the initial solution.
However it diverges quickly, unless the time step is severly reduced, thus falsely leading to the conclusion that the configuration is purely unstable and does not produce a limit cycle.

\subsection{Numerical experiment with time adaptation}
The previous study with constant steps has shown that high-order methods are interesting for the simulation of a limit cycle. 
We now compare the embedded ESDIRK schemes with time adaptation enabled, to see if additional computational gains can be obtained. Following the methodology exposed in Section \ref{section:RK:adaptation}, the time step is controlled by the relative integration error tolerance $rtol$, which is varied between $10^{-1}$ and $10^{-7}$. We first focus on the relative error achieved on the quantitative criteria used in the previous section. Finally, a comparison of the computational times and relative errors is presented, considering both fixed time step and adaptive simulations.

\subsubsection{Initial growth}

\begin{figure}[hbt!]
\centering
\begin{subfigure}[t]{0.45\textwidth}
\centering
\includegraphics[width=\textwidth]{convergence_amplification_factor.png}
\caption{Convergence of the initial growth factor}
\label{fig:limitcycle:dt_adaptatif:convergence_ampli_factor}
\end{subfigure}
~
\begin{subfigure}[t]{0.45\textwidth}
\centering
\includegraphics[width=\textwidth]{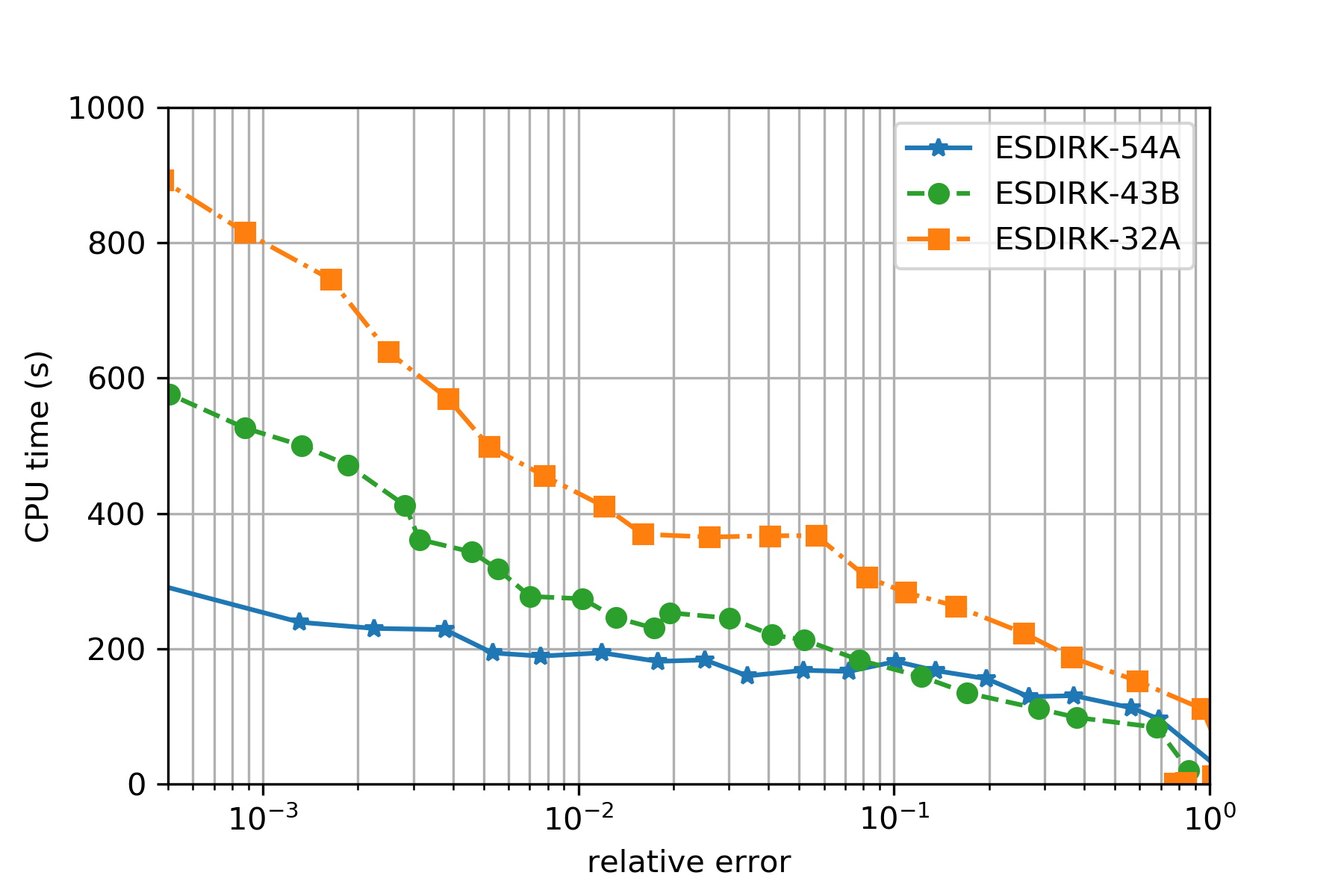}
\caption{Work-precision diagram}
\label{fig:limitcycle:dt_adaptatif:work_ampli_factor}
\end{subfigure}
\caption{Convergence and computational cost for the amplification factor}
\end{figure}

Figure \ref{fig:limitcycle:dt_adaptatif:convergence_ampli_factor} shows the convergence of the fitted amplification factor $\amplifactor$ when $rtol$ is decreased. We see that with fine tolerances, all three methods resolve the transient quite well.
The achieved relative error is proportional to $rtol$ when the latter is sufficiently small, which indicates that the error estimation based on the embedded schemes behaves correctly. However, we see that there is a difference in each method's error levels for a given value of $rtol$, which is expected as they have different leading truncation error coefficients.
In Figure \ref{fig:limitcycle:dt_adaptatif:work_ampli_factor}, we plot the accuracy achieved on the amplification factor (which is not equal to $rtol$) versus the computational time. We clearly see that, when the tolerance is sufficiently low, the computational cost decreases as the order of the method increases. 
For example, if we require a relative error of $10^{-3}$, ESDIRK-54A is twice as fast as ESDIRK-43B , and three times as fast as ESDIRK-32A. Only for high levels of error ($>10^{-1}$) is ESDIRK-43B slightly more efficient than ESDIRK-54A.

\subsubsection{Limit cycle}

\begin{figure}[hbt!]
\centering
\begin{subfigure}[b]{0.45\textwidth}
\centering
\includegraphics[width=\textwidth]{convergence_amplitude_harm0.png}
\caption{}
\label{fig:limitcycle:dt_adaptatif:convergence_amplitude0}
\end{subfigure}
~
\begin{subfigure}[b]{0.45\textwidth}
\centering
\includegraphics[width=\textwidth]{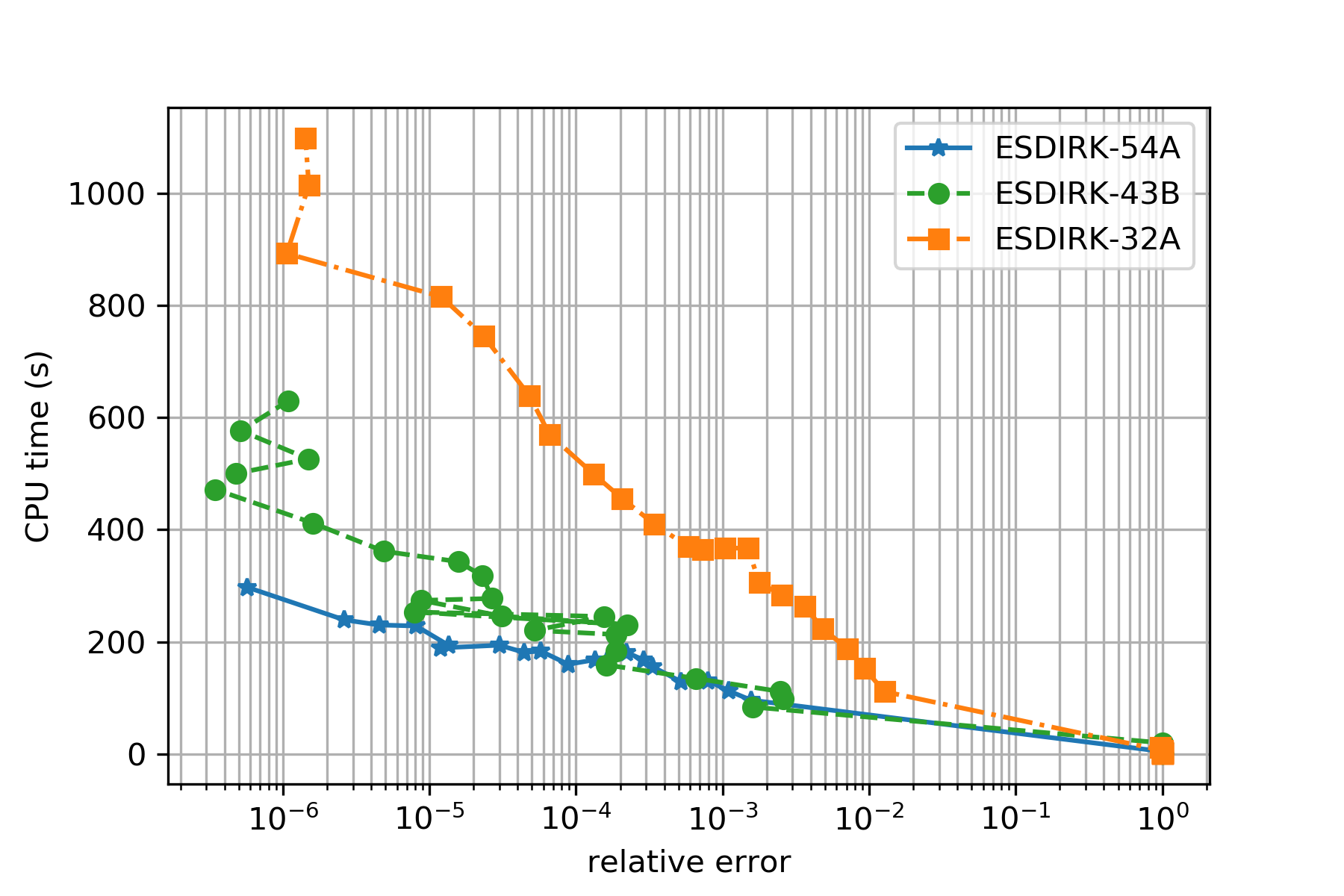}
\caption{}
\label{fig:limitcycle:dt_adaptatif:work_vs_error_amplitude0}
\end{subfigure}
\caption{Limit cycle with adaptive time stepping: (a) convergence of the fundamental amplitude, (b) work-precision diagram}
\end{figure}

We now compare the computational cost of each method when considering the resolution on the time interval $t\in [0, 1.5]$ s.
We have observed that the frequency of the fundamental and its amplitude converge equally well, therefore we only focus on the amplitude.
Figure \ref{fig:limitcycle:dt_adaptatif:convergence_amplitude0} shows how the relative error on the fundamental amplitude evolves with $rtol$. ESDIRK-43B has, for a given $rtol$, the lowest error, however there is an unexplained oscillation of the relative error.
Figure \ref{fig:limitcycle:dt_adaptatif:work_vs_error_amplitude0} is the corresponding work-precision diagram. ESDIRK-32A is the worst performer by far, whereas ESDIRK-54A is the most efficient method.

\subsubsection{Comparison with fixed time step results}
\label{sec:limitcycle:comparison_fixed_dt}
\begin{figure}[hbt!]
\centering
\begin{subfigure}[t]{\textwidth}
\centering
\includegraphics[width=0.6\textwidth]{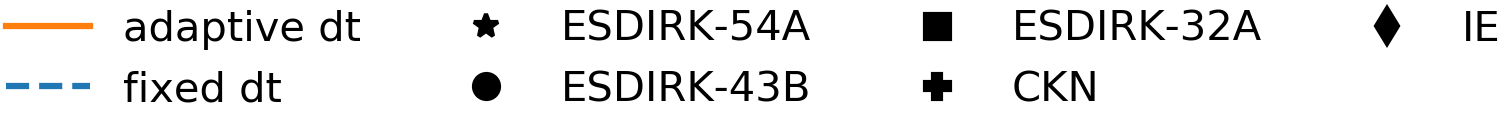}
\end{subfigure}

\begin{subfigure}[t]{0.45\textwidth}
\centering
\includegraphics[width=\textwidth]{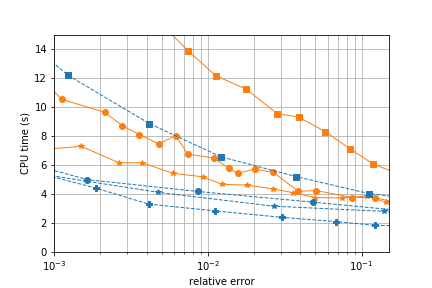}
\caption{Fitted amplification factor}
\label{fig:limitcycle:comparaison:amplification_factor}
\end{subfigure}
~
\begin{subfigure}[t]{0.45\textwidth}
\centering
\includegraphics[width=\textwidth]{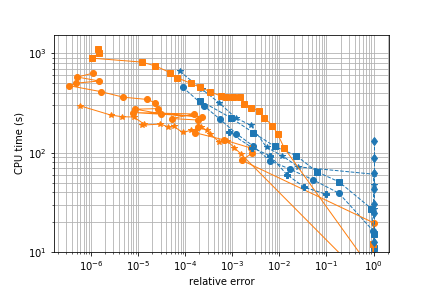}
\caption{Fundamental amplitude in established limit cycle}
\label{fig:limitcycle:comparaison:amplitude_harm0}
\end{subfigure}
\caption{Comparison of the computational cost for a given level of relative error}
\end{figure}

We now wish to assess the performance gain achieved with time adaptation. To this end, we compare the computational times between fixed time step simulations and adaptive simulations for a given level of relative error.
Figure \ref{fig:limitcycle:comparaison:amplification_factor} shows how the relative error on the amplification factor $\amplifactor$ during the initial growth evolves with computational time. The simulations were run for only 0.1 s of physical time, i.e. only for the initial growth.
We observe that fixed time step implementations of the ESDIRK methods are faster than their adaptive counterparts.
It has been assessed that this was due to an advantage in terms of number of Jacobian evaluations for the Newton method.
Indeed, for the short time range simulated, the instability remains linear, i.e. the instability is linked with constant positive eigenvalues of the Jacobian of the system. Hence, the Jacobian does not need to be updated if the time step is constant. In practice, fixed time steps solutions only required up to 3 evaluations of this matrix, while adaptive solutions required up to 100 evaluations due to the repeated changes in time step.
It is possible to mitigate this issue by adding a ``deadzone'' for the time step evolution, i.e. refuse time step increases that are lower than a given percentage.
In our testing, a 50 \% deadzone yields a 25\% CPU time improvement, dividing by 3 the number of Jacobian evaluations. However, such a deadzone approach did not perform very well in the ignition test case from Section \ref{section:allumage:simplifie}, hence we do not consider it further. Finally, the Jacobian update process could be improved.
Currently, as discussed in Section \ref{section:cost},
every time the Jacobian $J$ of the residuals $F$ needs to be updated, the Jacobians of \modif{$f$, $g$ and $Q$ are} also updated.
In the case of nearly linear dynamics, large computational gains could be expected from simply refactorising $J$ upon each time step change, and only updating the stored Jacobians of $f$, \modif{$g$ and $Q$ when the Newton fails to converge properly.} \modif{This would however require a number of modifications in our code, so that $f$, $g$ and $Q$ can be evaluated separately. Furthermore, in our typical use cases (e.g. ignition transients as in Section \ref{section:allumage:simplifie}), we observed that our code often performs successful time steps (i.e. with a sufficiently low estimated temporal error) during which the Newton algorithm has required one or more Jacobian updates for convergence for one or multiple stages of the Runge-Kutta step.
This means that nonlinearities encountered across a single step are already sufficiently large to require a full Jacobian update, hence we expect limited gains from being able to only refactor $J$ upon time step changes in these scenarios.}

The performance for the computation of the full limit cycle (initial growth and established cycle) is assessed in Figure \ref{fig:limitcycle:comparaison:amplitude_harm0}. The criterion is the relative error on the amplitude of the fundamental in the established limit cycle.
One may think that the lightweight second-order Crank-Nicolson method could outperform the other methods for the established limit cycle, as this method is known to have good damping properties for oscillating systems, while only requiring one stage to be computed per time step. Indeed the method is the fastest among the fixed time step ones, and the fastest overall for relative error levels around $10^{-2}$. For lower error levels however, the adaptive high-order methods ESDIRK-54A and ESDIRK-43B are the most efficient methods. ESDIRK-32A in adaptive mode is generally slower than in fixed time step mode, unless very low errors are sought. Although not shown here, adaptive simulations require more Jacobian evaluations due to the successive changes in time step, but they require many fewer steps and Newton iterations overall, which, for the complete simulation, far outweighs the drawback of evaluating the Jacobian more often.
This is supported by the observation that the slopes of CPU time versus relative error are smaller for adaptive methods compared to fixed time step implementations. Any additional cost, e.g. Jacobian evaluation, is overcome by the ability to take fewer steps.
Still, adaptive methods are at a slight disadvantage in this test case, as the solution oscillates smoothly: the characteristic time scale of the system stays roughly constant throughout the simulation, therefore fixed time step simulations, with $dt$ sufficiently low compared to this time scale, will be favoured by this consistency.

An interesting observation can be made: the adaptive methods always capture the correct limit cycle, unless $rtol$ is too high, leading to a stabilisation of the solution. This is seen in Figure \ref{fig:limitcycle:comparaison:amplitude_harm0}, as all adaptive methods have a jump from important errors ($\approx 1$) to much lower ones as $rtol$ is lowered. On the opposite, fixed time step implementations do not have such a jump in error and are more likely to capture a non-accurate limit cycle (typically with an error higher than $10^{-2}$) for an intermediate range of time step values.
Overall, ESDIRK-54A and ESDIRK-43B seem to be the most reliable methods in this comparison.

Another practical consideration is that the time step corresponding to $\mathrm{CFL}=1$ lies around $10^{-6}$ s, which is approximately 100 times smaller than the time step necessary to obtain a very precise simulation with the fifth-order method ESDIRK-54A (see Figures \ref{fig:limitcycle:convergence_amplification} and \ref{fig:limitcycle:convergence_amplitude_harm0}).
The CFL-controlled time step is based on a stability analysis of the convection operator with an explicit time integration, which is not relevant for implicit integration and does not guarantee any level of error on the solution. As already discussed for the ignition transient in Section \ref{section:allumage:simplifie}, use of a CFL-limitation would result in an important increase in computational time, without any valuable improvement on the solution accuracy.
Finally, the constant time step simulations that yield accurate results more efficiently than with adaptive methods correspond to $\mathrm{CFL} \approx 100$, which would usually not be expected to produce accurate unsteady results. One would rather safely choose a time step such that $\mathrm{CFL}=1$. This again highlights one practical benefit of time adaptation: precision is ensured based on a reliable mathematical criterion, and time step values can be used such that $\mathrm{CFL} \gg 1$, while still ensuring a precise solution.

\section{Application to unsteady combustion with detailed chemistry}
\label{section:chimie_complexe}
We have now verified that the high-order adaptive methods perform well on nonlinear problems involving simplified chemistry models.
Additional stiffness is usually observed when a complex kinetic mechanism is used, thus it is useful to check the behaviour of the proposed numerical strategy in this context.
The test case is the unsteady combustion of the AP monopropellant.
Gas-phase kinetics is based on the AP-HTPB mechanism developed by Jeppson \cite{Jeppson_mecanisme_AP} and Tanner \cite{these_Tanner}.
All reactions involving carbonated species were removed to account for the absence of HTPB, resulting in a pure AP combustion mechanism involving 25 species and 80 reactions.
Thermodynamic and transport properties are computed before-hand by CHEMKIN routines \cite{ref_CHEMKIN2} and stored as lookup tables.

The solid phase and the surface are handled as in \cite{theseShihab}: the solid is assumed inert, and all decomposition and gasification reactions occur at the surface. There are two global surface reactions: a direct dissociative sublimation, and a quasi-equilibrium decomposition. The regression speed is defined by a pyrolysis law taken from \cite{theseShihab}, and the proportions of gaseous products generated by the surface reactions are adjusted to obtain the experimentally measured regression rates at $20$ atm \cite{meynet2006}.
These modelling choices allow for the use of detailed combustion kinetics while remaining within the comparatively simpler framework of solid and surface representation used in the present paper.

The computational mesh has 49 cells for the solid phase, and 126 cells for the gas phase, distributed in a non-uniform manner so that steady-state gradients are well resolved.
Starting from a steady-state solution at $P=20.265 \times 10^5$ Pa, we study the transient occurring after a pressure step to $P=20\times 10^5$ Pa.

\subsection{Order of convergence}

The orders of convergence of the methods from Table \ref{table:rk_integrators} have already been verified for the simpler test case of Section \ref{section:simplemodel} and are presented in \ref{appendix:convergence}. We now verify that the orders are not affected by the additional complexity and stiffness induced by detailed kinetics.
We simulate the unsteady evolution for $t\in [0, 0.2]$ s and perform multiple integrations with various time steps.
To quantify the accuracy of the overall time integration, we define the following error: $\epsilon_{T_\surface} = ({1}/{t_f}) \int_0^{t_f} |T_{\surface}-T_{\surface,ref}| dt$,
with $ref$ designating the reference simulation and $t_f$ the final physical time.
The reference simulation is computed with ESDIRK-54A and $\Delta t=10^{-5}$ s.
Cubic interpolation is used to compare both solutions on the same time grid.

\begin{figure}[hbt!]
\centering
\begin{subfigure}[t]{0.45\textwidth}
\includegraphics[width=\textwidth]{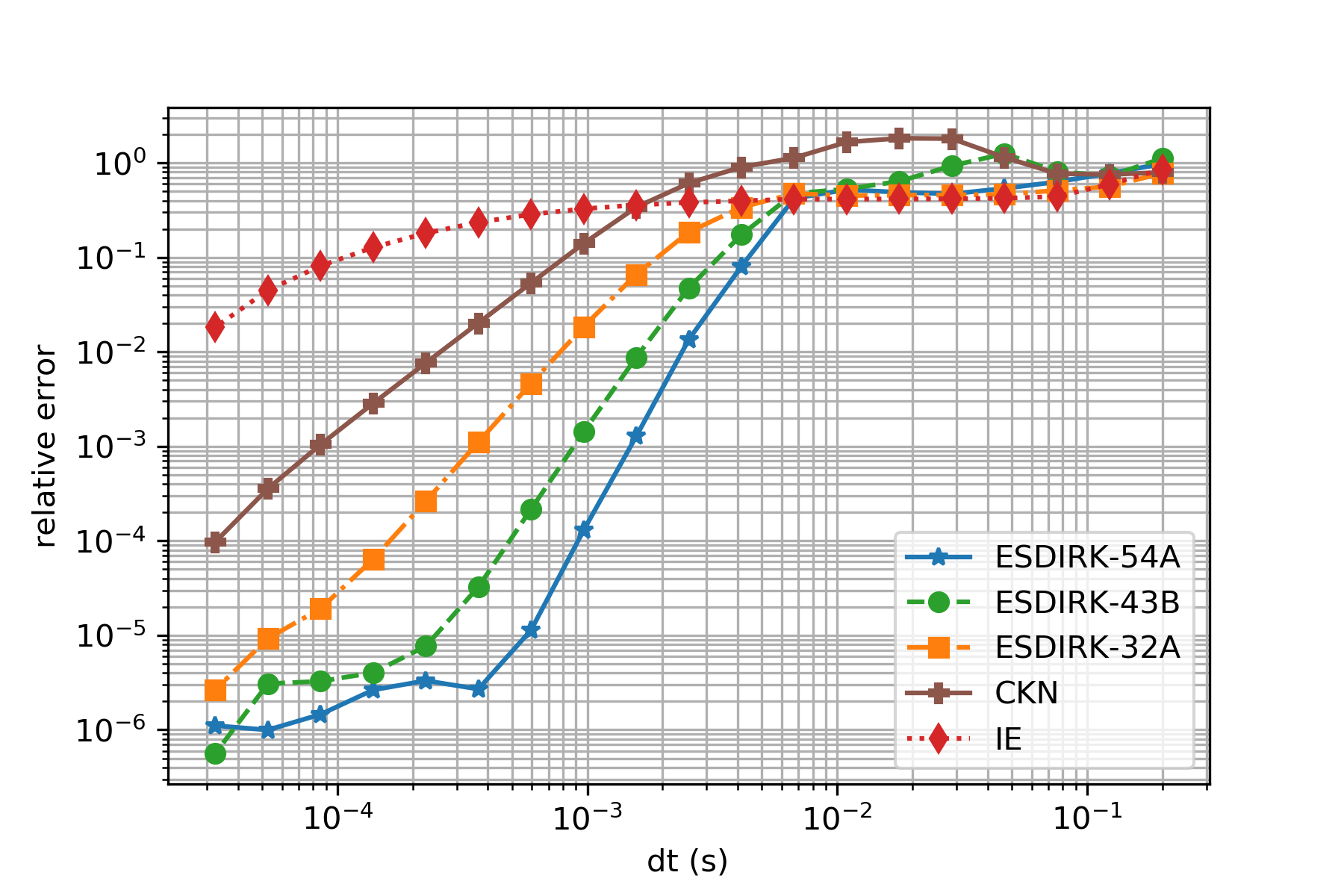}
\caption{}
\label{fig:chimie_complexe:dt_fixe:error_Ts_global}
\end{subfigure}
~
\begin{subfigure}[t]{0.45\textwidth}
\includegraphics[width=\textwidth]{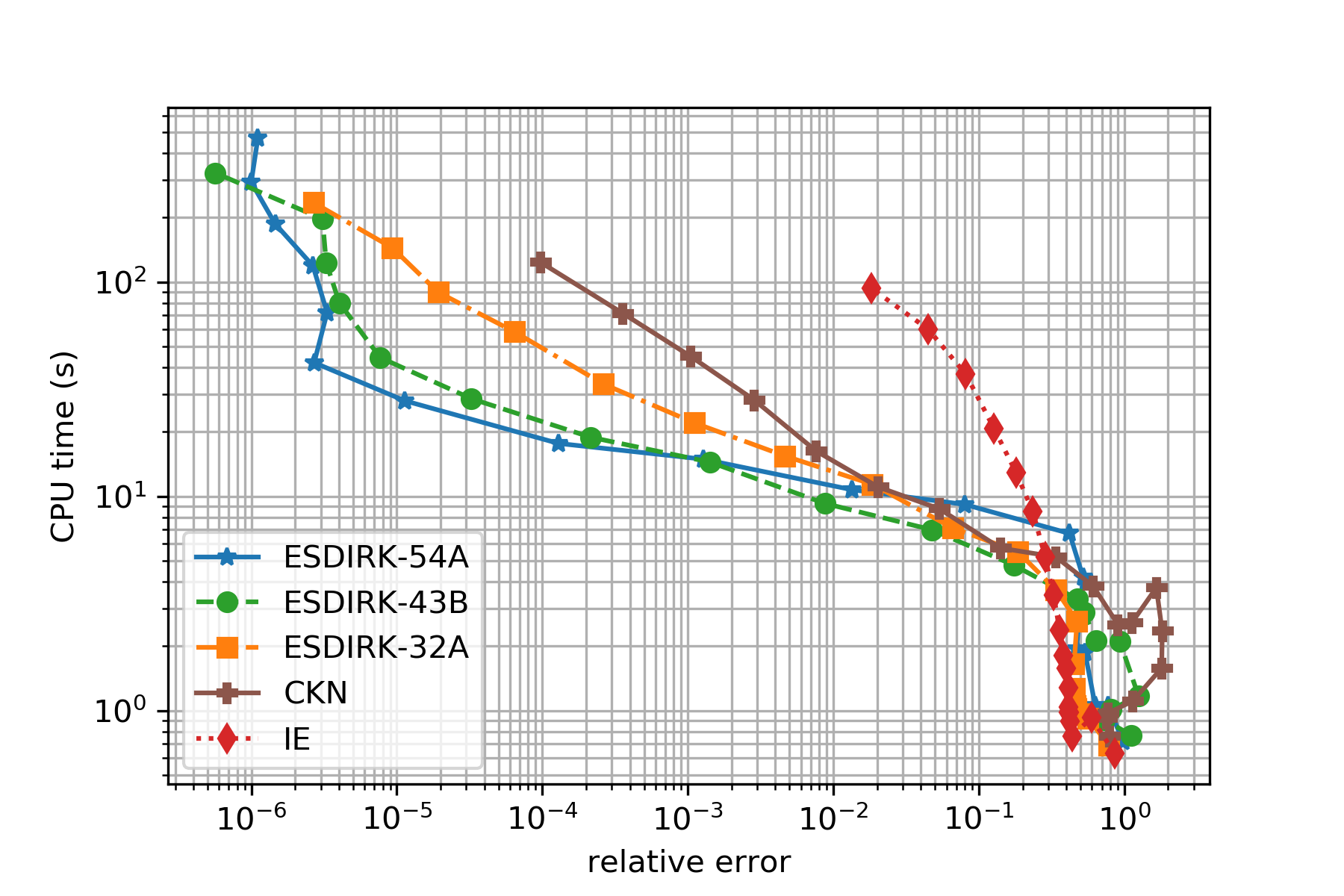}
\caption{}
\label{fig:chimie_complexe:dt_fixe:work_precision_Ts_global}
\end{subfigure}
\caption{Accuracy of the integration with fixed time steps: (a) Convergence of $\epsilon_{T_\surface}$,
(b) Work-precision diagram for $\epsilon_{T_\surface}$}
\end{figure}

Figure 
\ref{fig:chimie_complexe:dt_fixe:error_Ts_global} shows the evolution of this error when the time step is varied. We see that each method attains its theoretical order of convergence. Similar results have been obtained for the other variables, both differential or algebraic.
No order reduction is observed due to the stiffness of the chemical reactions with detailed kinetics. The ESDIRK methods perform  well in this more complex scenario.

Figure \ref{fig:chimie_complexe:dt_fixe:work_precision_Ts_global} shows for each method the computational time required to achieve a given integration error $\epsilon_{T_\surface}$. 
By analysing the different simulations, it was determined that the curve of $T_\surface$ is visually converged when $\epsilon_{T_\surface} \leq 10^{-2}$.
ESDIRK-43B and ESDIRK-54A are the most efficient methods on this error range, and the speed-up achieved by these high-order methods increases with the precision achieved. They are about 1.5 times faster than ESDIRK-32A and CKN for $\epsilon_{T_\surface} = 10^{-2}$, and approximately 5 times faster than CKN for $\epsilon_{T_\surface} = 10^{-4}$.
The poor performance of IE clearly advocates high-order methods.

\subsection{Computational performance with time step adaptation}
Now that we have verified that the convergence of the methods is not affected by the stiffness induced by complex kinetics, we use the ESDIRK methods with time step adaptation to see how they compare in terms of results. Different values of the relative integration tolerance $rtol$ are used between $10^{-1}$ and $10^{-7}$.

\begin{figure}[hbt!]
\centering
\begin{subfigure}[t]{0.45\textwidth}
\includegraphics[width=\textwidth]{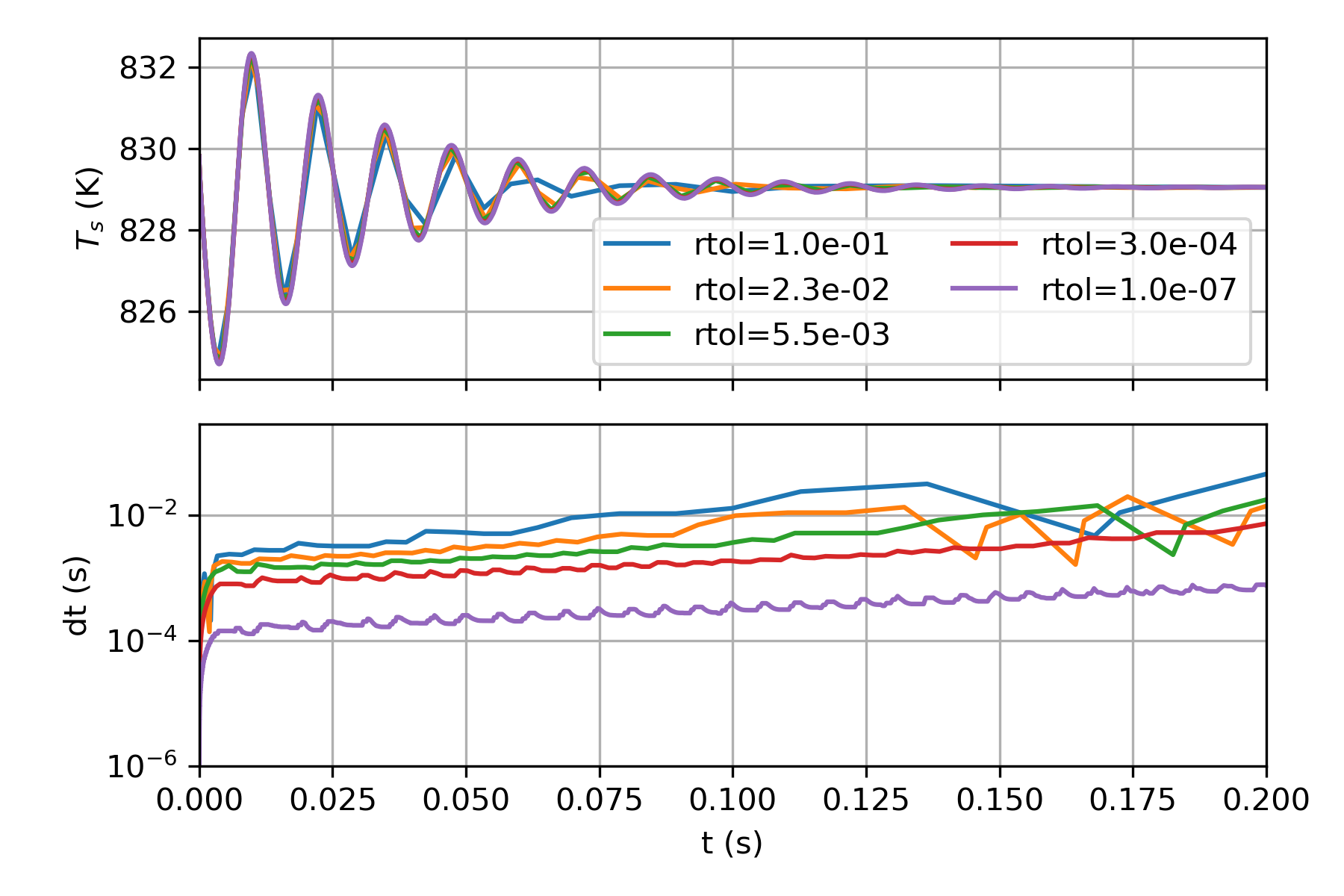}
\caption{}
\label{fig:chimie_complexe:adaptation}
\end{subfigure}
~
\begin{subfigure}[t]{0.45\textwidth}
\includegraphics[width=\textwidth]{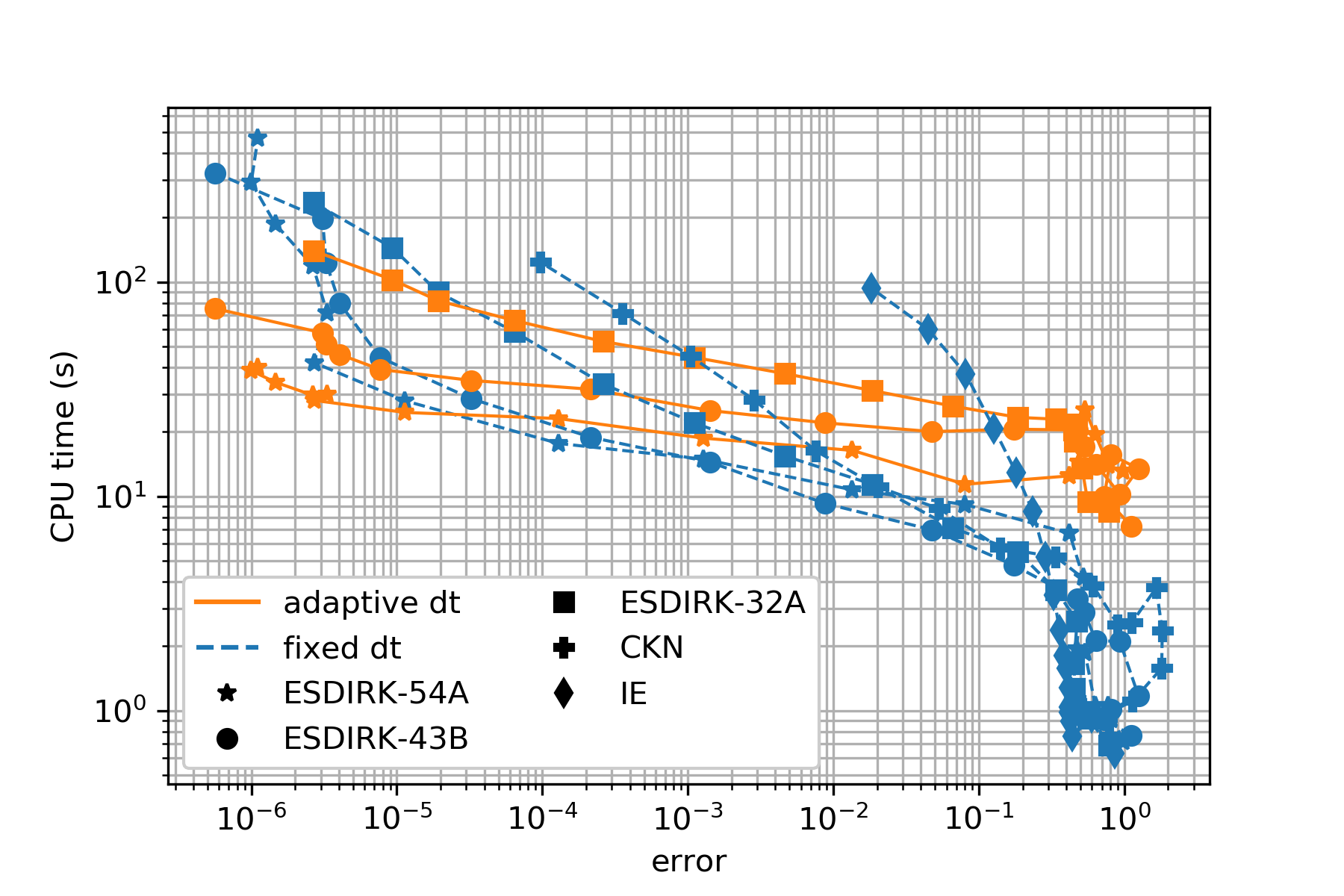}
\caption{}
\label{fig:chimie_complexe:comparaison:work-precision}
\end{subfigure}
\caption{Integration with adaptive time stepping: (a) time step evolution for ESDIRK-43B, (b) work-precision diagram for $\epsilon_{T\surface}$}
\end{figure}

Figure \ref{fig:chimie_complexe:adaptation} shows the complete transient for the surface temperature and the time step evolution for various values of $rtol$. We see that the change in time step is smooth, except for low tolerances when the time step becomes large, causing convergence issues. The temporal evolution of the surface temperature is well resolved even with relatively large values of $rtol$.
Figure \ref{fig:chimie_complexe:comparaison:work-precision} shows the comparison of the computational time required to achieve a given level of error $\epsilon_{T_\surface}$, both with fixed time steps (blue lines) and adaptive time stepping (orange lines).
Here, adaptive schemes do not seem to improve the performance globally. ESDIRK-54A is the best performing adaptive method, however it only becomes the fastest method overall for a very low level of error $\epsilon_{T_\surface} \leq 10^{-5}$. Its computational time is relatively close to the one of its fixed time step implementation.
We observe that, for a given increase in accuracy, adaptive methods have a lower increase in computational time compared to their fixed time step counterparts.

This comparison is however slightly unfair, because we had no a priori knowledge of the time step needed to properly resolve the transient. Also, the characteristic times of the dynamics does not vary much , as we can see from the small variations of the adapted time step in Figure \ref{fig:chimie_complexe:adaptation}. Furthermore, the pressure step is sufficiently low so that nonlinearities remain small, allowing for very few Jacobian updates to be performed when using fixed time steps.
Without dynamic adaptation, the time step would have typically been limited so that the CFL number is reasonably low, e.g. 1 to 10. A unitary $\mathrm{CFL}$ corresponds to $\Delta t \approx 7\times 10^{-8}$ s, which is much lower than the time step required to achieve a good accuracy with most methods.
A simulation has been carried out only on the first 0.03 s of the transient, with CKN and a time step set such that $\mathrm{CFL}=10$. The computational time was 660 s.
An equivalently well resolved transient can be obtained with ESDIRK-54A and $rtol=10^{-5}$ in only 11 s. As we can see in Figure \ref{fig:chimie_complexe:adaptation}, even a less stringent tolerance would also be sufficient.
From an engineering point of view, this represents a 60 times speed-up, due to the fact that the time adaptation will automatically choose the relevant time step values. This guarantees a proper accuracy, while reaching CFL numbers that one would usually never trust to yield accurate unsteady results. Finally, time adaptation based on embedded methods automatically detects  a slowdown of the dynamics as the solution stabilises and is able to increase the time step accordingly, whereas the CFL number stays roughly constant and cannot be an efficient time-step controlling criterion in that situation.
This gain in engineering time is not quantifiable precisely, however it is definitely important.

\section{Conclusion}
This contribution presents the development of a high-fidelity simulation tool for the one-dimensional low-Mach combustion, with a particular focus on solid propellant applications. The emphasis is on the numerical strategy for the time integration of the semi-discretised equations.
It has been shown that the system is differential-algebraic in nature. Multiple test cases show that stiffly accurate singly-diagonally implicit Runge-Kutta methods are highly efficient for the time integration of such a system, in particular the embedded ESDIRK methods presented in \cite{Kvaerno2004}.
Handling the continuity equation with a Runge-Kutta quadrature instead of the classic instantaneous reformulation leads to greater computational efficiency.
Applications have been presented for ignition transients and limit cycle with a simplified modelling, and appreciable computational gains have been observed. High-order methods reliably capture dynamics which are practically impossible to reproduce with traditional low-order methods.
It has also been verified that the proposed numerical strategy is robust and performs well when the modelling is much more complex, e.g. detailed kinetics in the gas phase.

Adaptive time stepping based on objective error estimates ensures the temporal evolution is well resolved.
From an engineering point of view, the single parameter that controls the time step is the relative integration error tolerance $rtol$. In all our test cases, we have observed that $rtol = 10^{-5}$ is sufficient to accurately resolve all unsteady phenomena. Using this value as standard tolerance liberates from the need of iterating over other practical criteria such as CFL limitation or maximum relative variation. We believe that a high-order adaptive method like ESDIRK-54A therefore allows for perceivable gains in computational time, trustworthiness of the results, and engineering time spent parametrising the time integration for a simulation. Furthermore, the proposed numerical strategy is easy to implement in an existing code if the latter already uses an implicit Euler, Crank-Nicolson or BDF scheme, as is often the case in the literature.

Future work includes applying the presented framework to one-dimensional models involving a multiphase foam layer at the surface, and adapting the approach to higher-dimensional detailed heterogeneous combustion codes for solid propellants,  e.g. \textsc{Compas}\cite{dmitry_compas}
from ONERA.
The code \textsc{Vulc1D} has already been successfully coupled as a dynamic boundary condition for combustion chamber simulations with the 3D multiphysics CFD tool \textsc{Cedre} from ONERA \cite{CEDREsoftware}, avoiding costly CFD mesh refinement near the propellant surface \cite{aiaa2021,maThese}. 
Eventually, the proposed time strategy is not dependent on the chosen spatial discretisation and other spatial schemes could be envisioned without changing the conclusion of our study.  As mentioned in the course of the paper, the numerical strategy can be applied to any other one-dimensional combustion problem in the low-Mach limit involving homogeneous or spray combustion. When index-2 algebraic variables are involved, the order obtained on such variables is limited by the stage order of the method \cite{Petzold_DAE,hairer_book2}. 
This situation occurs either for the strain rate eigenvalue for counterflow diffusion flames or in multi-dimensional Navier-Stokes equations, either incompressible or in the low-Mach limit \cite{nguessan2020,nguessan2019}; however the proposed strategy should be equally of interest in such cases.

\section*{Acknowledgments}
The present research was conducted thanks to a Ph.D grant co-funded by DGA, Ministry of Defence (E. Faucher, Technical Advisor), and ONERA.
The authors would like to thank Gilles Vilmart for interesting technical discussions.

\bibliography{references}

\appendix

\section{Generating configurations with various degrees of instability}
\label{appendix:explications_ligneRK}
 In order to highlight the benefit of the high-order adaptive time integration, we search for configurations which are linearly unstable around the corresponding steady-state solution. We generate such configurations by taking the stable simplified model from Section \ref{section:simplemodel:parameters} as baseline and varying its parameters. We use existing theoretical tools to approximately evaluate the stability of the steady-state solution.
 
 \subsection{Theoretical indicator of intrinsic instability}
\label{appendix:ZN}
The Zeldovich-Novozhilov (ZN) framework \cite{novozhilov1992} is a useful tool to study the stability of a steady-state solid propellant combustion .
It assumes that $\massflux$ and $T_\surface$ are linked in steady-state via laws of the form
$\massflux = \massflux(T_\initial, P)$ and $T_\surface = T_\surface(T_\initial, P)$.
Considering the steady-state temperature profile $T(x) = T_\initial + (T_\surface - T_\initial) \exp\left(x{\massflux c_c}/{\lambda_c}\right)$, the temperature gradient just below the surface is
$\partialdershort{T}{x}(0^-) = {\massflux c_c (T_\surface - T_\initial)}/{\lambda_c}$.
Replacing $T_\initial$ by $T_\surface - \phi {\lambda_c}/{\massflux c_c}$ (the ``apparent'' or ``initial'' temperature), the previous relations can now be used in the unsteady regime.
This is usually accepted, as long as the apparent initial temperature remains within acceptable bounds. It is also required that data for this initial temperature be available, or at least reasonably extrapolated.

The amplifications of compact perturbations can be studied by linearising the previous relations and the heat equation in the solid, leading to the definition of a stability criteria, which depends on two steady-state sensitivity coefficients: 
$\rZN = \left( \partialdershort{\overbar{T_\surface}}{T_\initial} \right)_{P}$, $\kZN =  ( \overbar{T_\surface} - T_\initial ) \left( \partialdershort{ln(\overbar{\massflux})}{T_\initial} \right)_{P}$,
with $\overbar{~\cdot^{~}}$ denoting steady-state values.
Steady-state combustion is always stable if $\kZN<1$.
If $\kZN>1$, the steady-state is stable only if $\rZN > {(\kZN-1)^2}/{(1+\kZN)}$. The line $\rZN={(\kZN-1)^2}/{(1+\kZN)}$ is the locus of a Hopf bifurcation, where the steady-state solution becomes linearly unstable in a oscillating manner, with the possibility of stabilising on a limit cycle.
If $\rZN > (\sqrt{\kZN}-1)^2$, the instability grows purely exponentially.
The associated stability diagram is shown in Figure \ref{fig:optim:ligneRK}: the leftmost parabola is the first stability limit, the second one is the onset of purely exponential instability.
This stability is called ``intrinsic'' because it is a property of the solid propellant as an isolated system,
and is not caused by the coupling with another system (e.g combustion chamber \cite{culickInstabilities}).

The original ZN approach assumes a quasi-steady gas phase.
Even though it has been shown that unsteady gas-phase phenomena tend to slightly widen the stability area, this first simplified analysis remains a good indicator of the stability bounds.
We refer the reader to \cite{novozhilov1992} for a detailed presentation of the ZN analysis and its extensions.

\subsection{Optimisation problem}
\newcommand{\varoptim}{\ensuremath X}
Let us denote as $\varoptim$ the vector
containing the parameters of the simple combustion model that we have chosen as free variables. For a given value of $\varoptim$, we can find the corresponding value of $(\rZN,\kZN)$ by performing three steady-state simulations: one baseline simulation, one simulation with a perturbed initial temperature $T_\initial$, and one simulation with a perturbed pressure $P$.
Then, by means of finite differences, $\rZN$ and $\kZN$ may be evaluated.
This process can be summarised as the function $f_{\rZN\kZN} : \varoptim \rightarrow (\rZN,\kZN)$.
We can then generate configurations with various values of $\rZN$ and $\kZN$ by solving the following optimisation problem:
\begin{subequations}
\begin{alignat}{2}
&\!\min_{\varoptim}        & \qquad & f_{obj}(\varoptim)       \label{eq:optProb}\\
&\text{subject to} &        & g(\varoptim) \leq 0 \label{eq:constraint1}\\
&                  &        & h(\varoptim) = 0      \label{eq:constraint2}
\end{alignat}
\end{subequations}
\noindent where the objective function is $f_{obj} : \varoptim \rightarrow || f_{\rZN\kZN}(\varoptim) - (\rZN,\kZN)_{target} ||_2^2$ with $(\cdot,\cdot)$ denoting a vector formulation. 
This is simply the constrained minimisation of the distance to the target $(\rZN,\kZN)$ coefficients.
Inequality constraints $g$ ensure the physical parameters remain within realistic bounds. They can be supplemented with equality constraints $h$ to enforce certain properties, e.g. steady-state surface temperature.
The problem is solved with SLSQP from the Scipy library \cite{2020SciPy-NMeth}, using finite-difference Jacobians for $f$, $q$ and $h$.

\subsection{Numerical assessment of intrinsic stability}
\label{section:appendix:instab_config}
\begin{figure}[htpb]
\begin{center}
	\begin{minipage}[t]{0.45\textwidth}
		\centering
		\adjincludegraphics[width=\textwidth, trim={{0\width} {0\height} {0\width} {0.1\height}}, clip]{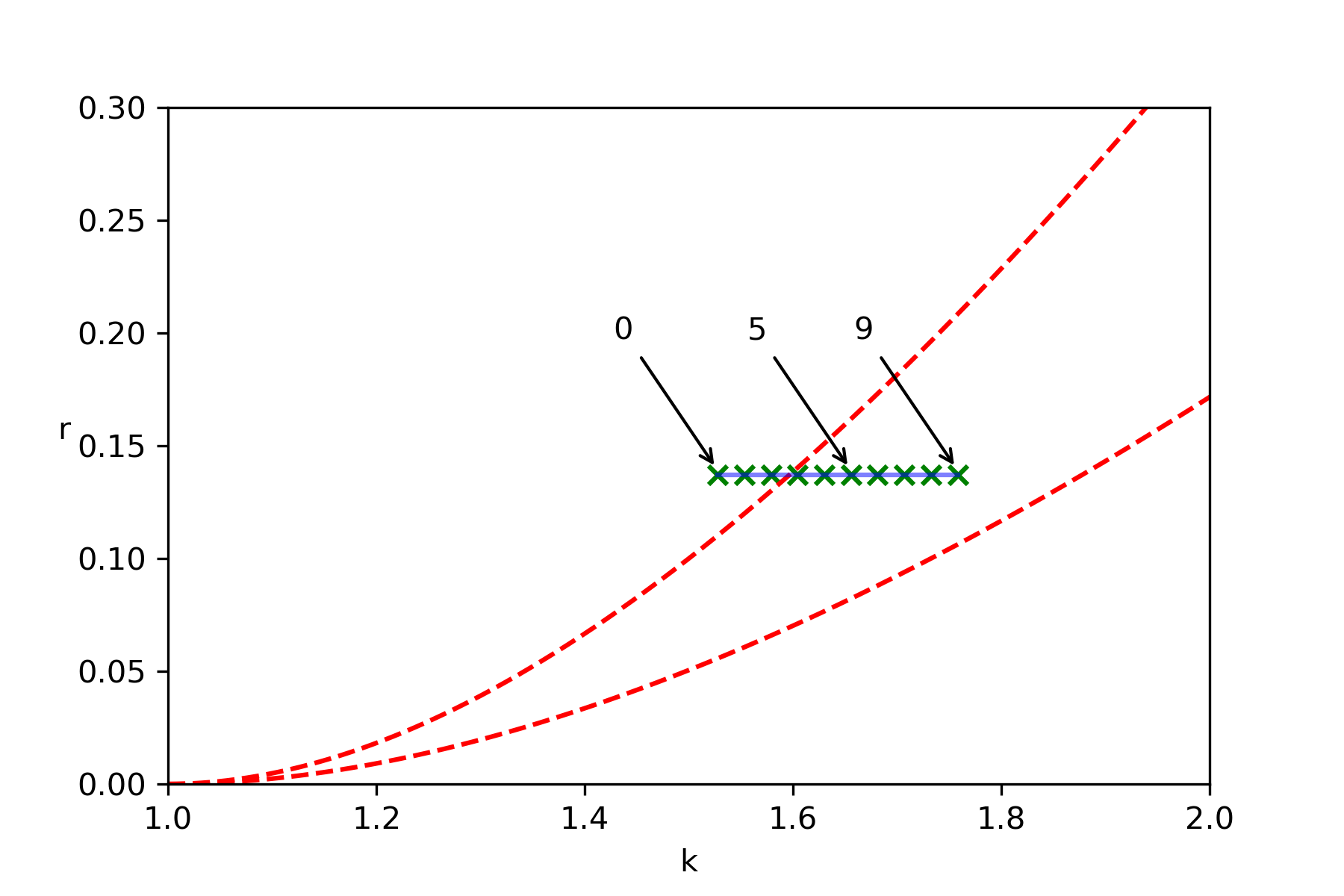}
		\captionof{figure}{Segment travelled in the $(\rZN,\kZN)$ stability diagram}
		\label{fig:optim:ligneRK}
	\end{minipage}
	~
	\begin{minipage}[t]{0.45\textwidth}
		\centering
		\adjincludegraphics[width=\textwidth, trim={{0\width} {0\height} {0\width} {0.1\height}}, clip]{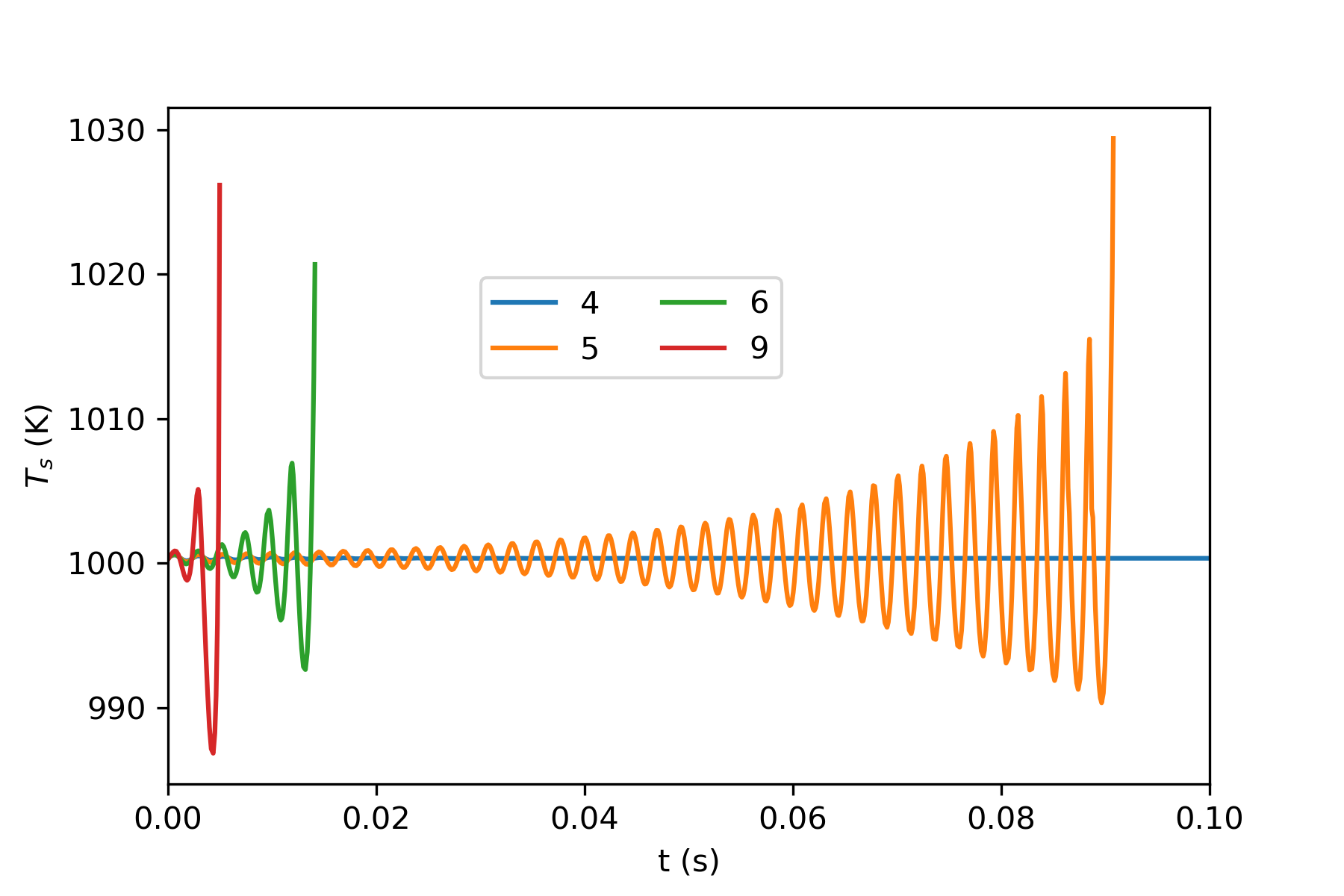}
		\captionof{figure}{Unsteady simulations}
		\label{fig:optim:main_transients}
	\end{minipage}
\end{center}
\end{figure}

We use the previous optimisation problem to generate configurations which have their sensitivity coefficients $\rZN$ and $\kZN$ distributed regularly on a segment defined
as $\rZN=0.137$ (baseline value) and $\kZN \in [1.5, 1.75]$ (see Figure \ref{fig:optim:ligneRK}),
thus crossing the ZN stability limit.
The optimisation is constrained to preserve physically sound characteristics (surface temperature at 1000 K, regression speed of 1 cm/s at 50 atm, 3540 K final flame temperature).

For each point, we numerically assess the stability of the corresponding steady-state combustion. 
A slight constant pressure perturbation (typically 0.1\% of the prescribed pressure) is applied and the one-dimensional tool is run with ESDIRK-54A and a relative error tolerance $rtol=10^{-6}$.
The stability of the combustion can then be assessed numerically by analysing whether the perturbation is damped out or not.

The unsteady simulations for a few points are shown in Figure \ref{fig:optim:main_transients}.
We see that instability appears between points 4 and 5, slightly further to the right than predicted analytically, as already discussed.
Refining the search between these points allows us to find a configuration that exhibits a limit cycle.
The corresponding model parameters are the same as in Section \ref{section:simplemodel:parameters}, except for the following changes:
$T_{ap}=14668$ K, $c_p=692.8$ J/kg/K, $c_c=1253$ J/kg/K, $T_\initial=182.4$K,
$\lambda_c=0.65$ W/m/K,
$\lambda=0.362$ W/m/K, $\molarmass=57.9$ g/mol,
$\Delta h_f^0(G_1)=-2.28 \times 10^5$ J/kg,
$\Delta h_f^0(G_2)=-2.22 \times 10^6$ J/kg,
$A=340.4$ s$^{-1}$.

\section{Detailed formulations of the mass flow rate constraint}
\label{section:appendix:drhodt_ordre_eleve}
Here we present the detailed discrete equations that can be used to account for the mass flow rate constraint arising from the continuity equation, following the discussion from Section \ref{section:RK:contrainte_debit}.

\subsection{Reformulation as an instantaneous constraint}
\label{appendix:contrainte:instant}
The term $-\partialdershort{\rho}{t}$ in the continuity equation $\partialdershort{\massflux}{x} = -\partialdershort{\rho}{t}$ can be considered as a source term which is a function of $\massflux$.
It can be obtained by differentiating the logarithm of the equation of state \cref{eq:base:idealgaslaw} with respect to time,
as classically done in the combustion community \cite{knio99,najm05,massot2014,MOTHEAU2016430}:

\begin{equation}
\partialdershort{\massflux}{x} = -\rho
\left( \dfrac{\partialdershort{P}{t}}{P} - \dfrac{\partialdershort{T}{t}}{T}
- \dfrac{\sum\limits_{\indexspecies=1}^{n_e} \dfrac{\partialdershort{Y_\indexspecies}{t}}{\molarmass_\indexspecies}}{\sum\limits_{\indexspecies=1}^{n_e}
\dfrac{Y_\indexspecies}{\molarmass_\indexspecies}} \right)
\label{eq:appendix:log_der_state}
\end{equation}
In the right hand-side, $\partialdershort{P}{t}$ is an input (prescribed pressure variation, or evolution based on a combustion chamber model).
The gas density $\rho$ is directly given by the equation of state \cref{eq:base:idealgaslaw}.
The other terms can be constructed based on our original gas phase system. For instance
$\partialdershort{Y_\indexspecies}{t} = \left(\partialdershort{\rho Y_\indexspecies}{t} - Y_\indexspecies \partialdershort{\rho}{t}\right)/{\rho}$,
where $\partialdershort{\rho Y_\indexspecies}{t}$ is given by Equation \cref{eq:base:species}, and $\partialdershort{\rho}{t}$ is replaced by $-\partialdershort{\massflux}{x}$ as per the continuity equation \cref{eq:base:continuity}.
The term $\partialdershort{T}{t}$ is slightly more involved.
\noindent We can write $\partialdershort{\hcomplet}{t} = 
\left( \partialdershort{\rho \hcomplet}{t} - \hcomplet \partialdershort{\rho}{t} \right)/\rho$,
where $\partialdershort{\rho h}{t}$ is given by Equation \cref{eq:base:enthalpy}.
Let us express $\partialdershort{T}{t}$ by differentiating with respect to time the definition of the enthalpy given in Section \ref{section:baseEquations}:
\begin{align}
	&\begin{aligned}
		& &\partialdershort{\hcomplet}{t} &=&
		&\sum\limits_{\indexspecies=1}^{n_e} \partialdershort{Y_\indexspecies}{t}
		\left( \Delta h_{f,\indexspecies}^0 + \int_{T_0}^{T} c_{p,\indexspecies}(a) da \right)
		+
		\sum\limits_{\indexspecies=1}^{n_e} Y_\indexspecies
		c_{p,\indexspecies}(T) \partialdershort{T}{t}
	\end{aligned}\\
	&\begin{aligned}
		&\Rightarrow
		&\partialdershort{T}{t}  &=&
		&\dfrac{
			\partialdershort{\hcomplet}{t} - \sum\limits_{\indexspecies=1}^{n_e} \partialdershort{Y_\indexspecies}{t}
			\left( \Delta h_{f,\indexspecies}^0 + \int_{T_0}^{T} c_{p,\indexspecies}(a) da \right)
		}
		{
			\sum\limits_{\indexspecies=1}^{n_e} Y_\indexspecies
			c_{p,\indexspecies}(T)
		}
	\end{aligned} \label{eq:instantconstraint:dtT}
\end{align}

Overall, in the $\indexcell$-th cell, the constraint \eqref{eq:appendix:log_der_state} can be discretised as:
\begin{equation}
 \dfrac{m_{\indexcell+1}-m_{\indexcell}}{\Delta x_\indexcell} =
 -\rho_\indexcell
\left( \dfrac{\partialdershort{P}{t}}{P} - \dfrac{\partialdershort{T_\indexcell}{t}}{T_\indexcell}
- \dfrac{\sum\limits_{\indexspecies=1}^{n_e} \dfrac{\partialdershort{Y_{\indexcell,\indexspecies}}{t}}{\molarmass_\indexspecies}}{\sum\limits_{\indexspecies=1}^{n_e}
\dfrac{Y_{\indexcell,\indexspecies}}{\molarmass_\indexspecies}} \right)
\label{eq:instant_constraint:discrete}
\end{equation}
\noindent where the various terms are evaluated following the previous formulas and the conservation equations \cref{eq:VF:continuity,eq:VF:species,eq:VF:enthalpy}. Equation \cref{eq:instant_constraint:discrete} is then used used in place of its alternative formulation \cref{eq:VF:continuity} in our DAE system.

Solving Equation \cref{eq:instant_constraint:discrete} for the discrete values $\massflux_i$ can be done by transforming it into a problem of the form $0=g(\massflux)$ and applying a Newton method. This system is linear in $\massflux$, thus a single Newton step would suffice. However we can see that the Jacobian $\partialdershort{g}{\massflux}$ involves the solution profiles $T$, $Y_\indexspecies$, therefore it will need to be updated often as the solution evolves in time and this may end up being costly. This is in agreement with our computational findings for the case where all fields are solved in a fully-coupled manner: the instantaneous formulation requires more Jacobian updates than the constraint formulation based on the Runge-Kutta quadrature presented hereafter.

\subsection{Natural formulation with the Runge-Kutta scheme}
\label{appendix:contrainte:RK}
An alternative and original approach of this paper is to apply the Runge-Kutta scheme to the semi-discrete continuity equation \cref{eq:VF:continuity} directly to obtain a
quadrature formula on $\int_t \partialdershort{\rho}{t} dt$. For the sake of readability, we momentarily change our notations to better distinguish the spatial and time representations. Let $q_\indexstep^{\indexcell}$ be the value of a variable $q$ in the $\indexcell$-th cell at time step $\indexstep$ (or at the $\indexcell$-th face for $\massflux$), and 
$q^{\indexcell}_{\indexstep,\indexstage}$ the same variable at the $\indexstage$-th stage of time step $\indexstep$.
For the $\indexstage$-th stage of any implicit Runge-Kutta  method, we obtain:
\begin{equation}
\label{eq:approx_drhodt_ordre_eleve:RKbase}
\rho^{\indexcell}_{\indexstep,\indexstage} - \rho_\indexstep^\indexcell = \int_{t_n}^{t_{ni}} (\derivshort{\rho^\indexcell}{t}) dt
\approx \Delta t \sum\limits_{\indexstagebis=1}^{\totalstep} a_{\indexstage \indexstagebis} \left(d_t \rho^\indexcell \right)_{\indexstep,\indexstagebis}
\end{equation}

\noindent where $\left(d_t \rho^\indexcell \right)_{\indexstep,\indexstagebis}$ is the time derivative of $\rho^\indexcell$ at time $t_\indexstep + c_\indexstagebis \Delta t$
(i.e. at the $\indexstagebis$-th stage). Based on the semi-discrete mass conservation equation \cref{eq:VF:continuity}, it is equal to the numerical approximation of the mass flow rate spatial gradient at this stage. Equation \cref{eq:approx_drhodt_ordre_eleve:RKbase} can then be interpreted as a constraint on $\massflux$:
\begin{equation}
  \label{eq:approx_drhodt_ordre_eleve:contrainte_RK}
-\dfrac{\massflux_{\indexstep \indexstage}^{\indexcell+1} - \massflux_{\indexstep \indexstage}^{\indexcell}}{\Delta x^\indexcell}
 = \dfrac{ \rho^{\indexcell}_{\indexstep \indexstage} - \rho_\indexstep^\indexcell }{a_{\indexstage \indexstage} \Delta t}
 + \sum\limits_{\indexstagebis=1, \indexstagebis \neq \indexstage}^{\totalstep}
      \dfrac{\massflux_{\indexstep \indexstagebis}^{\indexcell+1} - \massflux_{\indexstep \indexstagebis}^{\indexcell}}{\Delta x^\indexcell}     
\end{equation}
Comparing this equation to the semi-discrete continuity equation
$-({\massflux_{\indexstep \indexstage}^{\indexcell+1} - \massflux_{\indexstep \indexstage}^{\indexcell}})/{\Delta x^\indexcell}
= \left(d_t \rho^\indexcell \right)_{\indexstep \indexstage}$, we see that the right hand-side is the approximation of the source term $\left(d_t \rho^\indexcell \right)_{\indexstep \indexstage}$ that is naturally
constructed by the Runge-Kutta scheme and which can be entirely expressed in terms of the mass flow rates at various stages.

\section{Verification of the order of convergence}
\label{appendix:convergence}
\begin{figure}[htpb]
\begin{center}
	\begin{minipage}[t]{0.45\textwidth}
		\centering
		\adjincludegraphics[width=\textwidth, trim={{0\width} {0\height} {0\width} {0.1\height}}, clip]{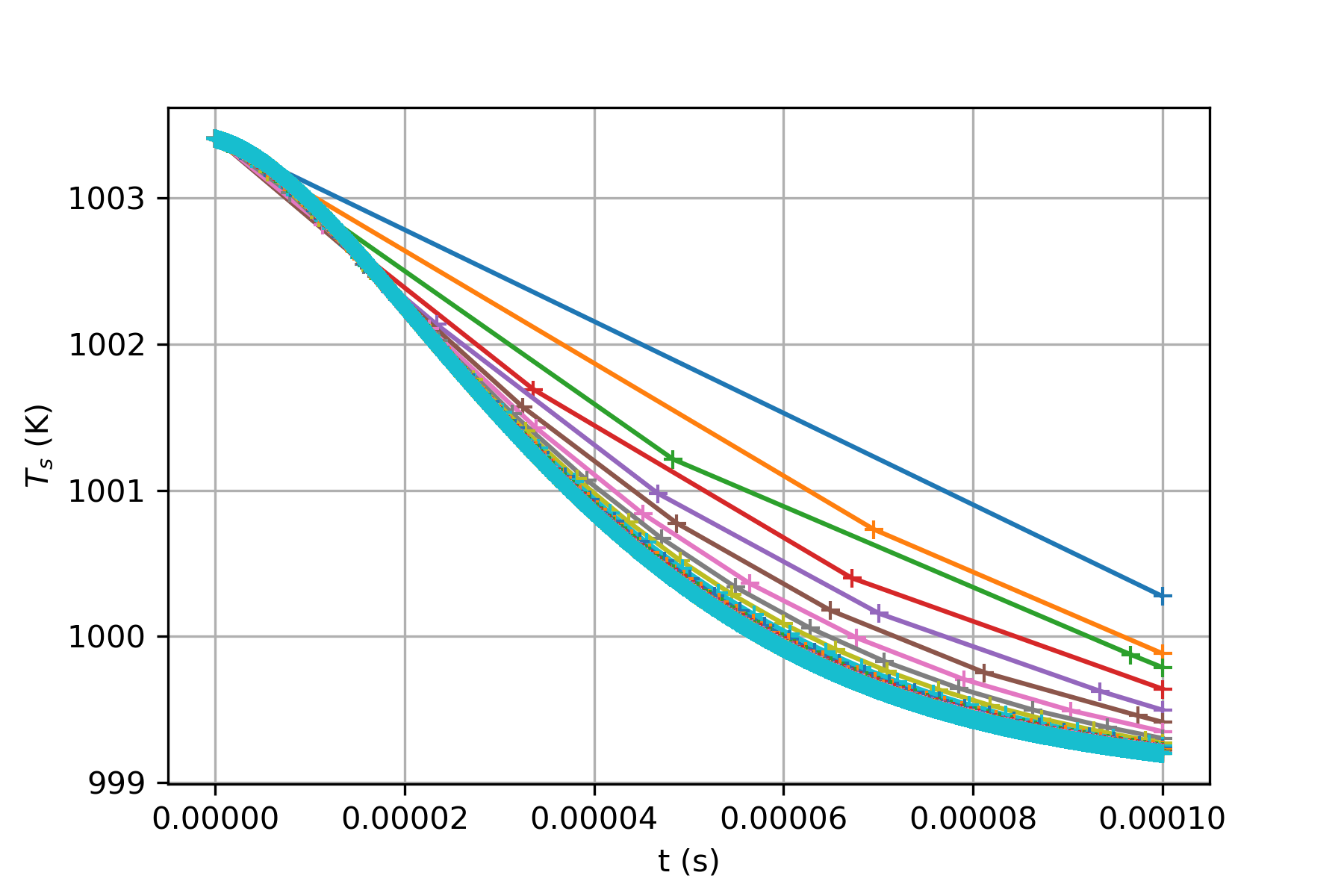}
		\captionof{figure}{Surface temperature histories obtained with IE when gradually lowering the time step}
		\label{fig:appendix:ordre:courbe_Ts_evolution}
	\end{minipage}
	~
	\begin{minipage}[t]{0.45\textwidth}
		\centering
		\adjincludegraphics[width=\textwidth, trim={{0\width} {0\height} {0\width} {0.1\height}}, clip]{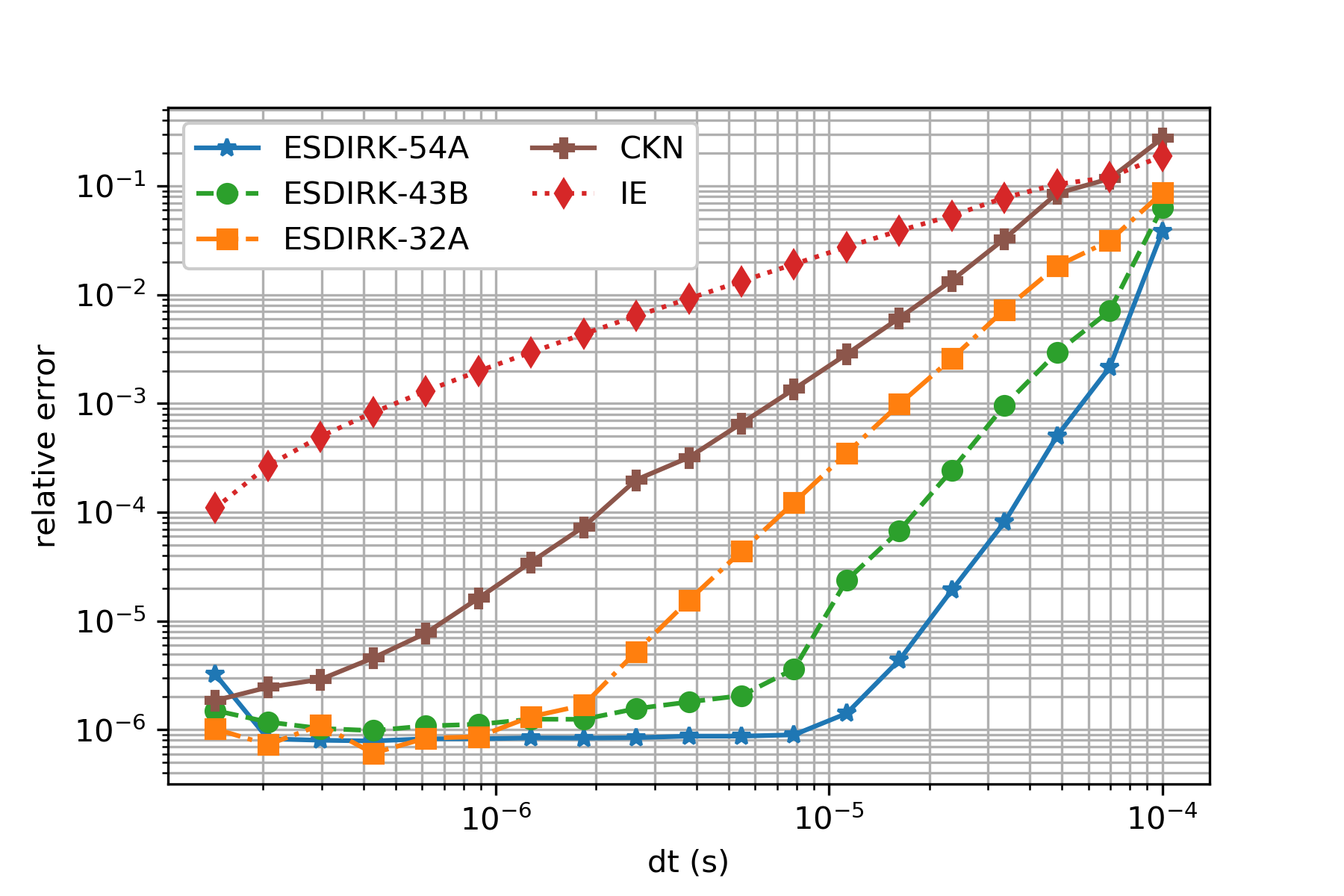}
		\captionof{figure}{Global error $\epsilon_{\rho u}$ on the mass flow rate field}
		\label{fig:appendix:ordre:convergence}
	\end{minipage}
\end{center}
\end{figure}

To verify the orders of convergence in time, a simple test case is set up, using the simplified model presented in Section \ref{section:verification}.
Starting from the steady-state solution at $P=5.5\times  10^{6}$ Pa, a pressure step to $P=5\times 10^{6}$ Pa is applied and the resulting transient is simulated for $10^{-4}$ s with a fixed time step.
The curves of $T_\surface$ obtained for various time step values are plotted in Figure \ref{fig:appendix:ordre:courbe_Ts_evolution}. The cyan curve represents the most refined solution.
A space-averaged relative error, plotted in Figure \ref{fig:appendix:ordre:convergence}, is computed on the final mass flow rate field:
\begin{equation}
	\epsilon_{\rho u} = \sqrt{ \sum\limits_{1}^{N+1} \frac{1}{N} \left( \dfrac{\massflux_\indexcell(t_f; \Delta t) - \massflux_\indexcell(t_f; dt_{ref})}{\massflux_\indexcell(t_f; \Delta t_{ref})} \right)^2}
\end{equation}
\noindent with $N$ the global number of cells and $ref$ denoting the reference simulation.
The theoretical orders of convergence are attained as long as the error is not limited by the precision of the Newton algorithm. Similar convergence rates have been observed for the other variables (point-wise or space-averaged errors).

\section{Jacobian for the Newton increment}
\label{appendix:jacobienne}
Applying an ESDIRK scheme (stage Equations \eqref{eq:RK:souspas_yni} and \eqref{eq:RK:souspas_zni}) to the system \eqref{eq:semiexplicit:conservative_ode}-\eqref{eq:semiexplicit:conservative_contraintes}, the $\indexstage$-th stage reads:
\begin{alignat}{2}
	Q(\Xdiffrksub_{\indexstep \indexstage}) &~=~& \Delta t a_{\indexstage\indexstage} &f(\Xdiffrksub_{\indexstep \indexstage}, \Xdaerksub_{\indexstep \indexstage})
	+\overbrace{ Q(\Xdiff^{n})
	+ \sum\limits_{\indexstagebis=0}^{\indexstage-1} \Delta t a_{\indexstage\indexstagebis} f(\Xdiffrksub_{\indexstep \indexstagebis}, \Xdaerksub_{\indexstep \indexstagebis})}^{b} \\
	0 &~=~& &g(\Xdiffrksub_{\indexstep \indexstage}, \Xdaerksub_{\indexstep \indexstage})
\end{alignat}

The residual vector is therefore
$F(\Xdiffrksub_{\indexstep \indexstage},\Xdaerksub_{\indexstep \indexstage})
=
\left(
Q(\Xdiffrksub_{\indexstep \indexstage}) - \Delta t a_{\indexstage\indexstage} f(\Xdiffrksub_{\indexstep \indexstage}, \Xdaerksub_{\indexstep \indexstage})
- b,~ g(\Xdiffrksub_{\indexstep \indexstage}, \Xdaerksub_{\indexstep \indexstage}) \right)^t$.
Therefore, the Jacobian \mbox{$J=\dfrac{\partial F}{\partial (\Xdiffrksub_{\indexstep \indexstage},~ \Xdaerksub_{\indexstep \indexstage})}$} needed for the Newton step is:
\begin{equation}
	J = \left(\begin{array}{cc}
		\partial_\Xdiffrksub Q -\Delta t a_{\indexstage\indexstage} \partial_\Xdiffrksub f & -\Delta t a_{\indexstage\indexstage} \partial_\Xdaerksub f\\
		\partial_\Xdiffrksub g & \partial_\Xdaerksub g
	\end{array}\right)
\end{equation}

\end{document}